\documentclass[preprint,review,12pt]{elsarticle}
\usepackage{geometry}   
\usepackage[fleqn]{amsmath}
\usepackage{xcolor}

\usepackage{natbib}
\biboptions{sort&compress}
\usepackage{graphicx}
\usepackage{epsfig}
\usepackage{epstopdf}
\usepackage{amssymb}
\usepackage{amsfonts}
\usepackage{amssymb}
\usepackage{booktabs}
\usepackage{color,soul}
\usepackage{float}
\usepackage[percent]{overpic}
\newcommand{\overbar}[1]{\mkern 1.5mu\overline{\mkern-1.5mu#1\mkern-1.5mu}\mkern 1.5mu}
\usepackage{bbm}
\usepackage{dsfont}
\usepackage{mathrsfs}
\usepackage[bbgreekl]{mathbbol}

\newcommand{\fre}{\mathfrak{e}}

\newcommand{\mbba}{\mathbbm{a}}
\newcommand{\mbbb}{\mathbbm{b}}

\newcommand{\mbbd}{\mathbbm{d}}
\newcommand{\mbbe}{\mathbbm{e}}

\newcommand{\mbbg}{\mathbbm{g}}

\newcommand{\mbbn}{\mathbbm{n}}

\newcommand{\mbbp}{\mathbbm{p}}

\newcommand{\mbbs}{\mathbbm{s}}
\newcommand{\mbbt}{\mathbbm{t}}
\newcommand{\mbbu}{\mathbbm{u}}
\newcommand{\mbbv}{\mathbbm{v}}
\newcommand{\mbbw}{\mathbbm{w}}
\newcommand{\mbbx}{\mathbbm{x}}
\newcommand{\mbby}{\mathbbm{y}}

\newcommand{\mbbA}{\mathbbm{A}}
\newcommand{\mbbB}{\mathbbm{B}}
\newcommand{\mbbC}{\mathbbm{C}}
\newcommand{\mbbD}{\mathbbm{D}}
\newcommand{\mbbE}{\mathbbm{E}}
\newcommand{\mbbF}{\mathbbm{F}}
\newcommand{\mbbG}{\mathbbm{G}}
\newcommand{\mbbH}{\mathbbm{H}}

\newcommand{\mbbK}{\mathbbm{K}}

\newcommand{\mbbM}{\mathbbm{M}}
\newcommand{\mbbN}{\mathbbm{N}}

\newcommand{\mbbP}{\mathbbm{P}}

\newcommand{\mbbR}{\mathbbm{R}}

\newcommand{\mbbU}{\mathbbm{U}}
\newcommand{\mbbV}{\mathbbm{V}}
\newcommand{\mbbW}{\mathbbm{W}}
\newcommand{\mbbX}{\mathbbm{X}}
\newcommand{\mbbY}{\mathbbm{Y}}

\newcommand{\calbfB}{\boldsymbol{\mathcal B}}
\newcommand{\calbfC}{\boldsymbol{\mathcal C}}

\newcommand{\calbfE}{\boldsymbol{\mathcal E}}

\newcommand{\bfb}{\bold{b}}

\newcommand{\bfm}{\bold{m}}

\newcommand{\bfB}{\bold{B}}

\newcommand{\bfF}{\bold{F}}

\newcommand{\bfH}{\bold{H}}
\newcommand{\bfI}{\bold{I}}
\newcommand{\bfJ}{\bold{J}}

\newcommand{\bfM}{\bold{M}}

\newcommand{\bfP}{\bold{P}}
\newcommand{\bfQ}{\bold{Q}}
\newcommand{\bfR}{\bold{R}}

\newcommand{\bfU}{\bold{U}}
\newcommand{\bfV}{\bold{V}}

\newcommand{\bfY}{\bold{Y}}

\newcommand{\bfalpha}{\boldsymbol{\alpha}}       

\newcommand{\bfomega}{{\boldsymbol{\omega}}}
\newcommand{\bfpsi}{\boldsymbol{\psi}}

\newcommand{\bfgamma}{\boldsymbol{\gamma}}
\newcommand{\bfsigma}{\boldsymbol{\sigma}}

\newcommand{\bfzeta}{\boldsymbol{\zeta}}

\newcommand{\bftheta}{\boldsymbol{\theta}}

\newcommand{\bfLambda}{\boldsymbol{\Lambda}}
\newcommand{\bfPhi}{\boldsymbol{\Phi}}

\newcommand{\bfUpsilon}{\boldsymbol{\Upsilon}}
\newcommand{\bfGamma}{\boldsymbol{\Gamma}}

\newcommand{\bfTheta}{\boldsymbol{\Theta}}

\newcommand{\calA}{\mathcal{A}}%
\newcommand{\calB}{\mathcal{B}}%
\newcommand{\calC}{\mathcal{C}}%

\newcommand{\calI}{\mathcal{I}}%
\newcommand{\calJ}{\mathcal{J}}%
\newcommand{\calK}{\mathcal{K}}%
\newcommand{\calL}{\mathcal{L}}%
\newcommand{\calN}{\mathcal{N}}%
\newcommand{\calP}{\mathcal{P}}%
\newcommand{\calS}{\mathcal{S}}%
\newcommand{\calU}{\mathcal{U}}%
\newcommand{\calV}{\mathcal{V}}%
\newcommand{\calW}{\mathcal{W}}%
\newfont{\sforf}{cmssq8 at 10pt}
\newfont{\tenss}{lcmss8 at 10pt}
\newfont{\ltenss}{lcmss8 at 8 pt}

\newfont{\forf}{cmssbx10}
\newfont{\tenbss}{cmssdc10}
\newfont{\ttenbss}{lcmssb8 at 10 pt}
\newfont{\ltenbss}{lcmssb8}

\newfont{\svnbss}{cmssbx10 scaled 900}

\newfont{\smsam}{msam5 scaled 1175}

\newfont{\ssmsam}{msam5 scaled 750}

\newcommand{\half}{{\textstyle{\frac{1}{2}}}}

\newcommand{\Lgrad}{\hbox{Grad}\mskip2mu}

\newcommand{\Egrad}{\hbox{grad}\mskip2mu}

\newcommand{\tr}{\hbox{tr}\mskip2mu}
\newcommand{\trans}{\scriptscriptstyle\mskip-1mu\top\mskip-2mu}

\newcommand{\lj}{\mbox{$[\kern-0.1478125em[$}}
\newcommand{\rj}{\mbox{$]\kern-0.1478125em]$}}

\newcommand{\la}{\mbox{$\langle\kern-0.2325em\langle$}}
\newcommand{\ra}{\mbox{$\rangle\kern-0.2325em\rangle$}}
\newcommand{\Blj}{\mbox{$\bigg[\kern-0.275em\bigg[$}}
\newcommand{\Brj}{\mbox{$\bigg]\kern-0.275em\bigg]$}}

\newcommand{\tendot}{\mskip1.75mu\raisebox{0.25ex}:\mskip2mu}

\newfont{\tenrm}{cmr10}%

\newfont{\sscmsy}{cmsy5 scaled 750}

\journal{JOURNAL}
\begin{document}
\begin{frontmatter}
\title{A micropolar shell model for hard-magnetic soft materials}
{\author[1]{Farzam~Dadgar-Rad}
\ead{dadgar@guilan.ac.ir}
\address[1]{Faculty of Mechanical Engineering, University of Guilan, Rasht, Iran}}
{\author[2]{Mokarram~Hossain\corref{cor2}}
\ead{mokarram.hossain@swansea.ac.uk}
\cortext[cor2]{Corresponding author.}
\address[2]{Zienkiewicz Centre for Computational Engineering, College of Engineering, Swansea University, SA1 8EN, UK}}
%
%
%
\vspace{0pt}
\begin{abstract}
Hard-magnetic soft materials (HMSMs) are particulate composites that 
particles with high coercivity are dispersed in a soft matrix. 
Since applying the magnetic loading induces a body couple in HMSMs, the resulting Cauchy stress is predicted to be asymmetric.
Therefore, the micropolar continuum theory can be employed to capture the deformation of these materials.
On the other hand, the geometries and structures made of HMSMs often possess small thickness compared to the overall dimensions of the body.
Accordingly, in the present contribution, a 10-parameter micropolar shell formulation to model the finite elastic deformation of thin hard-magnetic soft structures under magnetic stimuli is developed.
The proposed shell formulation allows for using three-dimensional constitutive laws without any need for modification to apply the plane stress assumption in thin structures.
A nonlinear finite element formulation is also presented for the numerical solution of the governing equations.
To alleviate the locking phenomenon, the enhanced assumed strain method is employed. 
Several examples are presented that demonstrate the performance and effectiveness of the proposed formulation.  
%
%
\end{abstract}
\end{frontmatter}

\vspace{-6pt}
\section{Introduction}
\label{intro}
\vspace{-2pt}
%
Magneto-active soft materials consist of magnetic particles dispersed into a soft elastomeric matrix and undergo large deformations under magnetic loading. 
This class of materials has been used in vibration absorbers,  sensors, actuators, soft robots, flexible electronics, and isolators   (see, e.g.,  
\cite{Ren2019, Wu2020, Bastola2021, Bastola2020, Lucarini2022,  Yarali2022} and references therein). 
Therefore, developing reliable theoretical models plays an essential role in the optimum and cost-effective design of the aforementioned devices and instruments.

Based on the type of the embedded particles, magneto-active soft materials are divided into two sub-classes, namely \textit{soft-magnetic soft materials} (SMSMs) and \textit{hard-magnetic soft materials} (HMSMs).
The former contains particles with low coercivity, such as iron or iron oxides, and their magnetization vector varies under external magnetic loading.
This sub-class has been the subject of a huge amount of research work in this century (e.g., \cite{Saxena2013, Miehe2016, Mehnert2017, Mukherjee2020, Bustamante2021, Hu2022, akbari2020}).
The latter sub-class is composed of particles of high coercivity, such as CoFe$_2$O$_4$ or NdFeB, so that their magnetization vector, or equivalently, their remnant magnetic flux, remains unchanged for a wide range of the applied external magnetic flux (e.g., \cite{Schumann2020, Lee2020}).
One of the main characteristics of HMSMs is that external magnetic induction of relatively small magnitude causes rapid finite deformations in these materials (e.g., \cite{Lum2016, Wu2019}).
Moreover, the 3D printing technologies have enabled the researchers to program the ferromagnetic domains in complex structures, which leads to the desired deformations \cite{Kim2018, Alapan2020, Kuang2021, Wang2021, Wu2021}.

Theoretical modeling of HMSMs has been the subject of a plethora of research articles in recent years (e.g., \cite{Kalina2017, Zhao2019, Garcia2019, Mukherjee2021,  rambausek2022, Zhang2020, Garcia2021a, Garcia2021b, Ye2021, DH2022IJSS}). 
In particular, Zhao et al.~\cite{Zhao2019} developed a continuum formulation with asymmetric Cauchy stress so that the external magnetic induction directly contributes to the expression of the stress tensor.
Their theory has been the foundation for the analysis of hard-magnetic soft beams (HMSBs) in Yan et al.~\cite{Yan2021},
Wang et al.~\cite{Wang2020}, 
Rajan and Arockiarajan~\cite{Rajan2021},
and Chen et al.~\cite{Chen2020a, Chen2021}
among others.
The same formulation has been employed to model the deformation of magneto-active shells by Yan et al.~\cite{Reis2021}.
Dadgar-Rad and Hossain~\cite{DH2022EML} enhanced the formulation of Zhao et al.~\cite{Zhao2019} to account for viscoelastic effects and analyzed the time-dependent dissipative response of HMSBs.
Some researchers have developed micromechanical and lattice models for HMSMs  (e.g., Refs.~\cite{Zhang2020,Garcia2021a, Garcia2021b,Ye2021}).
From a different point of view, Dadgar-Rad and Hossain~\cite{DH2022IJSS} focused on the well-known phenomenon that due to the presence of the remnant magnetic induction in HMSMs, applying an external magnetic loading produces a body couple that plays the role of the main driving load to generate mechanical deformation in the body (e.g., \cite{Dorfmann2014}).
Moreover, the Cauchy stress losses its symmetry, as had been previously pointed out by Zhao et al.~\cite{Zhao2019}. 
However, instead of following the methodology advocated in \cite{Zhao2019}, the authors developed a formulation based on micropolar continuum theory to predict the deformation of 3D hard-magnetic soft bodies.
Two significant differences between the results of the formulation of Zhao et al.~\cite{Zhao2019} and those based on the micropolar-enhanced formulation have been expressed by Dadgar-Rad and Hossain~\cite{DH2022IJSS}.

Eringen and his coworkers established the theoretical foundations of micropolar theory \cite{Kafadar1971, Eringen1976, Eringen1999}. 
In this theory, a microstructure at each material point is considered.
The microstructure can have arbitrary rigid rotations independent of the traditional motion field considered in the classical continuum.
Formulations of the micropolar theory to model localized elastic-plastic deformations (e.g., \cite{Borst1993,  Steinmann1994, Tsakmakis2005, Grammenoudis2007a, Grammenoudis2007b, Bauer2012a, Borst2022}) and size-dependent elastic deformations (e.g., \cite{Pietraszkiewicz2009,  Ramezani2008, Ramezani2009,  Bauer2010, Bauer2012b, Erdelj2020}) have been developed.
Some formulations to model micropolar shells have been also proposed (e.g., \cite{Sargsyan2020, Eremeyev2005, Eremeyev2017}).
Moreover, the theory has been used in the modeling of lattice structures,
crystal plasticity, phononic crystals, chiral auxetic lattices, phase-field fracture mechanics, and vertebral trabecular bone \cite{Yoder2018, Mayeur2011, Zapata2020, Spadoni2012, Goda2014, Suh2020}.


The current research is essentially the continuation of the previous work of the authors, namely Dadgar-Rad and Hossain \cite{DH2022IJSS}, which had been developed for three-dimensional bodies.
However, most bodies made of HMSMs are thin structures, and using three-dimensional elements is computationally expensive. 
Accordingly, the purpose of this research is to develop a micropolar-based shell model to predict the deformation of thin HMSMs.
To do so, the 7-parameter shell formulation of Sansour (e.g., \cite{Sansour98, SansKoll2000}) has been extended to a 10-parameter one that involves the micro-rotation of the microstructure.
On the other hand, the enhanced assumed strain method (EAS) is a widely-used strategy to eliminate locking in shell structures, e.g., \cite{SimoArmero92, Simo93, KW96, Glaser97}.
Therefore, this method is adopted here to circumvent locking effects in the present micropolar shell formulation.

The next sections of this paper are as follows: 
The basic kinematic and kinetic relations of the micropolar continuum theory are summarized in Section~\ref{micropolar}.
In Section~\ref{HMSMs}, the main characteristics of HMSMs are presented.
In Section~\ref{shell}, the kinematic equations describing a 10-parameter micropolar shell model are provided.
Section~\ref{PVW} presents the variational formulation, followed by a FE formulation in Section~\ref{FEM}.
Numerical examples are studied in Section~\ref{Examples}, and the paper concludes in Section~\ref{conc}.

\vspace{10pt}
\textbf{Notation:} 
In this work, Greek indices take $1$ and $2$.
All upper-case and lower-case Latin indices take $1$, $2$, and $3$.
Upper-case indices with calligraphic font, e.g., $\calK$ and $\calL$, take the values specified in the corresponding equations. 
The repeated Latin and Greek indices obey Einstein's summation convention.
If $\bfP$ and $\bfQ$ are two 2nd-order tensors, 
the tensorial products defined via the symbols $\otimes$, $\odot$, and $\boxtimes$ generate 4th-order tensors so that the corresponding components are given by 
$(\calbfC)_{ijkl}=(\bfP \otimes \bfQ)_{ijkl}=P_{ij}Q_{kl}$,  
$(\calbfB)_{ijkl}=(\bfP \odot \bfQ)_{ijkl} =P_{ik}Q_{jl}$, and 
$(\calbfC)_{ijkl}=(\bfP \boxtimes \bfQ)_{ijkl} =P_{il}Q_{kj}$, respectively.
%
%
%
For numerical simulations, the notation $\mbbU=\{U_{11}, U_{22}, U_{33}, U_{12}, U_{21} , U_{13} , U_{31} , U_{23}, U_{32} \}^{\trans}$ will be used as the $9 \times 1 $ vectorial representation of the arbitrary 2nd-order tensor $\bfU$.
%

\section{A brief review of the micropolar theory}
\label{micropolar}
\vspace{-5pt}
The purpose of this section is to introduce some concepts and relations of the micropolar theory. 
The interested reader may refer to the pioneering works developed in Refs.~\cite{Eringen1976, Eringen1999, Steinmann1994} for more details and discussions.

In this section, two coincident Cartesian coordinates  $\{X_I\}$ and $\{x_i\}$, with  $\{\mbbE_I\}$ and $\{\mbbe_i\}$ as the corresponding basis vectors, are considered.
%
%
%
%
The center of a macro-element in the reference configuration $\calB_0$ is denoted by $\mbbX$.
After deformation by the deformation mapping $\bfpsi$, the center of the macro-element in the current configuration $\calB$ at the time $t$ is denoted by $\mbbx$, so that $\mbbx=\bfpsi(\mbbX,t)$.
The deformation gradient $\bfF$ is given by
\begin{equation}\label{bfF}
\bfF=\Lgrad \bfpsi, \quad   F_{iJ}=\frac{\partial x_i}{\partial X_I},
\quad J=\det \bfF >0,
\end{equation}
which can be uniquely decomposed as $\bfF=\bfR \bfU = \bfV \bfR$. Here, $\bfR$ is the macro-rotation tensor, and $\bfU$ and $\bfV$ are the symmetric positive definite right and the left stretch tensors, respectively.
For later use, the variation of the deformation gradient is written as follows:
\begin{equation}\label{dbfF}
\begin{split}
\delta \bfF=  \delta \bfY  \bfF
\quad \text{with} \quad
(\delta \bfY)_{ij} =
(\Egrad \delta \hat{\mbbu})_{ij}
=\frac{ \partial }{\partial x_j}\delta \hat{u}_i.
\end{split}
\end{equation}
%
%
%
Moreover, $\hat{\mbbu} = \mbbx-\mbbX$ and $\delta \hat{\mbbu}= \delta \mbbx$ are the actual and virtual displacement fields, respectively.
%
%


As a basic assumption in the micropolar theory, there exists a microstructure inside each macro-element so that it experiences rigid micro-rotations independent of the macro-motion $\mbbx$.
Let $\bftheta=\theta_i \mbbe_i$ denote the micro-rotation pseudo-vector, and $\theta=(\theta_i \theta_i)^{1/2}$ be its magnitude.
The corresponding micro-rotation tensor, denoted by  $\tilde{\bfR}$, can be expressed via the Euler--Rodriguez formula, namely (e.g., \cite{Steinmann1994, Ramezani2009})
\begin{equation}\label{Rtilde}
\tilde{\bfR} (\bftheta)  =\frac{1}{\theta^2} 
\big[ \theta^2  \bfI 
+ \theta   \sin \theta  \hat{\bftheta}
+(1- \cos \theta ) \hat{\bftheta}^2 \big],
\end{equation}
where $\hat{\bftheta}=  - \calbfE  \bftheta$, or $\hat{\theta}_{ij}= -\epsilon_{ijk} \theta_k$, is the skew-symmetric tensor corresponding to $\bftheta$.
Moreover, $\epsilon_{ijk}$ are the components of the alternating symbol $\calbfE$.
By defining $\delta \bftheta$ as the virtual micro-rotation pseudo-vector, the variation of $\tilde{\bfR}$ may be expressed via the following relations \cite{Ramezani2009, DH2022IJSS}
\begin{equation}\label{dbfR}
\begin{split}
\delta \tilde{\bfR}= - \calbfE  \bfLambda \delta \bftheta \tilde{\bfR}
\quad \text{with} \quad
\bfLambda= \frac{1}{\theta^3} 
\big[ \theta^2 \sin \theta   \bfI
+\theta (1-\cos \theta)  \hat{\bftheta}
+(\theta - \sin \theta)  \bftheta \otimes \bftheta \big].
\end{split}
\end{equation}
The deformation gradient, in the micropolar theory,  is decomposed as 
$\bfF=\tilde{\bfR} \tilde{\bfU} = \tilde{\bfV} \tilde{\bfR}$, from which it follows that (e.g., \cite{Steinmann1994}):
\begin{equation}\label{UVtilde}
\begin{split}
\tilde{\bfU}=\tilde{\bfR}^{\trans} \bfF, \quad 
\tilde{U}_{IJ}=\tilde{R}_{pI} F_{pJ}, \quad 
\tilde{\bfV}=\bfF \tilde{\bfR}^{\trans}, \quad 
\tilde{V}_{ij}=F_{iQ} \tilde{R}_{jQ}.
\end{split}
\end{equation}
%
%
%

To take the gradient of the micro-rotation into account, the material wryness tensor $\bfGamma$ and its spacial counterpart $\bfgamma$, are defined by \cite{Eringen1976, Eringen1999, Steinmann1994, Pietraszkiewicz2009}
\begin{equation}\label{Gamma}
\bfGamma= -\frac{1}{2} \calbfE \tendot ( \tilde{\bfR}^{\trans} \Lgrad \tilde{\bfR}  ), 
\quad
\Gamma_{IJ}=\frac{1}{2} \epsilon_{IPK} \tilde{R}_{iK} \tilde{R}_{iP,J},
\quad
\bfgamma = \tilde{\bfR} \bfGamma \tilde{\bfR}^{\trans},
\quad
\gamma_{ij}=\tilde{R}_{iP} \Gamma_{PQ} \tilde{R}_{jQ}.
\end{equation}
The deformation measures $\tilde{\bfU}$ and $\bfGamma$ are the main kinematic tensors to develop a formulation in material framework (see also, Refs.~\cite{Eringen1976,Steinmann1994}). 
Combinations of Eqs.~\eqref{dbfF}, \eqref{dbfR}, \eqref{UVtilde}$_1$, and \eqref{Gamma}$_1$, furnishes the following relations for the virtual kinematic tensors $\delta \tilde{\bfU}$ and $\delta \bfGamma$  \cite{DH2022IJSS}:
\begin{equation}\label{delUG}
\delta \tilde{\bfU}= \tilde{\bfR}^{\trans} (\delta \bfY - \delta \hat{\bfomega}) \bfF, \quad
\delta \bfGamma =  \tilde{\bfR}^{\trans} \Egrad \delta \bfomega \bfF,
\end{equation}
where $\delta \hat{\bfomega}=- \calbfE   \delta \bfomega$ is the skew-symmetric tensor corresponding to $\delta {\bfomega}$.

Next, in the current configuration $\calB$, let $\text{d} \calA$ and $\mbbn$ be an infinitesimal area element and its corresponding outward unit normal vector, respectively.
In the micropolar theory, the traction $\mbbt^{(\mbbn)}$ and the couple vector $\mbbs^{(\mbbn)}$ (as the moment per unit area) act on  $\text{d} \calA$.
Let $\bfsigma$ and $\bfm$ be the asymmetric Cauchy stress and the asymmetric couple stress corresponding to $\mbbt^{(\mbbn)}$ and $\mbbs^{(\mbbn)}$, respectively.
Accordingly, the well-known Cauchy's stress principle is extended as follows (e.g., \cite{DH2022IJSS, Eringen1976}):
\begin{equation}\label{tractions1}
\begin{split}
\mbbt^{(\mbbn)} =\bfsigma \mbbn, \quad 
t^{(\mbbn)}_i=\sigma_{ij} n_j, \quad
\mbbs^{(\mbbn)} =\bfm \mbbn, \quad
s^{(\mbbn)}_i=m_{ij} n_j.
\end{split}
\end{equation}
For later use, the first Piola--Kirchoff stress $\bfP$, the material stress $\tilde{\bfP}$, the first Piola--Kirchoff couple stress $\bfM$, and the material couple stress $\tilde{\bfM}$ are defined by
\begin{equation}\label{tractions2}
\begin{split}
\{ \bfP,\bfM  \}=J \{ \bfsigma,\bfm  \} \bfF^{-\trans}, \quad
\{ \tilde{\bfP}, \tilde{\bfM}  \} =
\tilde{\bfR}^{\trans} \{ \bfP , \bfM  \}
= J \tilde{\bfR}^{\trans} \{ \bfsigma,\bfm  \}  \bfF^{-\trans}.
\end{split}
\end{equation}
%
%
%
%
%

\section{Basic relations of HMSMs}
\label{HMSMs}
\vspace{-5pt}
The main property of hard-magnetic soft materials is the existence of a remnant magnetic flux density, that remains almost unchanged under a wide range of the applied external magnetic flux $\mbbB^{\text{ext}}$ (e.g., \cite{Lum2016, Wu2019, Zhao2019}). 
Let $\tilde{\mbbB}^{\text{rem}}$ and ${\mbbB}^{\text{rem}}$ be the remnant magnetic flux in the reference and current configurations, respectively.
The relation between $\tilde{\mbbB}^{\text{rem}}$ and ${\mbbB}^{\text{rem}}$ is as follows \cite{Zhao2019}:
\begin{equation}\label{Br}
{\mbbB}^{\text{rem}} = J^{-1} \bfF \tilde{\mbbB}^{\text{rem}},
\quad  
B^{\text{rem}}_i=J^{-1} F_{iJ} \tilde{B}^{\text{rem}}_J.
\end{equation}
The action of $\mbbB^{\text{ext}}$ on ${\mbbB}^{\text{rem}}$ leads to a body couple (moment per unit volume) in HMSMs.
The relations for the body couple per unit current volume $\mbbp$, and the body couple per unit reference volume $\mbbp^{\ast}$ may be written as (e.g., \cite{Zhao2019, Dorfmann2014})
\begin{equation}\label{couple}
\mbbp^{\ast}= J \mbbp =  \frac{J}{\mu_0} {\mbbB}^{\text{rem}} \times \mbbB^{\text{ext}}
=\frac{1}{\mu_0} \big(\bfF \tilde{\mbbB}^{\text{rem}} \big) \times \mbbB^{\text{ext}}, 
\end{equation}
where the constant $\mu_0= 4 \pi \times 10^{-7} \frac{N}{A^2}$ is the free space magnetic permeability.
For HMSMs, the external magnetic flux density $\mbbB^{\text{ext}}$ is often assumed to remain constant in space (e.g., Refs.~\cite{Zhao2019, Wang2020, Chen2020a, Chen2021, Rajan2021}).
Using this point, it has been proven that the following Maxwell equations   are satisfied in HMSMs (e.g., \cite{Zhao2019, Dorfmann2014}):
\begin{equation}\label{maxwell}
\text{Div}\mbbB =B_{I,I}=0,\quad 
\text{Curl}\mbbH=\epsilon_{IJK} H_{J,K} \mbbE_I= {\bf0},   
\end{equation}
where $\mbbB$ is the referential magnetic flux density, 
$\mbbH$ is the referential magnetic field, 
{\color{black}$\mbbE_I$ is the basis vector of the referential coordinate system $\{ X_I \}$ defined in the previous section}, and "$\text{Curl}$" is the referential curl operator.

\vspace{0pt}
\section{Kinematics of a 10-parameter micropolar shell model}
\label{shell}
\begin{figure}
\centering
\includegraphics[width=15cm]{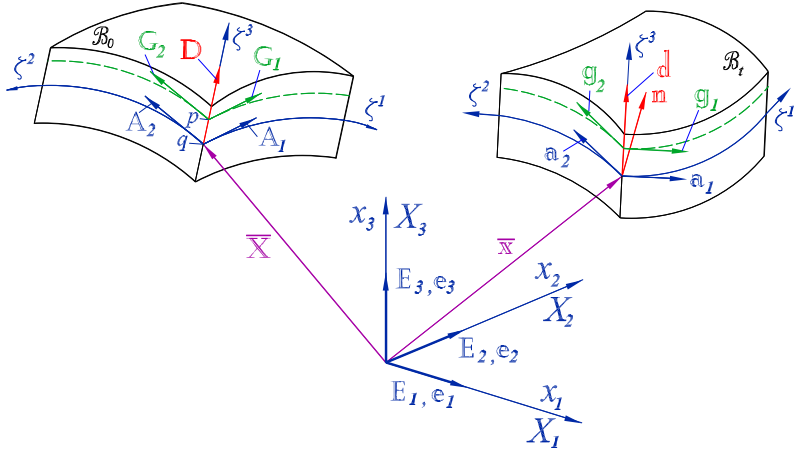} 
\vspace{-15pt}
\caption{Geometry of the shell}
\label{fig1geo} 
\end{figure}

\vspace{-5pt}
The geometry of a part $\calP$ of a shell in the reference and current configurations is displayed in Fig.~\ref{fig1geo}.
Let $\calS_0$ be the mid-surface of the shell in the reference configuration, which deforms into the surface $\calS$ in the current one. 
As shown in Fig.~\ref{fig1geo}, in addition to the two common-frame Cartesian coordinates $\{X_I\}$ and $\{x_i\}$, described in the previous section, the convective coordinate system  $\{\zeta^i\}$ at each material particle $q$ of the reference mid-surface $\calS_0$ is also constructed.
The coordinate lines $\zeta^i$ deform during the motion of the shell in space so that the coordinated lines $\zeta^1$  and $\zeta^2$  are tangent to both $\calS_0$ and  $\calS$. 
Moreover, the coordinated line 
$\zeta^3 \in [-\half h, \half h]$, with  $h$ as the initial thickness of the shell, is considered to be perpendicular to $\calS_0$ in the reference configuration. 
However, it does not remain perpendicular to $\calS$ in the current configuration, in general.
In the sequel, for the sake of simplicity, the coordinate $\zeta^3$ may be replaced by $z$. 
%

\vspace{0pt}

The position of the material particle $q$ on the mid-surface $\calS_0$ may be described by the vector $\overbar{\mbbX}(\zeta^1,\zeta^2)$. 
Let $\{\mbbA_{\alpha},\mbbA^{\alpha}, A_{\alpha\beta},A^{\alpha\beta}, \mbbD,\bfB\}$ be, respectively, the covariant and contravariant basis vectors, covariant and contravariant components of the metric tensor, outward unit normal vector, and the curvature tensor  on the undeformed mid-surface $\calS_0$. 
Then the following relations from the differential geometry of surfaces hold (e.g.,~\cite{I19}):
\begin{equation}\label{onM0}
\left.
\begin{split}
\mbbA_\alpha =\frac{ \partial \overbar{\mbbX}} {\partial \zeta^{\alpha}}, \quad
A_{\alpha\beta}=\mbbA_{\alpha} \cdot \mbbA_{\beta}, \quad
A^{\alpha\eta} A_{\eta\beta}= \delta^{\alpha}_{\beta}, \quad
\mbbA^{\alpha}=A^{\alpha\beta} \mbbA_{\beta}, \quad
\mbbA_{\alpha} \cdot \mbbA^{\beta} = \delta_{\alpha}^{\beta}\\
A=\det[A_{\alpha\beta}], \quad
\mbbD=\mbbA_3=\mbbA^3
=\frac  { \mbbA_1   \times   \mbbA_2}  {| \mbbA_1   \times   \mbbA_2 |} 
= \frac  { \mbbA_1   \times   \mbbA_2}  {\sqrt A}, \quad
\bfB= - \mbbD_{,\alpha} \otimes \mbbA^{\alpha}
\end{split}
\right \},
\end{equation}
where $\delta_{\alpha}^{\beta}$ is the two-dimensional Kronecker delta. 
For later use, the surface contravariant basis vectors may be written as
$\mbbA^{\alpha}=A^{\ast \alpha J} \mbbE_J$, where $A^{\ast \alpha J}$ are the Cartesian components of $\mbbA^{\alpha}$.
The position of the material particle $p$ located at the elevation $z$ with respect to  $\calS_0$ is described by
\begin{equation}\label{bfX}
\mbbX(\zeta^1,\zeta^2,z)
=\overbar{\mbbX}(\zeta^1,\zeta^2)
+z \mbbD(\zeta^1,\zeta^2),
\end{equation}
from which the covariant basis vectors $\mbbG_i$  are obtained to be
\begin{equation}\label{bfG}
\mbbG_{\alpha}=\frac {\partial \mbbX}{\partial \zeta^{\alpha}}
=\mbbA_{\alpha}+z \mbbD_{,\alpha}, \quad
\mbbG_{3}=\frac{\partial \mbbX}{\partial z}=\mbbD.
\end{equation}
Motivated by Eqs.~\eqref{onM0}$_8$ and \eqref{bfG}, the symmetric shifter tensor $\bfQ=\bfQ^{\trans}=\mbbG_i \otimes \mbbA^i= \bfI - z \bfB$, with ${\bfI}$  as the identity tensor, is defined.
The shifter tensor can be used to map the covariant and contravariant basis vectors from $z \neq 0$ to the mid-surface with $z=0$, and vice versa.
More precisely, the following relations hold:
\begin{equation}\label{ShifterS}
\mbbG_i=\bfQ \mbbA_i, \quad 
\mbbA_i = \bfQ^{-1} \mbbG_i, \quad
\mbbG^i=\bfQ^{-1} \mbbA^i, \quad 
\mbbA^i = \bfQ \mbbG^i, \quad
\end{equation}
where $\mbbG^i$ are the contravariant basis vectors at $\mbbX$, and use has been made of the symmetry property of $\bfQ$.
For later use, the three-dimensional material gradient operator $\text{Grad}_{\bfzeta}$, with respect to the convective coordinates $\{\zeta^i\}$ in the reference configuration, and the material surface gradient operator $\text{Grad}_{\calS_0}$ with respect to $\{\zeta^ {\alpha} \}$ are defined as follows:
\begin{equation}\label{Gradians}
\Lgrad_{\bfzeta} \{\bullet\}
= \frac{\partial \{\bullet\}} {\partial \zeta^i} \otimes \mbbG^i, \quad
\Lgrad_{\calS_0} \{\bullet\}
= \frac {\partial \{\bullet\}} {\partial \zeta^{\alpha}} \otimes \mbbA^{\alpha}.
\end{equation}
By assuming that a straight material fiber perpendicular to $\calS_0$  remains straight during deformation, the following macro deformation field is considered (e.g., \cite{Sansour98, SansKoll2000}):
\begin{equation}\label{bfpsi}
\mbbx=\bfpsi(\zeta^1,\zeta^2,z,t)
=\overbar{\mbbx}(\zeta^1,\zeta^2,t)
+z[1+z \phi(\zeta^1,\zeta^2,t)] \mbbd(\zeta^1,\zeta^2,t),
\end{equation}
where $\overbar{\mbbx}$ is the image of $\overbar{\mbbX}$ on  $\calS$, and $\mbbd$ is a director vector along the deformed $z$-axis. 
Moreover, the scalar field $\phi$  describes through the thickness stretching of the shell. 
Similar to the quantities defined on $\calS_0$ in Eq.~\eqref{onM0}, let $\{\mbba_{\alpha},\mbba^{\alpha}, a_{\alpha\beta},a^{\alpha\beta}, \mbbn,\bfb\}$ be the surface quantities defined  on $\calS$. It then follows that
\begin{equation}\label{onM}
\left.
\begin{split}
\mbba_\alpha =\frac{ \partial \overbar{\mbbx}} {\partial \zeta^{\alpha}}, \quad
a_{\alpha\beta}=\mbba_{\alpha} \cdot \mbba_{\beta}, \quad
a^{\alpha\eta} a_{\eta\beta}= \delta^{\alpha}_{\beta}, \quad
\mbba^{\alpha}=a^{\alpha\beta} \mbba_{\beta}, \quad
\mbba_{\alpha} \cdot \mbba^{\beta} = \delta_{\alpha}^{\beta}\\
a=\det[a_{\alpha\beta}], \quad
\mbbn=\frac  { \mbba_1   \times   \mbba_2}  {| \mbba_1   \times   \mbba_2 |} = \frac  { \mbba_1   \times   \mbba_2}  {\sqrt a}, \quad
\bfb= - \mbbn_{,\alpha} \otimes \mbba^{\alpha}
\end{split}
\right \}.
\end{equation}
Moreover, based on Eq.~\eqref{bfpsi}, the covariant basis vectors $\mbbg_i$ at $\mbbx$ are as follows:
\begin{equation}\label{bfg}
\mbbg_{\alpha}=\frac {\partial \mbbx}{\partial \zeta^{\alpha}}
=\mbba_{\alpha}+z^2 \phi_{,\alpha} \mbbd + z(1+z \phi) \mbbd_{,\alpha}, \quad
\mbbg_{3}=\frac{\partial \mbbx}{\partial z}=(1+2z \phi) \mbbd.
\end{equation}
It is observed from Eq.~\eqref{bfg}$_2$  that the director $\mbbd$ and the basis vector $\mbbg_{3}$ are in the same direction. However, the normal vector $\mbbn$ and $\mbbg_{3}$ are not in the same direction, in general.
Next, the vectors $\mbbu=u_i \mbbe_i$ and $\mbbw=w_i \mbbe_i$  are defined as the mid-surface and the director displacements, respectively.
This allows one to write
\begin{equation}\label{uw}
\overbar{\mbbx}=\overbar{\mbbX}+\mbbu, \quad
\mbbd=\mbbD+\mbbw.
\end{equation}
From Eqs.~\eqref{ShifterS} and \eqref{Gradians}$_1$, the deformation gradient tensor  $\bfF$, described  in the convective coordinate system   $\{ \zeta^i \}$, takes the form
\begin{equation}\label{bfF1}
\bfF=\Lgrad_{\bfzeta}{\mbbx}
=\frac{\partial \mbbx}{\partial \zeta^i} \otimes \mbbG^i 
= \mbbg_i \otimes \mbbG^i
= (\mbbg_i \otimes \mbbA^i) \bfQ^{-1}.
\end{equation}
In the present shell model, using Eqs.~\eqref{bfG}, \eqref{ShifterS}$_3$, \eqref{bfg},  and \eqref{bfF1}, and neglecting the higher-order terms involving $z^2$, the deformation gradient is approximated as follows:
\begin{equation}\label{bfF2}
\bfF \approx \tilde{\bfF} \bfQ^{-1}
\quad \text{with} \quad \tilde{\bfF} = \bfF^{[0]} + z \bfF^{[1]},
\end{equation}
and the tensors $\bfF^{[0]}$ and $\bfF^{[1]}$, with the aid of  Eqs.~\eqref{onM0}$_8$ and \eqref{uw}, are given by
\begin{equation}\label{F0F1_1}
\left.
\begin{split}
\bfF^{[0]}=\mbba_i \otimes \mbbA^i
=\mbba_{\alpha} \otimes \mbbA^{\alpha} + \mbbd \otimes \mbbD
=\bfI + \Lgrad_{\calS_0} \mbbu +\mbbw\otimes \mbbD\\
\bfF^{[1]}=\mbbd_{,\alpha} \otimes \mbbA^{\alpha} + 2 \phi \ \mbbd \otimes \mbbD
= \Lgrad_{\calS_0} \mbbw +2 \phi (\mbbD+\mbbw)\otimes \mbbD-\bfB
\end{split}
\right\}.
\end{equation}
To circumvent numerical difficulties in finite element solution, the in-plane deformation gradient term $\bfF^{[0]}$ is enhanced by the second-order tensor $\bar{\bfF}$, to be introduced in Section~\ref{FEM}.
Accordingly, the term $\bfF^{[0]}$ in Eq.~\eqref{bfF2}$_2$ is replaced by $\bfF^{[0]}+\bar{\bfF}$.
Moreover, following Ramezani and Naghadabadi~\cite{RamezBeam} in the context of the micropolar Timoshenko beam model, it is assumed that the micro-rotation pseudo-vector is constant along the shell thickness, namely 
$\bftheta = \tilde{\bftheta}(\zeta^1,\zeta^2)$.
Accordingly, the micro-rotation tensor $\tilde{\bfR}(\bftheta)$ is independent of the $z$ coordinate.
Keeping this in mind and using Eqs.~\eqref{UVtilde}$_1$,  \eqref{Gamma}$_1$, \eqref{ShifterS}$_3$,  \eqref{Gradians}$_2$, and \eqref{bfF2}, the micropolar deformation measures $\tilde{\bfU}$ and $\bfGamma$, in the present shell formulation, may be written as
\begin{equation}\label{UtildeGamma}
\begin{split}
\tilde{\bfU}=\tilde{\bfR}^{\trans} \bfF^{\ast} \bfQ^{-1}
=\tilde{\bfR}^{\trans} (\bfF^{[0]} + \bar{\bfF} + z \bfF^{[1]}) \bfQ^{-1}, \quad
\bfGamma= -\frac{1}{2} \calbfE \tendot 
\big[ \tilde{\bfR}^{\trans}  (\Lgrad_{\calS_0}\tilde{\bfR}) \bfQ^{-1} \big]. 
\end{split}
\end{equation}
where $\bfF^{\ast} = \tilde{\bfF} + \bar{\bfF}$ is the enhanced form of $\tilde{\bfF}$. 
The present formulation with $\{ \mbbu, \mbbw, \phi, \bftheta\}$ as the unknown field variables may be regarded as a 10-parameter micropolar shell model.
In other words, the present formulation is the extension of the classical 7-parameter shell model, with  $\{ \mbbu, \mbbw,  \phi\}$ as its unknowns, introduced by Sansour (e.g., \cite{Sansour98, SansKoll2000}).

\section{Variational formulation}
\label{PVW}
\vspace{-5pt}
Let $\delta \calU$ be the virtual internal energy and $\delta \calW$ denote the virtual work of external loads.
The principle of virtual work states that 
$\delta \calU-\delta \calW=0$ \cite{Wriggers2008}.
In what follows, the expressions for  $\delta \Psi$ and $\delta \hat{\calW}$, as, respectively, $\delta \calU$ and $\delta \calW$ per unit reference volume, are derived.
Moreover, for the linearization purpose to be used in the next section, the increments of $\delta \Psi$ and $\delta \hat{\calW}$ are also calculated.

By neglecting thermal effects, assuming that the material is hyperelastic, and to develop a material formulation, the internal energy per unit reference volume may be written as $\Psi=\tilde{\Psi}(\tilde{\bfU}, \bfGamma)$ \cite{Eringen1976, Steinmann1994}.
Using this point and Eq.~\eqref{delUG} furnishes
\begin{equation}\label{delPsi2}
\begin{split}
\delta \Psi &= \frac{\partial \Psi} {\partial \tilde{\bfU}} \tendot \delta \tilde{\bfU}
+ \frac{\partial \Psi} {\partial \bfGamma } \tendot \delta \bfGamma
=(\tilde{\bfR} \frac{\partial \Psi} {\partial \tilde{\bfU}} \bfF^{\trans}) \tendot (\delta \bfY - \delta \hat{\bfomega}) 
+(\tilde{\bfR} \frac{\partial \Psi} {\partial \bfGamma }\bfF^{\trans})  \tendot \Egrad \delta \bfomega.
\end{split}
\end{equation}
Moreover, the constitutive equations for the pairs 
$\{\tilde{\bfP} , \tilde{\bfM} \}$ and
$\{ \bfP,  \bfM \}$ are as follows \cite{DH2022IJSS}: 
\begin{equation}\label{strcop}
\begin{split}
\{\tilde{\bfP} , \tilde{\bfM} \} =\bigg\{  \frac{\partial \Psi} {\partial \tilde{\bfU}} ,\frac{\partial \Psi} {\partial \bfGamma} \bigg\}, \quad
\{ \bfP,  \bfM \} = \tilde{\bfR}  \bigg\{  \frac{\partial \Psi} {\partial \tilde{\bfU}} ,\frac{\partial \Psi} {\partial \bfGamma} \bigg\}.
\end{split}
\end{equation}
It is noted that Eqs.~\eqref{delPsi2} and \eqref{strcop} hold for all three-dimensional micropolar hyperelastic solids. 
For the present shell model, first the quantities denoted by $\delta \bfUpsilon^{[\calJ]}$ ($\calJ=1,2$) are defined by
\begin{equation}\label{Y12}
\begin{split}
\delta \bfUpsilon^{[1]} 
=\delta \bfF^{[0]}+\delta \bar{\bfF} 
+z\delta \bfF^{[1]}
- \delta \hat{\bfomega} \bfF^{\ast}, \quad
\delta \bfUpsilon^{[2]} =\Lgrad_{\calS_0} \delta \bfomega.
\end{split}
\end{equation}
Next, after replacing $\tilde{\bfF}$ by the enhanced form $\bfF^{\ast}$, combination of Eqs.~\eqref{UtildeGamma}, \eqref{delPsi2}$_1$,  \eqref{strcop}$_{1,2}$, and \eqref{Y12} leads to the following expression for $\delta \Psi$:
\begin{equation}\label{PsiPK1}
\begin{split}
\delta \Psi &= \bfP^{[0]} \tendot \delta \bfUpsilon^{[1]} 
+\bfM^{[0]} \tendot \delta \bfUpsilon^{[2]}
\quad \text{with} \quad
\{\bfP^{[0]}, \bfM^{[0]} \}=\{\bfP  , \bfM  \} \bfQ^{-1}.
\end{split}
\end{equation}
For linearization purpose, the increment of $\delta \Psi$ under the increment of the field variables $\Delta \mbbu$, $\Delta \mbbw$, $\Delta \phi$, and $\Delta \bftheta$ is needed.
Accordingly, from Eqs.~\eqref{delPsi2}$_1$, \eqref{strcop}, and \eqref{PsiPK1} it is deduced that
\begin{equation}\label{DdPsi}
\begin{split}
\Delta \delta \Psi =&
\bfP^{[0]} \tendot \Delta \delta \bfH^{[1]}
+\bfM^{[0]} \tendot \Delta \delta \bfH^{[2]}
+\delta\bfUpsilon^{[1]} \tendot (\calbfC^{[1]} \tendot \Delta \bfUpsilon^{[1]} +  \calbfC^{[2]} \tendot \Delta \bfUpsilon^{[2]})\\
&+\delta\bfUpsilon^{[2]} \tendot (\calbfC^{[3]} \tendot \Delta \bfUpsilon^{[1]}+ \calbfC^{[4]} \tendot \Delta \bfUpsilon^{[2]}),
\end{split}
\end{equation}
where $\Delta \delta \bfH^{[1]}$ and $\Delta \delta \bfH^{[2]}$ are as follows:
\begin{equation}\label{Q12}
\left.
\begin{split}
\Delta \delta \bfH^{[1]} = \frac{1}{2} 
(\Delta \hat{\bfomega} \delta \hat{\bfomega}
+\delta \hat{\bfomega} \Delta \hat{\bfomega}) \bfF^{\ast}
-(\delta \hat{\bfomega} \Delta \bfF^{\ast}
+\Delta \hat{\bfomega} \delta \bfF^{\ast}) \\
\Delta \delta \bfH^{[2]} = -\frac{1}{2}(
\Delta \hat{\bfomega} \Lgrad_{\calS_0} \delta \bfomega
+\delta \hat{\bfomega} \Lgrad_{\calS_0} \Delta \bfomega)
\end{split}
\right\}.
\end{equation}
Moreover, the fourth-order tensors $\calbfC^{[\calI]}$ ($\calI=1,2,3,4$) have the following components:
\begin{equation}\label{A14}
\begin{split}
\calC^{[\calI]}_{iJkL} = \tilde{R}_{iP} \tilde{R}_{kQ} Q^{-1}_{JR} Q^{-1}_{LS} 
\tilde{\calC}^{[\calI]}_{PRQS}
\quad \text{with} \quad 
Q^{-1}_{JR}= (\bfQ^{-1})_{JR},
\end{split}
\end{equation}
and $\tilde{\calC}^{[\calI]}_{PJQL}$ are the components of the following fourth-order tensors:
\begin{equation}\label{A142}
\begin{split}
 \tilde{\calbfC}^{[1]} = \frac{\partial^2 \Psi} 
{\partial \tilde{\bfU} \partial \tilde{\bfU} }, \quad
\tilde{\calbfC}^{[2]}=\frac{\partial^2 \Psi} 
{\partial \tilde{\bfU} \partial \bfGamma },\quad
\tilde{\calbfC}^{[3]}=\frac{\partial^2 \Psi} 
{\partial \bfGamma \partial \tilde{\bfU} } , \quad
\tilde{\calbfC}^{[4]} = \frac{\partial^2 \Psi} 
{\partial \bfGamma \partial \bfGamma }.
\end{split}
\end{equation}
In this work, the micropolar-enhanced neo-Hookean constitutive model proposed in Ref.~\cite{DH2022IJSS} is employed, according to which the free energy density is given by
\begin{equation}\label{nH1}
\Psi=(\eta+\half \mu) \tr (\tilde{\bfU} \tilde{\bfU}^{\trans})
-\eta \tr (\tilde{\bfU} ^2)
+\frac{1}{2}\lambda (\ln J)^2-\mu \ln J
+ \frac{1}{2} \mu l^2 \tr ( \bfGamma \bfGamma^{\trans} ),
\end{equation}
where $\eta$ is a material constant, and $l$ is the well-known material length-scale parameter.
%
%
From Eqs.~\eqref{strcop}$_1$ and \eqref{nH1}, the expressions for  $\tilde{\bfP}$ and $\tilde{\bfM}$ are then calculated to be~\cite{DH2022IJSS}
\begin{equation}\label{PK11}
\tilde{\bfP} = (\mu+\eta)\tilde{\bfU} -\eta \tilde{\bfU}^{\trans}
-(\mu - \lambda \ln J) \tilde{\bfU}^{-\trans}, \quad
\tilde{\bfM}= \mu l^2 \bfGamma.
\end{equation}
Moreover, Eqs.~\eqref{A142} and \eqref{nH1} lead to the following 4th-order tensors 
$\tilde{\calbfC}^{[\calI]}$:
\begin{equation}\label{A12}
\left.
\begin{split}
\tilde{\calbfC}^{[1]}
=(\mu+\eta) \bfI \odot \bfI   
+\lambda  \tilde{\bfU}^{-\trans} \otimes \tilde{\bfU}^{-\trans}
-\eta \bfI \boxtimes \bfI
-(\lambda \ln J-\mu) \tilde{\bfU}^{-\trans} \boxtimes \tilde{\bfU}^{-\trans}\\
\tilde{\calbfC}^{[2]} = \tilde{\calbfC}^{[3]} = {\bf0}, \qquad
\tilde{\calbfC}^{[4]}
=\mu l^2 \bfI \odot \bfI
\end{split}
\right\}.
\end{equation}
Next, it is recalled that $\mbbB^{\text{ext}}$  induces the body couple density $\mbbp^{\ast}$ on an HMSM (c.f. Eq.~\eqref{couple}$_2$). 
Since $\mbbp^{\ast}$ and $\bftheta$ are work-conjugate quantities, the virtual external work density $\delta \hat{\calW}$ and its increment may be written as
\begin{equation}\label{dWext}
\left.
\begin{split}
\delta \hat{\calW} 
=\frac{1}{\mu_0} [(\bfF \tilde{\mbbB}^{\text{rem}}) \times \mbbB^{\text{ext}}] \cdot \delta \bftheta\\
\Delta \delta \hat{\calW} = \frac{1}{\mu_0}  [(\Delta \bfF^{\ast}  \bfQ^{-1} \tilde{\mbbB}^{\text{rem}}) \times \mbbB^{\text{ext}}] \cdot \delta \bftheta
\end{split}
\right\}.
\end{equation}
As will be shown in the next section, Eq.~\eqref{dWext}$_2$ leads to the expression for the load stiffness matrix.
The expressions for $\delta \Psi$,   $\Delta \delta \Psi$, $\delta \hat{\calW}$, and $\Delta \delta \hat{\calW}$ given in Eqs.~\eqref{PsiPK1}, \eqref{DdPsi}, and \eqref{dWext} are the basic relations for the FE formulation presented in the next section.


\section{FE formulation}
\label{FEM}
\vspace{-5pt}
A nonlinear finite element formulation for the present shell model is developed in this section.
Let $\calS_0^{\fre}$ be a typical element in the referential mid-surface  $\calS_0^{\fre}$.
To perform numerical integration, the typical element is mapped to the two-dimensional parent square element $\square =[-1,1] \times [-1,1]$ in the  $\{ \xi, \eta \}$ space, with $\xi , \eta  \in [-1,1]$.
The field variables $\{u_i, w_i, \theta_i, \phi\}$, over the parent element $\calS_0^{\fre}$, are interpolated as follows:
\begin{equation}\label{interpol}
\begin{split}
u_i= \mbbN_{u} \mbbU_i, \quad
w_i=\mbbN_{w} \mbbW_i, \quad
\theta_i=\mbbN_{w}      \bfTheta_{i},  \quad
\phi=\mbbN_{\phi} \bfPhi, 
\end{split}
\end{equation}
where $\mbbN_{u}=\{N_u^1,N_u^2, ..., N_u^{n_u} \}$ is a row vector containing the shape  functions that interpolate the mid-surface displacement $u_i$ over the element.
Here, $n_u$ is the number of nodes of the element that possess the $u_i$-DOF.
Let $U_i^{\calI}$ be the displacement component $u_i$ at the $\calI$'th node ($\calI=1,2,...,n_u$) of the element.
Accordingly, $\mbbU_i=\{U_i^1, U_i^2,..., U_i^{n_u}\}^{\trans}$ is a column vector that involves all $U_i^{\calI}$'s over the element.
Similar definitions hold for the other quantities in Eq.~\eqref{interpol}.
Moreover, similar relations hold for the increment 
$\{\Delta u_i,\Delta w_i,\Delta \theta_i,\Delta \phi\}$ and variation 
$\{\delta u_i,\delta w_i,\delta \theta_i,\delta \phi\}$ of the field variables.
The generalized displacement vector ${\mbbv}^{\fre}$ involving all nodal DOFs of the typical element may be written as
\begin{equation}\label{mbbv}
{\mbbv}^{\fre}_{n^{\fre} \times 1}
=\{\mbbU_1^{\trans},\mbbU_2^{\trans}, \mbbU_3^{\trans},
\mbbW_1^{\trans},\mbbW_2^{\trans},\mbbW_3^{\trans}, 
\bfTheta_1^{\trans},
\bfTheta_2^{\trans},
\bfTheta_3^{\trans},
\bfPhi^{\trans} \}^{\trans},
\end{equation}
where $n^{\fre}=3(n^u+n^w+n^{\theta})+n^{\phi}$ is the number of nodal DOFs.
Based on Esq.~\eqref{dbfF}, \eqref{F0F1_1}, \eqref{Y12}, and \eqref{interpol}, 
the following relations hold:
\begin{equation}\label{PolFs}
\left.
\begin{split}
\delta F^{[0]}_{iJ}
= A^{\ast \alpha J} \mbbN_{u,\alpha} \delta \mbbU_{i} 
+ \mbbN_w D_J \delta \mbbW_{i}
=\mbbb^{[0]}_{iJ} \delta  {\mbbv}^{\fre}\\
\delta F^{[1]}_{iJ}=
  (A^{\ast \alpha J} \mbbN_{w,\alpha} 
+2 \phi  D_J \mbbN_w ) \delta \mbbW_{i}
+2 d_i D_j    \mbbN_{\phi}  \delta \bfPhi
= \mbbb^{[1]}_{iJ} \delta  {\mbbv}^{\fre}\\
(\delta \hat{\bfomega} \bfF^{\ast})_{iJ}=
\epsilon_{ijk} \Lambda_{kp}  F^{\ast}_{jJ}
\mbbN_{\theta} \delta {\bfTheta}_p
=\mbbb^{[\omega F]}_{iJ} \delta  {\mbbv}^{\fre}\\
\delta Y^{[2]}_{iJ}=
A^{\ast \alpha J}  (\Lambda_{ip}  \mbbN_{\theta} )_{,\alpha}   
\delta \bfTheta_{p}
= \mbbb^{[2]}_{iJ} \delta  {\mbbv}^{\fre}\\
\end{split}
\right\}.
\end{equation}
Here, the last equality in each relation indicates that all components can be expressed in terms of the generalized virtual displacement vector 
$\delta {\mbbv}^{\fre}$.
Next, the enhanced deformation gradient tensor $\bar{\bfF}$ is considered.
Let  $\bfalpha= \{\alpha_1, \alpha_2, ..., \alpha_{\calP^{\ast} } \}^{\trans}$ be the column vector of enhanced parameters $\alpha_{\calP}$'s ($\calP=1,2,...,\calP^{\ast}$), where $\calP^{\ast}$ is the total number of enhanced parameters.
The components of $\bar{\bfF}$ and its variation/increment depend linearly on $\alpha_{\calP}$'s (see, e.g., Refs.~\cite{SimoArmero92, Simo93, KW96, Glaser97}). Here, following the notation used in Eq.~\eqref{PolFs}, one may write
\begin{equation}\label{Fbarcomp}
\begin{split}
\{\bar{F}_{iJ}, \delta \bar{F}_{iJ}, \Delta \bar{F}_{iJ} \} 
= \bar{\mbbb}_{iJ} \{ \bfalpha, \delta \bfalpha, \Delta \bfalpha \}
\quad \text{or} \quad
\{\bar{\mbbF}, \delta \bar{\mbbF}, \Delta \bar{\mbbF} \} 
=\bar{\mbbB}  \{ \bfalpha, \delta \bfalpha, \Delta \bfalpha \},
\end{split}
\end{equation}
where $\bar{\mbbb}_{iJ}$ are the $\calP^{\ast}  \times 1$ row vectors, $\bar{\mbbF}$ is the $9 \times 1$ vectorial representation of $\bar{\bfF}$, and $\bar{\mbbB}$ is a $9 \times \calP^{\ast}$ matrix the rows of which are $\bar{\mbbb}_{iJ}$.

Now, let $\mbbB^{(\calN)}$ ($\calN=0,1,2$) and $\mbbB^{[\omega F]}$ be the 
$9 \times n^{\fre}$ matrices whose rows are $\mbbb^{(\calN)}_{iJ}$ and $\mbbb^{[\omega F]}_{iJ}$, respectively.
From Eqs.~\eqref{Y12}, \eqref{PolFs}, and  \eqref{Fbarcomp} it then follows that
\begin{equation}\label{Y12vec}
\delta \mbbY^{[1]}=\tilde{\mbbB} \delta {\mbbv} 
+ \bar{\mbbB}  \delta \bfalpha 
\quad \text{and} \quad
\delta \mbbY^{[2]}=\mbbB^{[2]} \delta {\mbbv} 
\quad \text{with} \quad
\tilde{\mbbB}= \mbbB^{[0]}+\mbbB^{[\omega F]}+z \mbbB^{[1]}.
\end{equation}
The differential volume element  $\text{d}\calV^{\fre}_0$ located at the elevation $z$ with respect to the typical element $\calS^{\fre}_0$ is given by (e.g., \cite{Sansour98})
\begin{equation}\label{dV0}
\text{d}\calV^{\fre}_0=Q \text{d}\calS^{\fre}_0 \text{d}z
\quad \text{with} \quad
Q=\det \bfQ
\quad \text{and} \quad
\text{d}\calS^{\fre}_0 =\sqrt{A} \text{d} \zeta^1  \text{d} \zeta^2.
\end{equation}
Now, the element virtual internal energy is given by $\delta \calU^{\fre} = \int_{\calV_0^{\fre}} \Psi  \text{d} \calV_0^{\fre}$.
Similarly, the expression for the virtual work over the element can be calculated via $\delta \calW^{\fre} = \int_{\calV_0^{\fre}} \hat{\calW} \text{d} \calV_0^{\fre}$.
Using Eqs.~\eqref{Y12} and \eqref{Y12vec}, the expressions for $\delta \calU^{\fre}$ and $\delta \calW^{\fre}$ may be written as
\begin{equation}\label{UintWext1}
\delta \calU^{\fre} = \delta {\mbbv}^{\trans} \mbbF_{\text{int}}^{v } 
 +\delta \bfalpha^{\trans} \mbbF_{\text{int}}^{\alpha}, \quad
\delta \calW^{\fre}=
\delta \bfTheta_i^{\trans} \mbbF^{\theta}_{ \text{ext} i} 
=\delta {\mbbv}^{\trans} \mbbF_{\text{ext}}^{v } , 
\end{equation}
where the internal force vectors $\mbbF^{\text{int} v }$ and 
$\mbbF^{\text{int} \alpha }$  are as follows:
\begin{equation}\label{Fint1}
\begin{split}
\mbbF_{\text{int}}^{v } = \int_{\calV_0^{\fre}} 
\big( \tilde{\mbbB}^{\trans} \mbbP^{[0]}
+{\mbbB}^{[2]\trans} \mbbM^{[0]} \big) \text{d} \calV_0^{\fre}, \quad
\mbbF_{\text{int}}^{\alpha} = \int_{\calV_0^{\fre}} 
\bar{\mbbB}^{\trans} \mbbP^{[0]} \text{d} \calV_0^{\fre}.
\end{split}
\end{equation}
Moreover, the external force  vector $\mbbF^{\theta}_{ \text{ext} i}$, work conjugate to $\bfTheta_i$, is given by
\begin{equation}\label{Fext1}
\mbbF^{\theta}_{ \text{ext} i}  =\frac{1}{\mu_0} \int_{V_0^{\fre}} 
\epsilon_{imj}  F_{mJ} \tilde{B}^{\text{rem}}_J B^{\text{ext}}_j  \mbbN_{\theta}^{\trans}
\text{d} \calV_0^{\fre}.
\end{equation}
Next, the linearized equations resulting from Eqs.~\eqref{DdPsi}, \eqref{dWext}$_2$, and \eqref{UintWext1} may be written as 
\begin{equation}\label{DVW}
\Delta \delta \calU^{\fre} - \Delta \delta \calW^{\fre} 
=  -(\delta \calU^{\fre} - \delta \calW^{\fre}).
\end{equation}
The system of algebraic equations extracted from Eq.~\eqref{DVW} may be written as
\begin{equation}\label{kuf}
\begin{split}
\begin{bmatrix}
\mbbK^{v v}_{\text{mat}}
+\mbbK^{v v}_{\text{geo}}
-\mbbK^{v v}_{\text{load}} &  
\mbbK^{v \alpha}_{\text{mat}}
+\mbbK^{v \alpha}_{\text{geo}}
-\mbbK^{v \alpha}_{\text{load}}\\
\mbbK^{v \alpha \trans}_{\text{mat}}
+\mbbK^{\alpha v}_{\text{geo}} &  
\mbbK^{\alpha \alpha}_{\text{mat}}\\
\end{bmatrix}
\begin{Bmatrix}
\Delta {\mbbv} \\
\Delta \bfalpha
\end{Bmatrix}=-
\begin{Bmatrix}
\mbbF_{\text{int}}^{v } - \mbbF_{\text{ext}}^{v } \\
\mbbF_{\text{int}}^{\alpha}
\end{Bmatrix},
\end{split}
\end{equation}
where the subscripts "mat", "geo", and "load", represent the material, geometric, and load part of the element stiffness matrix. 
In particular, the material sub-matrices $\mbbK^{v v}_{\text{mat}}$, 
$\mbbK^{v \alpha}_{\text{mat}}$, and 
$\mbbK^{\alpha \alpha}_{\text{mat}}$ in Eq.~\eqref{kuf} are as follows:
%
%
\begin{equation}\label{Kmatuu}
\left.
\begin{split}
\mbbK^{v v}_{\text{mat}} = \int_{V_0^{\fre}}
[\tilde{\mbbB}^{\trans}(\mbbC^{[1]}\tilde{\mbbB}+\mbbC^{[3]}\mbbB^{[2]})
+\mbbB^{[2] \trans} (\mbbC^{[2]}\mbbB^{[2]}+\mbbC^{[4]} \tilde{\mbbB})]
\text{d} \calV_0^{\fre}\\
\mbbK^{v \alpha}_{\text{mat}} 
= \int_{V_0^{\fre}}
(\tilde{\mbbB}^{\trans} \mbbC^{[1]}
+ \mbbB^{[2] \trans} \mbbC^{[4]}) \bar{\mbbB} 
\text{d} \calV_0^{\fre}, \quad
\mbbK^{\alpha \alpha}_{\text{mat}} 
= \int_{V_0^{\fre}}
\bar{\mbbB} ^{\trans} \mbbC^{[1]} \bar{\mbbB} 
\text{d} \calV_0^{\fre}
\end{split}
\right\},
\end{equation}
where  $\mbbC^{[\calJ]}$ are the matrix forms of $\calbfC^{[\calJ]}$.
Moreover, the load sub-matrices $\mbbK^{v v}_{\text{load}}$, 
$\mbbK^{v \alpha}_{\text{load}}$ are given by
\begin{equation}\label{Kmatuu}
\left.
\begin{split}
\mbbK^{v v}_{\text{load}} = \int_{V_0^{\fre}}
\epsilon_{ijk} Q^{-1}_{JN} \tilde{B}^{\text{rem}}_N B^{\text{ext}}_j 
\mbby_k (\mbbb^{[0]}_{iJ}+z \mbbb^{[1]}_{iJ}) 
\text{d} \calV_0^{\fre}\\
\mbbK^{v \alpha}_{\text{load}} = \int_{V_0^{\fre}}
\epsilon_{ijk} Q^{-1}_{JN} \tilde{B}^{\text{rem}}_N B^{\text{ext}}_j 
\mbby_k \bar{\mbbb}_{iJ} \text{d} \calV_0^{\fre}
\end{split}
\right\},
\end{equation}
where $\mbby_k=\{ {\bf0}_{1 \times 3(n^u+n^w)}, 
{\bf0}_{1 \times (k-1)n^{\theta}},
\mbbN_{\theta}, 
{\bf0}_{1 \times (3-k)n^{\theta}}, 
{\bf0}_{1 \times n^{\phi}} \}^{\trans}$, with $k \in \{1,2,3\}$, is a column vector whose nonzero entry is $\mbbN_{\theta}$.
The expressions for the geometric sub-matrices, resulting from the terms 
$\bfP^{[0]} \tendot \Delta \delta \bfH^{[1]}$ and $\bfM^{[0]} \tendot \Delta \delta \bfH^{[2]}$ in Eq.~\eqref{DdPsi}, are too lengthy and are not presented here.
The assembled system of equations is of the form 
$\tilde{\mbbK} \Delta \tilde{\mbbV}=-\tilde{\mbbR}$, where 
$\tilde{\mbbK}$,  $\Delta \tilde{\mbbV}$, and $\tilde{\mbbR}$ are the assembled forms of the stiffness matrix, incremental generalized displacement, and residual vector, respectively.
%
%
After finding $\Delta \tilde{\mbbV}$, the non-rotational quantities are update via the relations $\mbbu+\Delta \mbbu \rightarrow \mbbu$, 
$\mbbw+\Delta \mbbw \rightarrow \mbbw$,
and $\phi+\Delta \phi \rightarrow \phi$.
However, the update procedure for the rotation pseudo-vector is completely different. 
Let $\Delta \bftheta$ be the increment of the rotation pseudo-vector.
The updated rotation pseudo-vector ${\bftheta}^{\ast}_{\text{updated}}$ resulting from the two subsequent rotations $\bftheta$ and $\Delta \bftheta$  is then calculated via the following relations~\cite{Argyris1982}:
\begin{equation}\label{newth}
{\bftheta}^{\ast}_{\text{updated}} = \frac{ {\bftheta}^{\ast} + \Delta {\bftheta}^{\ast} + ({\Delta \hat{\bftheta}}^{\ast}) {\bftheta}^{\ast} }
{1-{\bftheta}^{\ast} \cdot {\Delta \bftheta}^{\ast} },
\end{equation}
where ${\Delta \hat{\bftheta}}^{\ast} = - \calbfE \Delta {\bftheta}^{\ast}$  and ${\bftheta}^{\ast}= \bfalpha \tan \frac{\theta}{2}$. 
Moreover, $\bfalpha= {\bftheta} / {\theta}$ is the unit vector along $\bftheta$.
The proof of Eq.~\eqref{newth} is lengthy and is available in, e.g., Argyris \cite{Argyris1982}.


\vspace{10pt}
\section{Numerical examples}
\label{Examples}
\vspace{-5pt}
To examine the applicability of the developed formulation, six examples are solved in this section.
The formulation has been implemented in our in-house finite element code.
The 10-parameter micropolar shell element designed for the present numerical simulations is an eight-node quadrilateral.
All eight nodes contain the three displacement components $u_i$.
However, only the corner nodes contain the $w_i$, $\phi$, and $\theta_i$ DOFs.
In other words, the DOF parameters defined after Eqs.~\eqref{interpol} and \eqref{mbbv} are $n^u=8$ and $n^w=n^{\theta}=n^{\phi}=4$. 
Following Korelc and Wriggers \cite{KW96}, the enhancing deformation gradient $\bar{\bfF}$ is considered to be of the following form:
\begin{equation}\label{barF}
\bar{\bfF}= \bfJ^{-\trans} \bar{\bfF}^{\text{ref}} \bfJ^{-1},
\end{equation}
where $\bfJ$ is the Jacobi matrix between the physical and parent elements.
Moreover, $\bar{\bfF}^{\text{ref}}$ is the enhancing deformation  gradient defined in the parent $\{ \xi, \eta \}$ space. 
In this work, the nonzero components of $\bar{\bfF}^{\text{ref}}$ are considered  as follows:
\begin{equation}\label{Fbarref}
\bar{F}^{\text{ref}}_{13}=\alpha_1 \xi + \alpha_2 \eta, \quad
\bar{F}^{\text{ref}}_{23}=\alpha_3 \xi + \alpha_4 \eta, \quad
\bar{F}^{\text{ref}}_{33}=\alpha_5 \xi + \alpha_6 \eta,
\end{equation}
which are linear functions in terms of the parent coordinates $\xi$ and $\eta$.
This indicates that $\bar{\bfF}$ contains six enhanced parameters, namely $\calP^{\ast}=6${\footnote{It is also possible to include the nonlinear terms involving $\{ \xi^2, \eta^2, \xi \eta^2,\xi^2 \eta \}$ or 
$\{ 1-3\xi^2, 1-3\eta^2, \xi (1-3\eta^2),\eta(1-3\xi^2)  \}$ in the components of $\bar{\bfF}^{\text{ref}}$, which increases the number of the enhanced parameters to $39$.
However, our numerical simulations reveal that the change in the results is negligible.}.
To evaluate the integrals over the element surface, the $2 \times 2$ Gauss--Legendre integration is used.
Moreover, the two-point rule is employed for integration along the shell thickness.
%
\subsection{VERIFICATION EXAMPLE: bending of beam-like strips}
\label{Cantilever}
\vspace{-5pt}
To examine the validity of the results of the proposed formulation, the flexural deformation of four beam-like strips under magnetic loading is studied in this example.
Extensive experiments on these structures have been previously conducted by Zhao et al.~\cite{Zhao2019}.
The values of the mechanical properties $\lambda$ and $\mu$ are, respectively, $7300$ and $303$ (in kPa).
As can be seen from Fig.~\ref{Ex1Fig1}, the length of the undeformed strips is considered to be along the $X_1$ axis.
The referential remnant magnetic flux can be described by the vector $\tilde{\mbbB}^{\text{rem}}= 143 \mbbE_1$ (mT).
The width of all strips is $5$ mm.
The length $L$ (in mm), the height $h$ (in mm), and the aspect ratio $AR=L/h$ of the strips are given by ($L,h,AR$)= ($11,1.1,10$), ($19.2,1.1,17.5$), ($17.2,0.84,20.5$), and ($17.2, 0.42,41$).
The strips are clamped at $X_1=0$, and the maximum applied magnetic loading is $\mbbB^{\text{ext}}_{\text{max}}= 50 \mbbe_3$ (mT).
Convergence analysis reveals that the minimum required number of elements along the length of the strips is $10$, $15$, $30$, and $40$, respectively.
Additionally, two elements in the width direction are necessary for the four strips. 
Furthermore, for the micropolar parameter $\eta = \mu/10$ and the material length-scale $l=h/10$, the present results are in good agreement with the available data obtained in \cite{Zhao2019}.
Therefore, these relations will be employed for the next examples, as well.

Fig.~\ref{Ex1Fig1}(\textbf{a}) displays the nonrmalized deflection   $u_3^{\text{T}}/L$ at the tip of strips versus the nondimensional load $\frac{10^3}{\mu \mu_0}|\mbbB^{\text{ext}}||\tilde{\mbbB}^{\text{rem}}|$.
From the figure, it is clear that the results based on the present shell formulation are close to the numerical as well as experimental data reported in Ref.~\cite{Zhao2019}.
The deformation patterns of the strips for four values of the external magnetic flux are displayed in Figs.~\ref{Ex1Fig1}(\textbf{b,c,d,e}).
To have a comparison between the deformation of the strips for a specific value of $|\mbbB^\text{ext}|$,  the four strips are plotted in the same figure.
The importance of the aspect ratio can be observed in Fig.~\ref{Ex1Fig1}(\textbf{b}), where the strip with $AR=41$ experiences considerable large deformation even for $|\mbbB^\text{ext}|=2$ (mT), which is a small value for the applied magnetic flux.
\begin{figure*}
\centering
\begin{minipage}[b]{.5\textwidth}
\includegraphics[width=80mm]{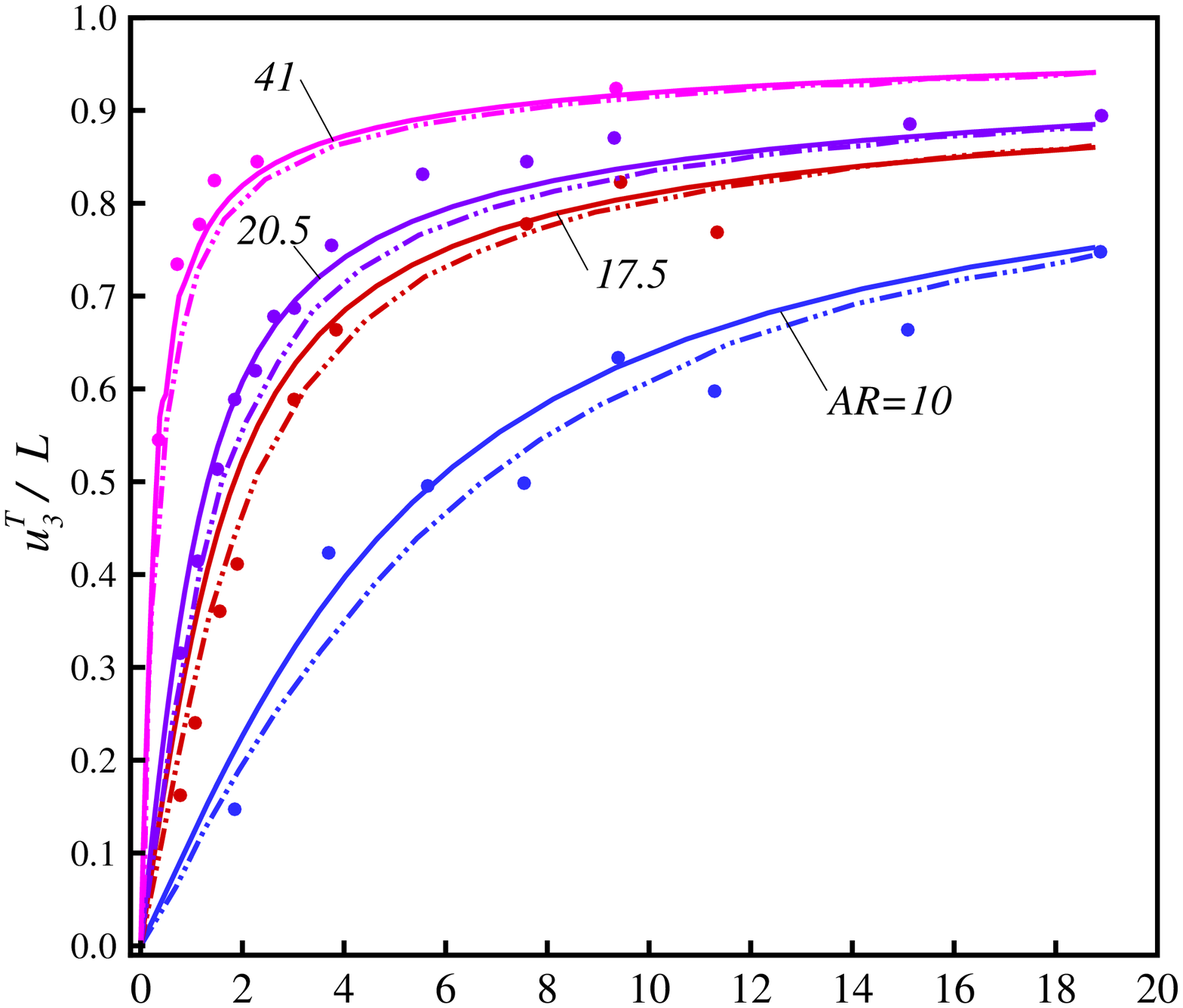} 
\put(-50,90){\scriptsize(\textbf{a})}
\put(-150,4){\scriptsize $\frac{1}{\mu \mu_0} |\mbbB^{\text{ext}}||\tilde{\mbbB}^{\text{rem}}| \times 10^3$}
\put(-130,55){\scriptsize $\bullet$ Experiment, Ref.~\cite{Zhao2019} }
\put(-130,40){\scriptsize $- \cdot \cdot-$ FEM, Ref.~\cite{Zhao2019} }
\end{minipage}
\begin{minipage}[b]{.4\textwidth}
\includegraphics[trim=0cm -3cm  0cm  0cm, width=50mm]{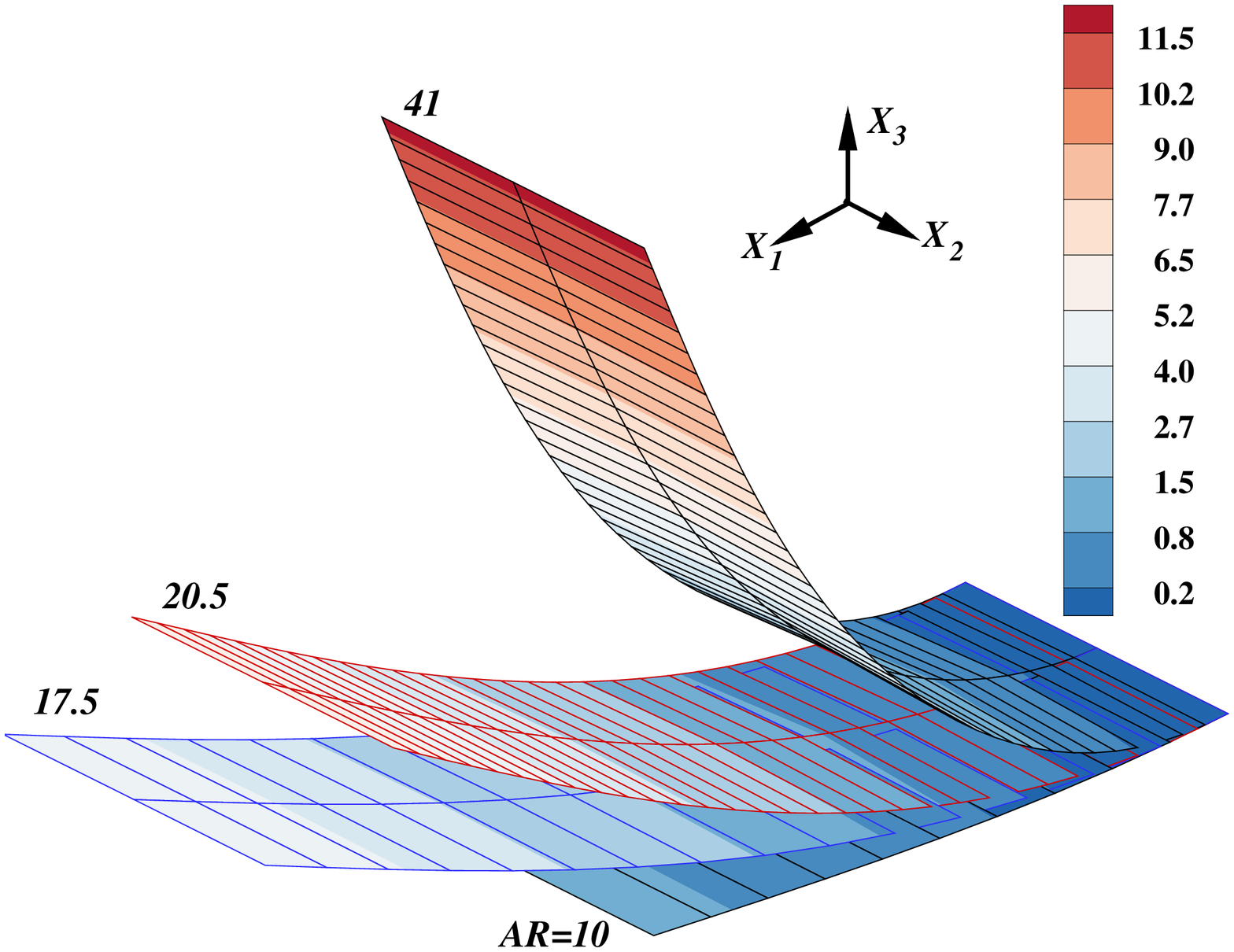} 
\put(-70,140){\scriptsize(\textbf{b})}
\end{minipage} \\
\vspace{10pt}
\begin{minipage}[b]{.32\textwidth}
\includegraphics[trim=0cm 2cm  0cm  0cm, width=50mm]{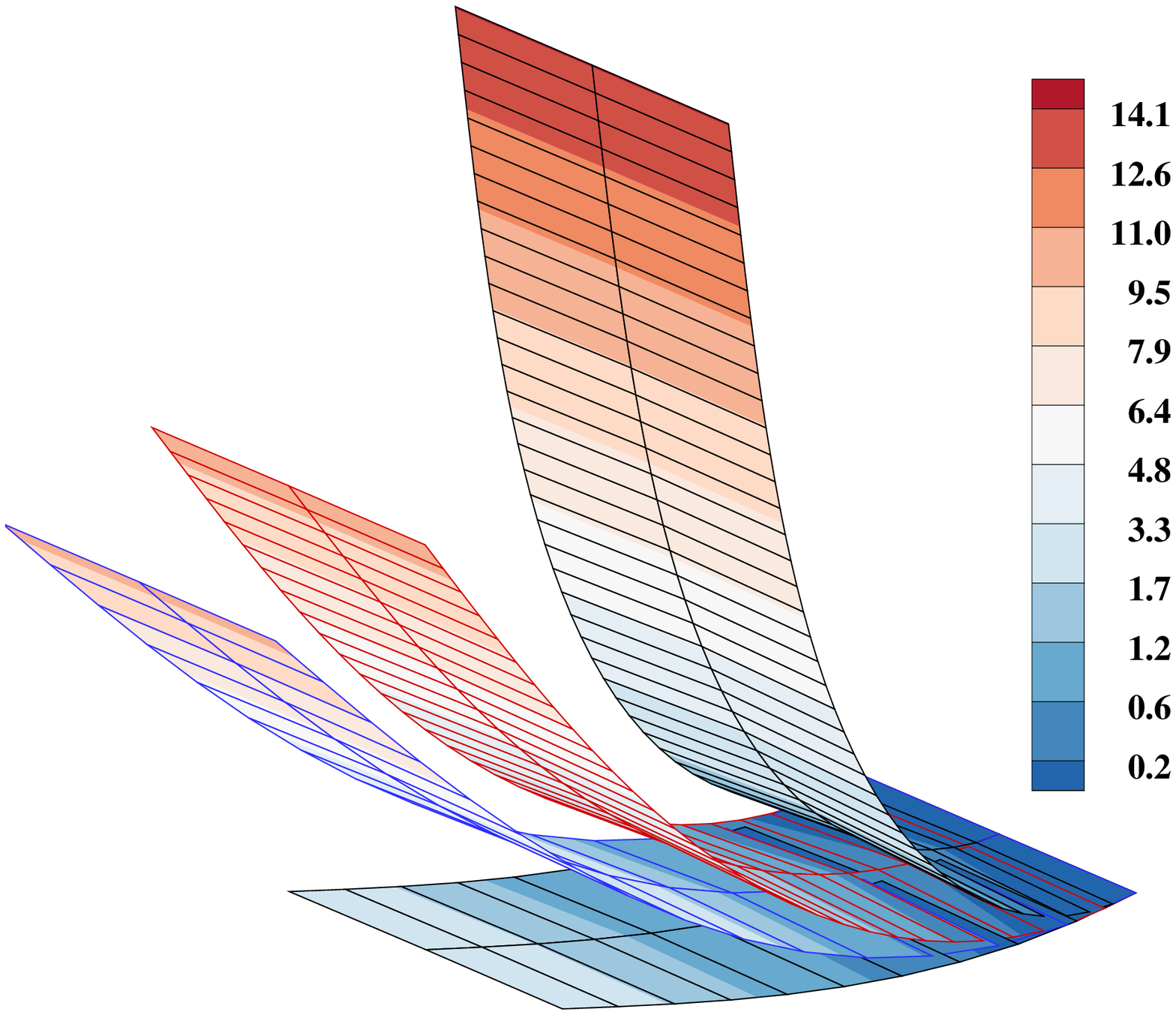} 
\put(-60,120){\scriptsize(\textbf{c})}
\end{minipage}   
\begin{minipage}[b]{.32\textwidth}
\includegraphics[trim=0cm 5mm  0cm  0cm, width=50mm]{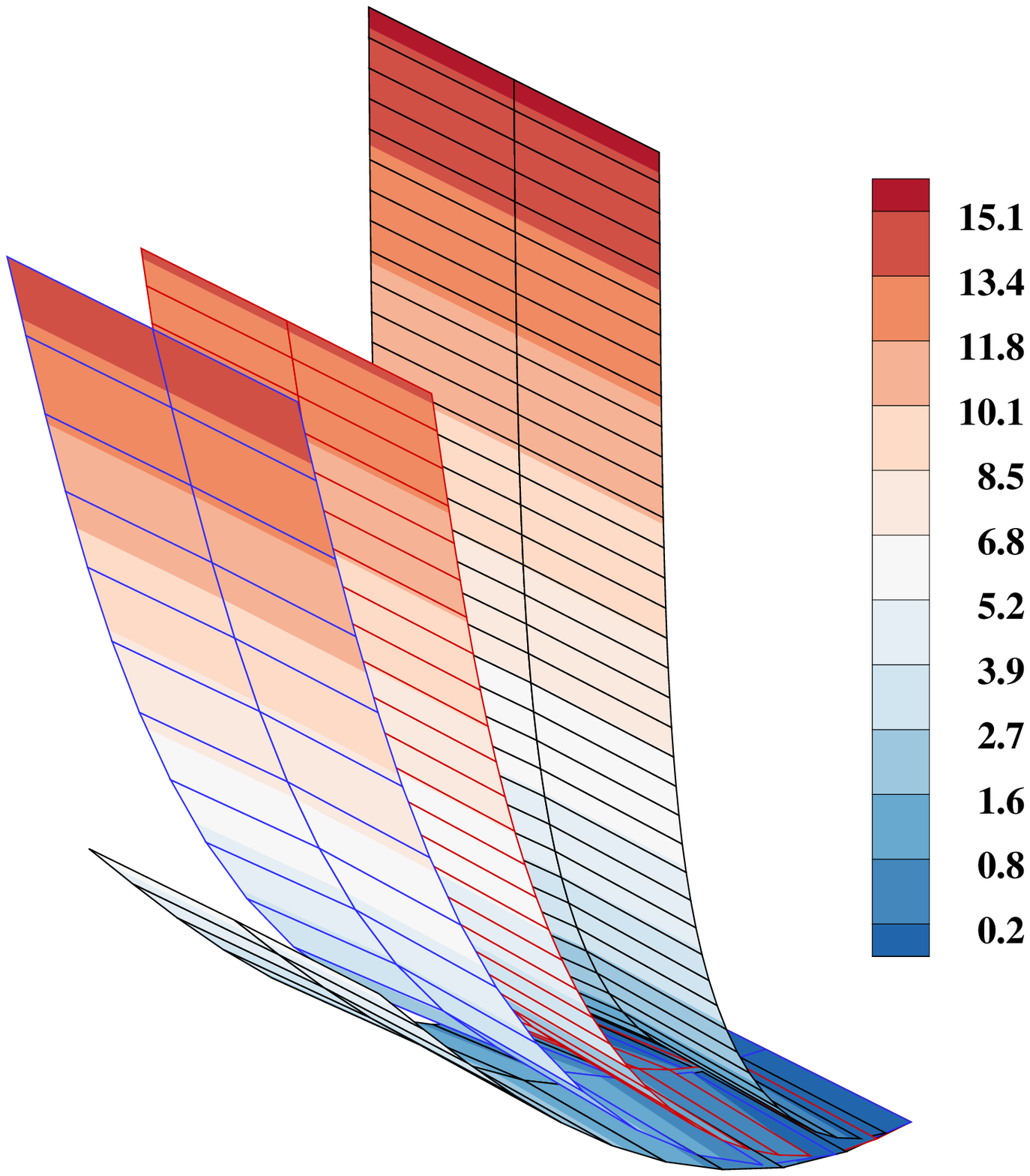} 
\put(-65,120){\scriptsize(\textbf{d})}
\end{minipage}
\begin{minipage}[b]{.32\textwidth}
\includegraphics[trim=0cm 0cm  0cm  0cm, width=50mm]{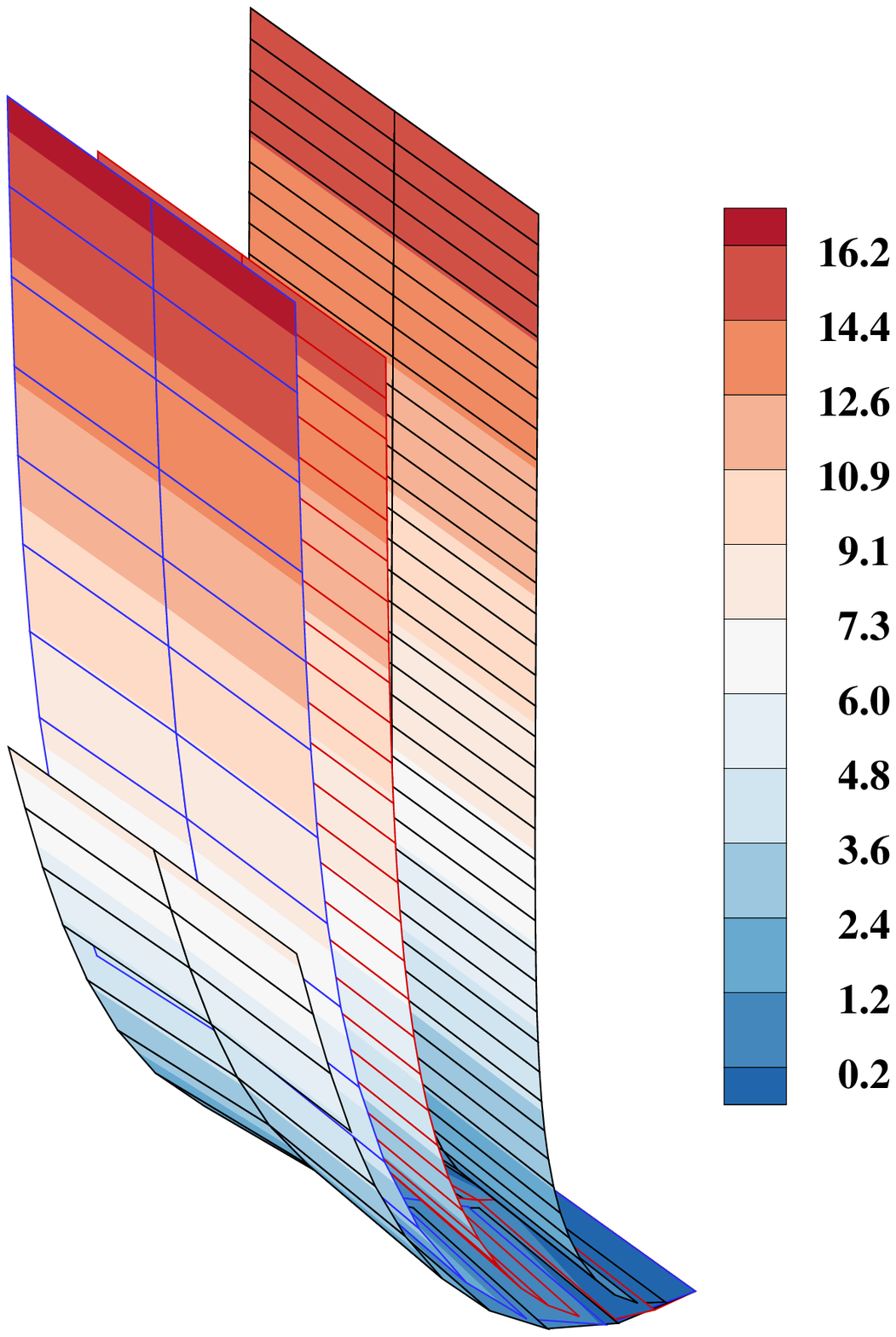} 
\put(-65,120){\scriptsize(\textbf{e})}
\end{minipage}
\vspace{5pt}
\caption{Beam-like strips under magnetic loading,
(\textbf{a}): the curves of $\frac{u_3^{\text{T}}}{L}$ vs 
$\frac{10^3}{\mu \mu_0}|\mbbB^{\text{ext}}||\tilde{\mbbB}^{\text{rem}}|$, 
(\textbf{b},\textbf{c},\textbf{d},\textbf{e}): sequences of deformation (with the contours of $u_3$ in mm) for $|\mbbB^\text{ext}| \in \{2,5, 15, 50 \}$ (mT)}
\label{Ex1Fig1} 
\end{figure*}

It is recalled from Eq.~\eqref{bfpsi} that the present formulation employs the through-the-thickness stretching parameter $\phi$, which leads to linear shear strain as well as linear normal strain in the thickness direction. 
To show the effect of this parameter, two new cases are considered. In the first case, the condition $\phi=0$ is enforced in the formulation, while the $3$D constitutive equations are still employed. 
In the second case, the plane stress assumption $P_{33}=0$ is enforced and $\phi$ has not been considered in the formulation. Then the constitutive equation is modified to include the plane stress assumption. 
For the thick beam with $AR=10$ and the thin one with $AR=40.5$ the results are displayed in Figs.~\ref{Ex1Fig2}(a,b).
It is noted that for the beams with $AR=17.5$ and $AR=20.5$, similar results are obtained, which have not been shown in the figures.
It is observed the new cases exhibit locking phenomenon in the resulting elements. 
The second new case is better than the first one, however, the improvement is negligible.
In other words, including the through-the-thickness stretching $\phi$ in the formulation and employing the $3$D constitutive equations is an effective method for improving the performance of the present micropolar shell element.

\begin{figure*}
\centering
\begin{minipage}[b]{.48\textwidth}
\includegraphics[width=80mm]{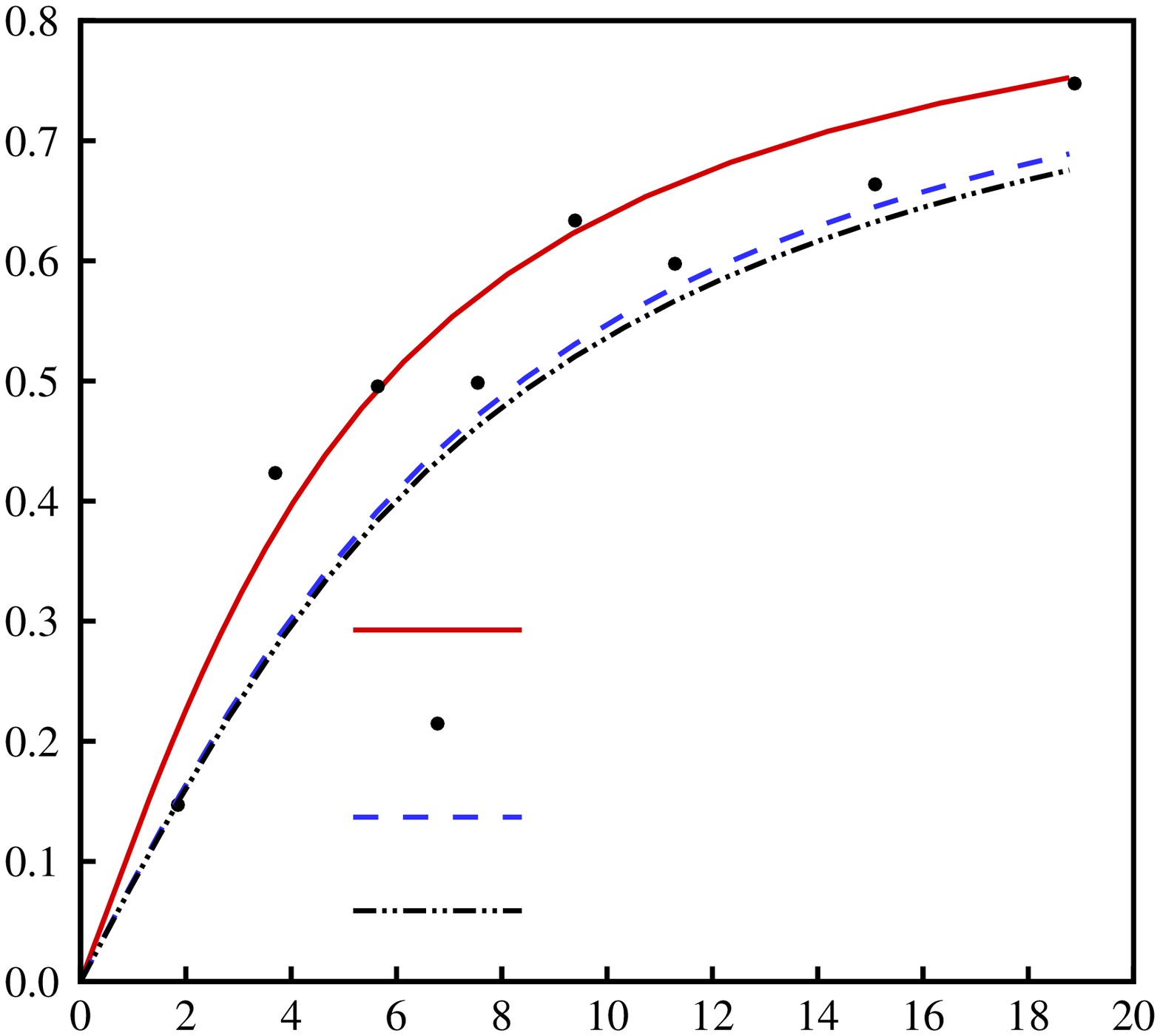} 
\put(-50,120){\scriptsize(\textbf{a})}
\put(-150,4){\scriptsize $\frac{1}{\mu \mu_0} |\mbbB^{\text{ext}}||\tilde{\mbbB}^{\text{rem}}| \times 10^3$}
\put(-120,77){\scriptsize $\phi \neq 0$, 3D constitutive eq.}
\put(-120,62){\scriptsize Experiment, Ref.~\cite{Zhao2019}}
\put(-120,47){\scriptsize $\phi = 0$, 3D constitutive eq.}
\put(-120,32){\scriptsize $\phi = 0$, plane stress}
\put(-225,90){\rotatebox{90}{\footnotesize $u_3^T/L$}}
\end{minipage} \quad
\begin{minipage}[b]{.48\textwidth}
\includegraphics[width=80mm]{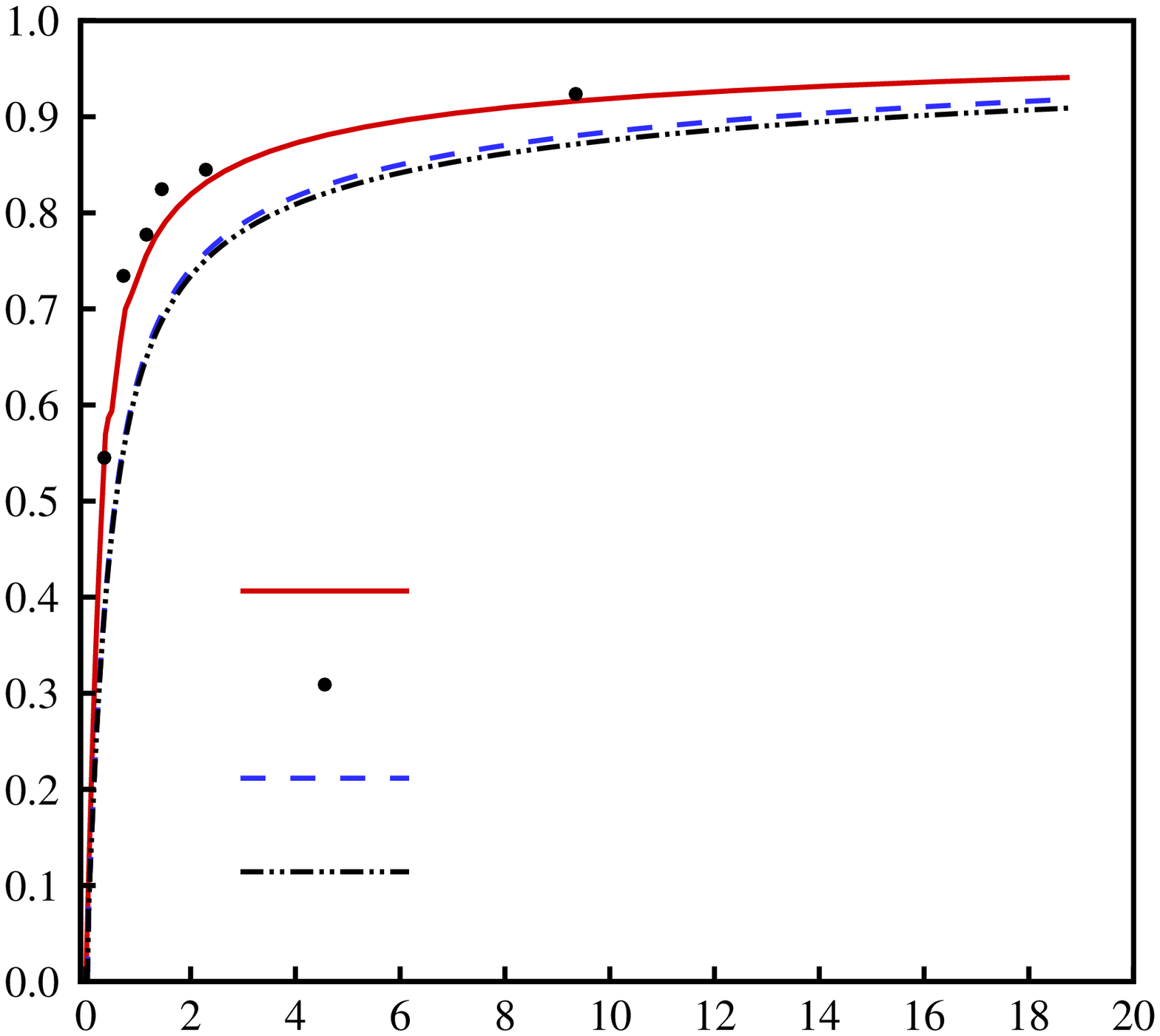} 
\put(-50,120){\scriptsize(\textbf{b})}
\put(-150,4){\scriptsize $\frac{1}{\mu \mu_0} |\mbbB^{\text{ext}}||\tilde{\mbbB}^{\text{rem}}| \times 10^3$}
\put(-140,84){\scriptsize $\phi \neq 0$, 3D constitutive eq.}
\put(-140,69){\scriptsize Experiment, Ref.~\cite{Zhao2019}}
\put(-140,54){\scriptsize $\phi = 0$, 3D constitutive eq.}
\put(-140,39){\scriptsize $\phi = 0$, plane stress}
\put(-225,90){\rotatebox{90}{\footnotesize $u_3^T/L$}}
\end{minipage}
\vspace{-5pt}
\caption{The effect of the through-the-thickness stretching parameter $\phi$ on the load-deflection curve,(\textbf{a}): the thick beam with $AR=10$,  (\textbf{b}): the thin beam with $AR=41$}  
\label{Ex1Fig2} 
\end{figure*}

%
\subsection{Deformation of a hollow cross}
\label{hollow cross}
\vspace{-5pt}
In the present example, the finite deformation of a hollow cross under magnetic loading is simulated.
Following Kim et al.~\cite{Kim2018} and Zhao et al.~\cite{Zhao2019}, the geometry of the hollow cross is composed of $24$ trapezoidal blocks (Fig.~\ref{Ex2Fig1}(\textbf{a})). 
The thickness is $0.41$ mm, and all dimensions in the $X_1X_2$ plane are displayed in the figure.
The values of the mechanical properties are identical to those given in the previous example.
As can be seen in Fig.~\ref{Ex2Fig1}(\textbf{a}), the direction of $\tilde{\mathbbm{B}}^{\text{r}}$ is constant in each block, but it varies in different blocks.
The constant value of $|\tilde{\mathbbm{B}}^{\text{r}}|=102$ (mT) has been considered for the referential remnant magnetic flux \cite{Kim2018,Zhao2019}. 
The maximum external magnetic flux density 
$\mathbbm{B}^{\text{ext}}_{\text{max}}= -200 \mathbbm{e}_3$ (mT) 
acts on the body.
The symmetry of the geometry allows us to discretize merely $1/4$ of the body in the $X_1X_2$ plane.
Moreover, the displacement component $u_3$ at the points $A$ and $G$ is assumed to be zero.

By performing various numerical simulations, it is found that a  $6 \times 6$ mesh of shell elements in each trapezoidal block leads to convergent results.
Fig.~\ref{Ex2Fig1}(\textbf{a}) displays the displacement component $u_3$ against the normalized loading parameter 
$\frac{10^3}{\mu \mu_0}|\mbbB^{\text{ext}}||\tilde{\mbbB}^{\text{rem}}|$
at some material points.
At the final stage of deformation, the lateral displacement at the points $E$ and $C$ is very close to each other.
More precisely, the maximum $u_3$ displacement achieved at the point $C$ is $10.39$ mm.
Fig.~\ref{Ex2Fig1}(\textbf{b}) depicts the fully deformed hollow cross observed in the experiments of Kim et al.~\cite{Kim2018}.
Moreover, the deformed shapes of the hollow cross under four different values of the external magnetic flux are displayed in Figs.~\ref{Ex2Fig1}(\textbf{c,d,e,f}).
By comparing figures \ref{Ex2Fig1}(\textbf{b}) and \ref{Ex2Fig1}(\textbf{e}), it is deduced that the final deformed shape obtained by the present formulation is qualitatively similar to that reported in the experimental studies of Kim et al.~\cite{Kim2018}.
\begin{figure*}
\centering
\begin{minipage}[b]{.5\textwidth}
\includegraphics[width=85mm]{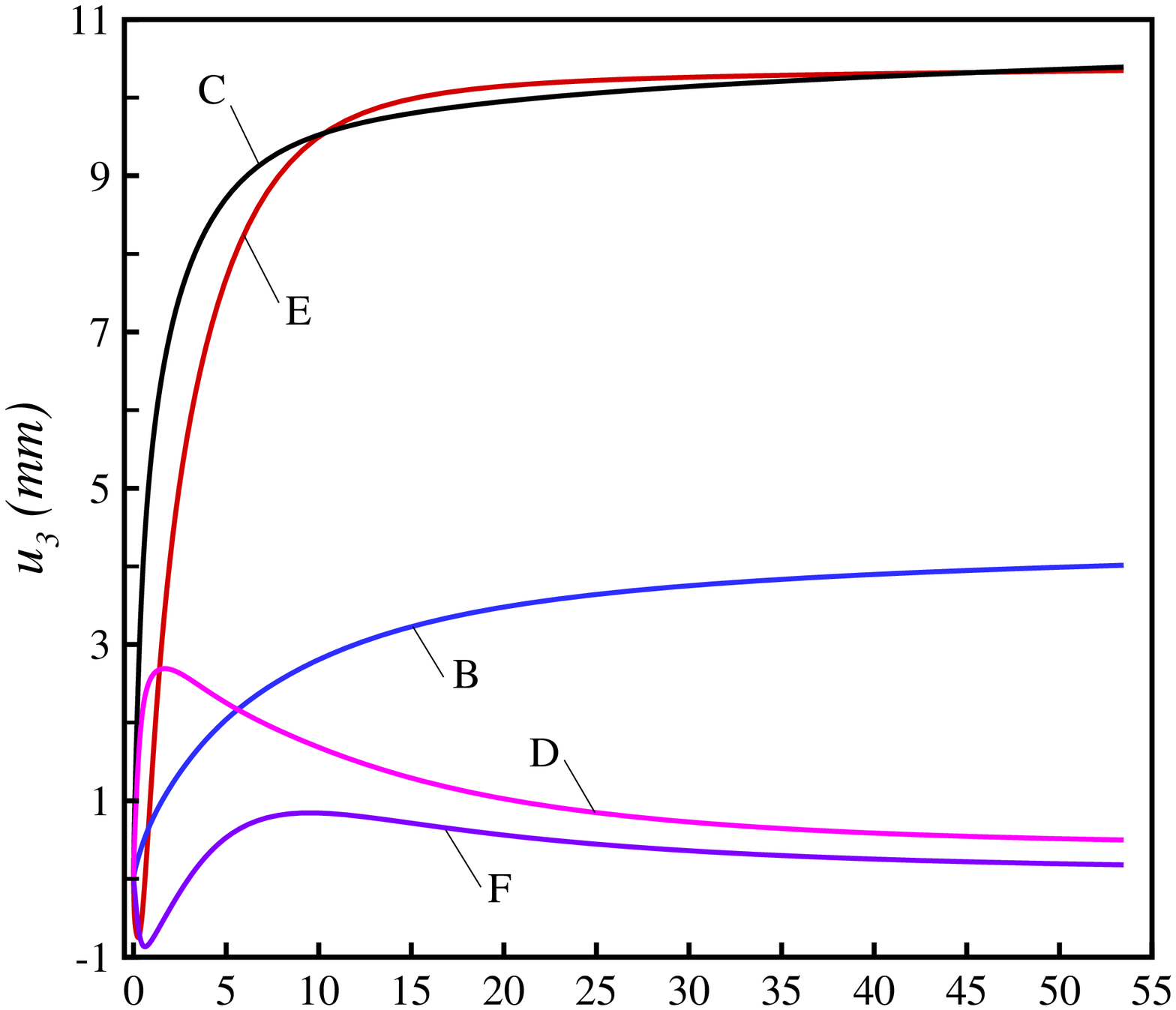} 
\put(-205,94){\includegraphics[width=70mm]{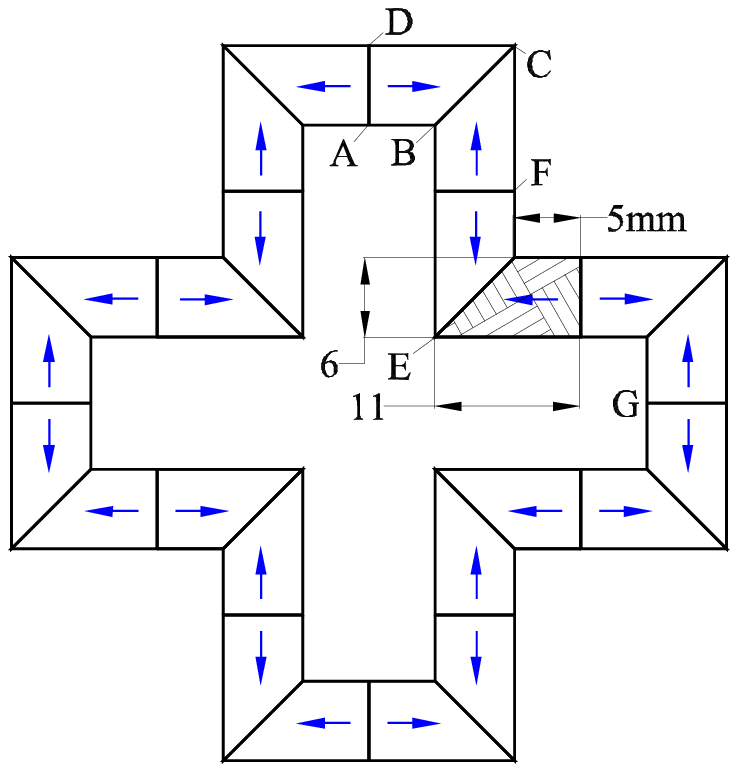}}
\put(-50,140){\scriptsize(\textbf{a})}
\put(-150,4){\scriptsize $\frac{1}{\mu \mu_0} |\mbbB^{\text{ext}}||\tilde{\mbbB}^{\text{rem}}| \times 10^3$}
\end{minipage}
\begin{minipage}[b]{.48\textwidth}
\includegraphics[trim=-2cm -1cm  0cm  0cm, width=50mm]{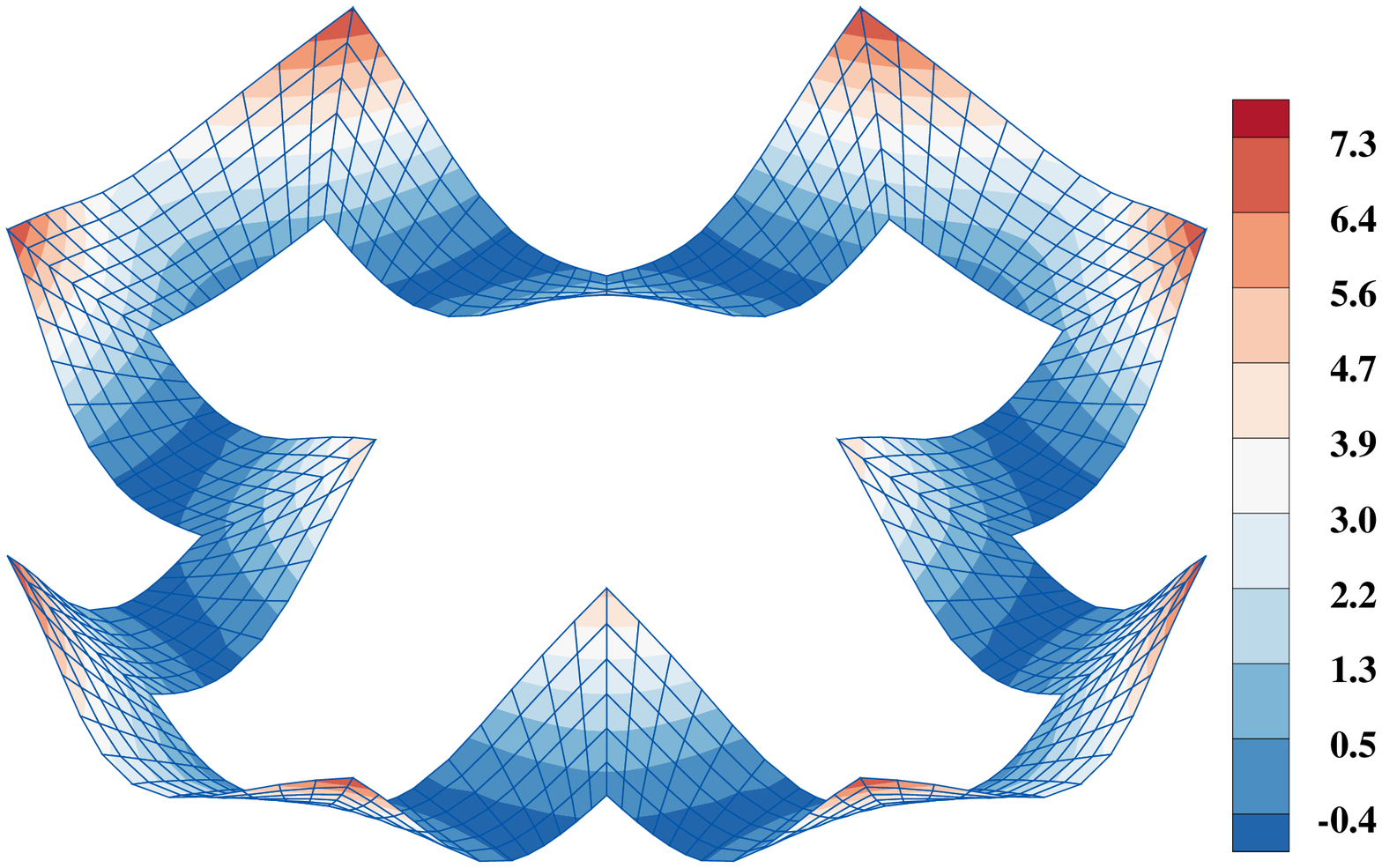} 
\put(-110,120){\includegraphics[width=25mm]{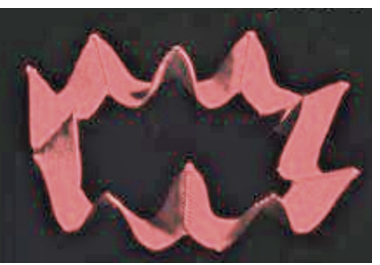}}
\put(-30,140){\scriptsize(\textbf{b})}
\put(-30,110){\scriptsize(\textbf{c})}
\end{minipage} \\
\vspace{10pt}
\begin{minipage}[b]{.32\textwidth}
\includegraphics[trim=0cm 1cm  0cm  0cm, width=45mm]{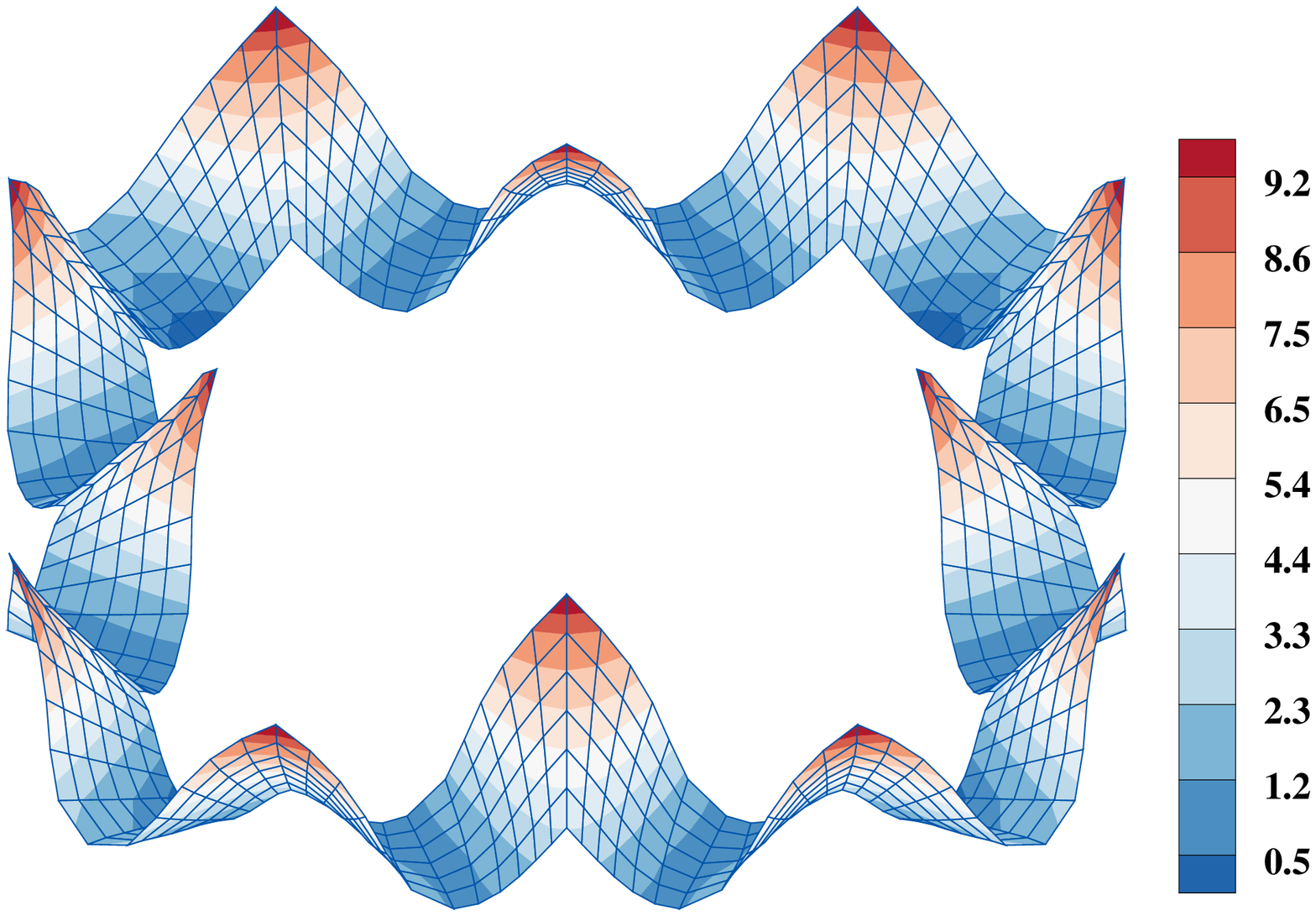} 
\put(-85,100){\scriptsize(\textbf{d})}
\end{minipage}   
\begin{minipage}[b]{.32\textwidth}
\includegraphics[trim=0cm 3mm  0cm  0cm, width=45mm]{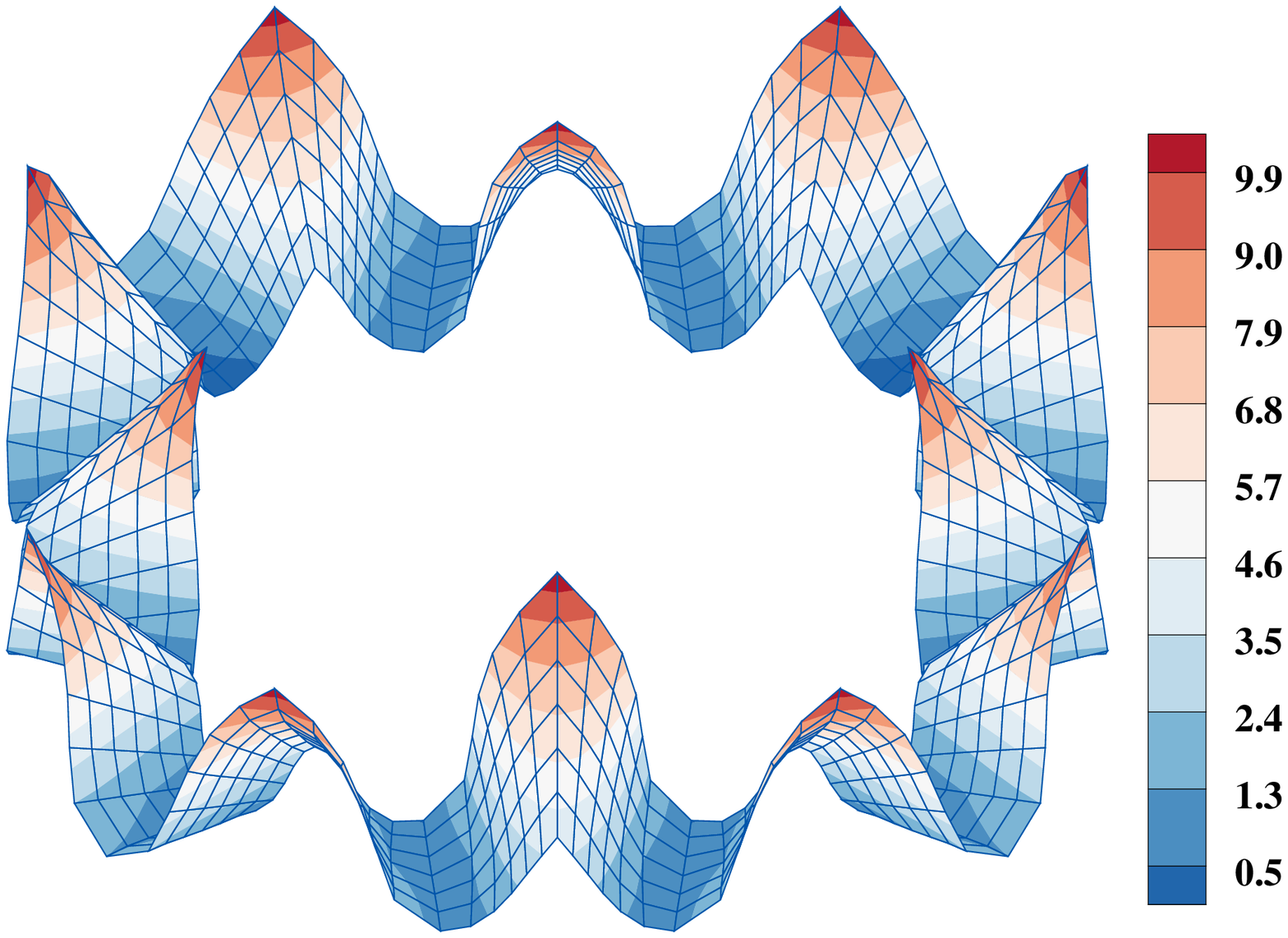} 
\put(-85,100){\scriptsize(\textbf{e})}
\end{minipage}
\begin{minipage}[b]{.32\textwidth}
\includegraphics[trim=0cm 0cm  0cm  0cm, width=45mm]{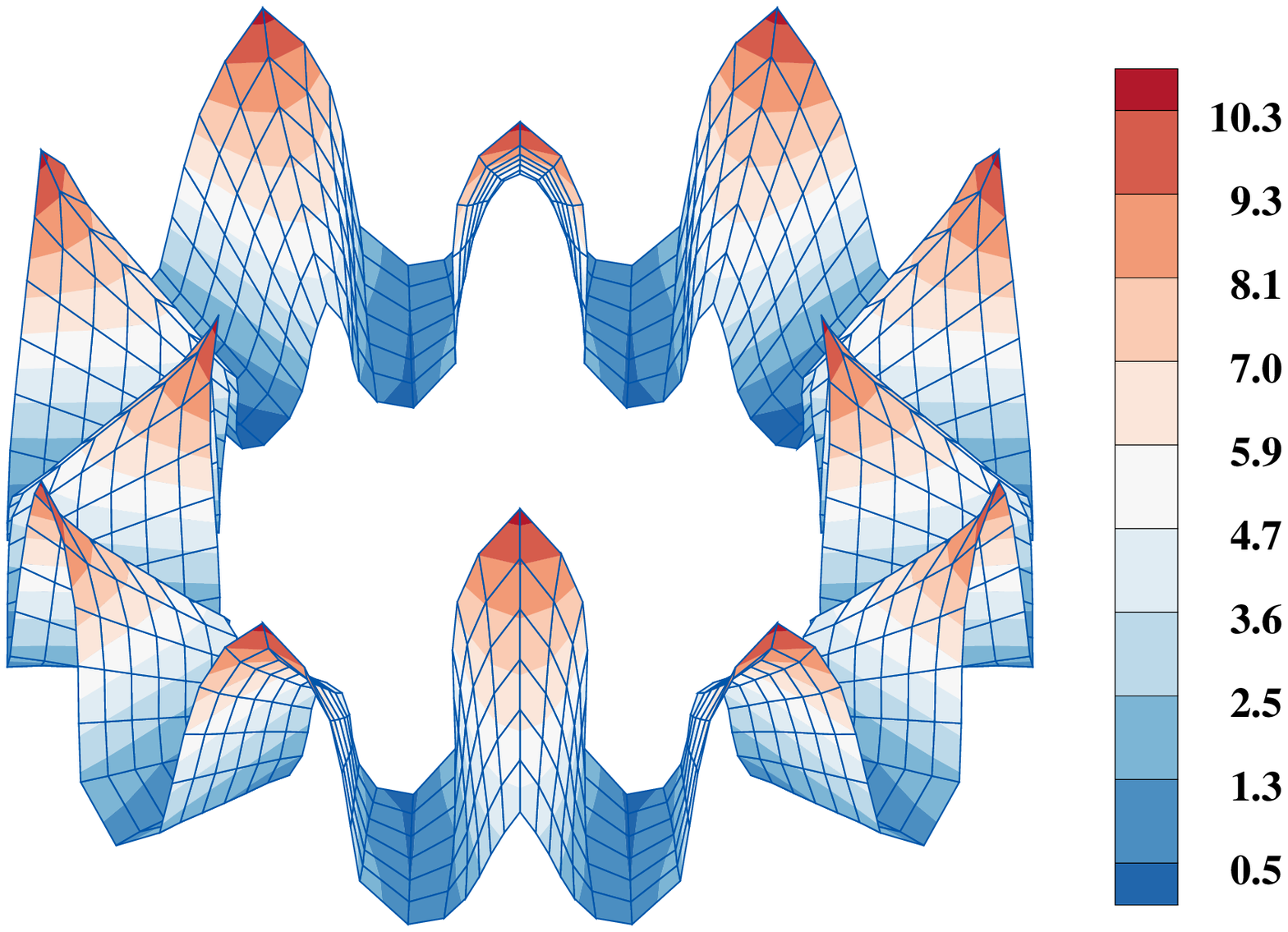} 
\put(-85,100){\scriptsize(\textbf{f})}
\end{minipage}
\vspace{-5pt}
\caption{A hollow cross under magnetic loading, 
(\textbf{a}): load-displacement curves, 
(\textbf{b}): experiment \cite{Kim2018}, 
(\textbf{c},\textbf{d},\textbf{e},\textbf{f}): sequences of deformation (with the contours of $u_3$ in mm) for $|\mbbB^\text{ext}| \in \{10,50, 100, 200 \}$ (mT)}
\label{Ex2Fig1} 
\end{figure*}
\subsection{Deformation analysis of a thin cross}
\label{Cross}
\vspace{-5pt}
The finite elastic response of a thin cross made of HMSMs is simulated in this example.
The geometry of the cross involves nine welded $6 \times 6$ (mm) square-shaped blocks (Fig.~\ref{Ex3Fig1}(\textbf{a})), and the thickness is $0.9$ mm.
%
%
%
The magnitude of $\tilde{\mathbbm{B}}^{\text{r}}$ is the constant value of $94$ mT.
%
%
To deform the body by magnetic loading, the maximum value of $\mathbbm{B}^{\text{ext}}$ is considered to be $40$ mT, which is applied perpendicular to the plane of the cross in the $X_3$ direction.
Moreover, the mechanical properties are $\mu=135$ and $\lambda=3250$ (kPa).
Due to symmetry in the $X_1X_2$ plane, only $1/4$ of the cross is used in the simulations.

From numerical experiments, it is found that a mesh containing $15$ shell elements along $AC$ and $3$ elements along $AA'$ leads to convergence in the results.
The displacement component $u_3$ at some points is plotted in Fig.~\ref{Ex3Fig1}(\textbf{a}).
As usual, the horizontal axis is considered to be the nondimensional loading 
$\frac{10^3}{\mu \mu_0}|\mbbB^{\text{ext}}||\tilde{\mbbB}^{\text{rem}}|$.
It is noted that the $u_3$ displacement of the point $D$ has been considered to be zero, and the maximum lateral displacement $22.78$ mm is predicted at the point $A$.
The fully deformed shape of the cross observed experimentally in Ref.~\cite{Kuang2021} is displayed in Fig.~\ref{Ex3Fig1}(\textbf{b}).
Moreover, the deformed shapes of the cross under four different values of the external magnetic flux are illustrated in Figs.~\ref{Ex3Fig1}(\textbf{c,d,e,f}).
Obviously, the final deformed shape in Fig.~\ref{Ex3Fig1}(\textbf{e}), predicted by the present formulation, is qualitatively similar to that observed in the experiments of Kuang et al.~\cite{Kuang2021} in Fig.~\ref{Ex3Fig1}(\textbf{b}).
\begin{figure*}
\centering
\begin{minipage}[b]{.52\textwidth}
\includegraphics[width=80mm]{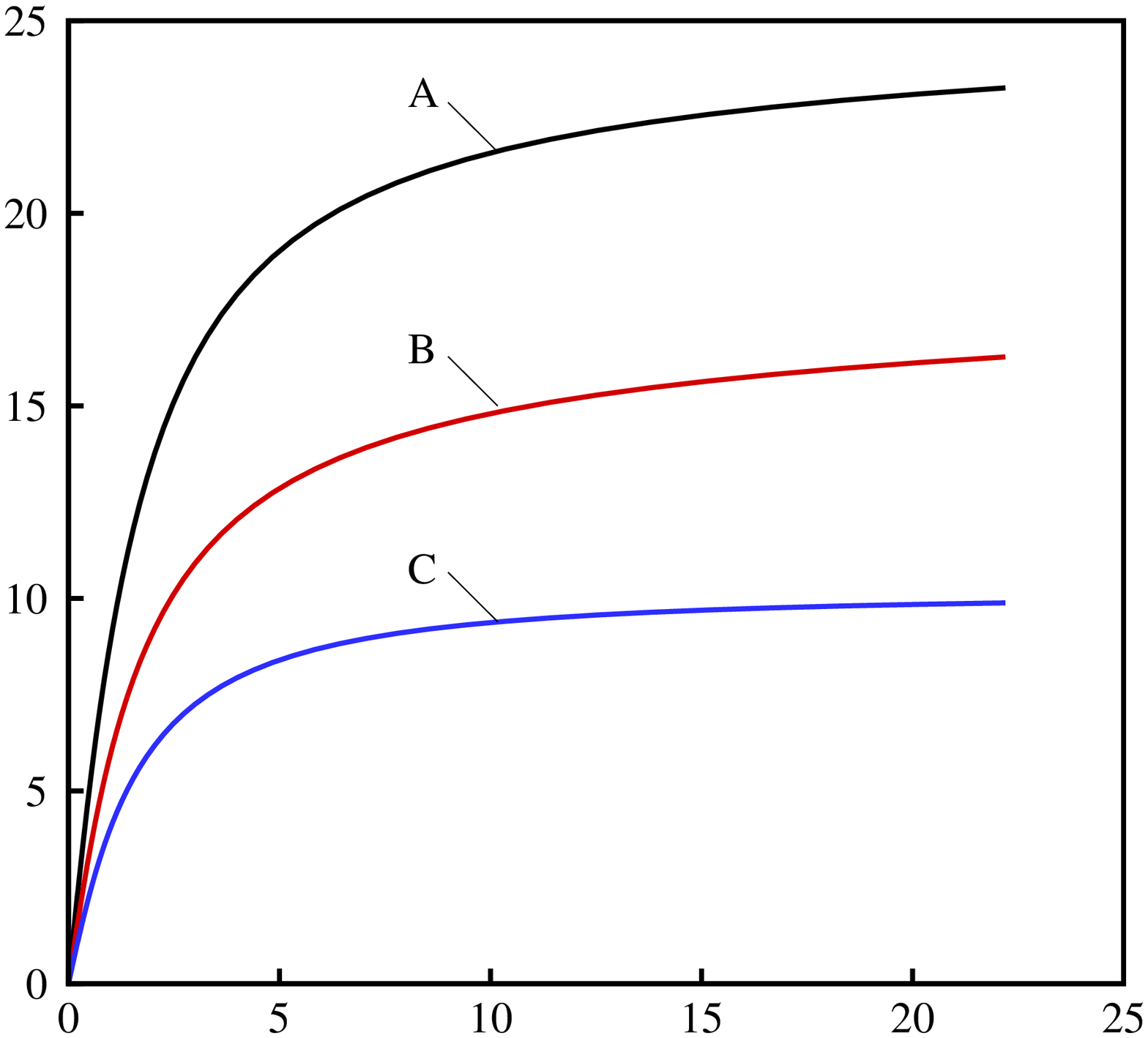} 
\put(-165,24){\includegraphics[width=45mm]{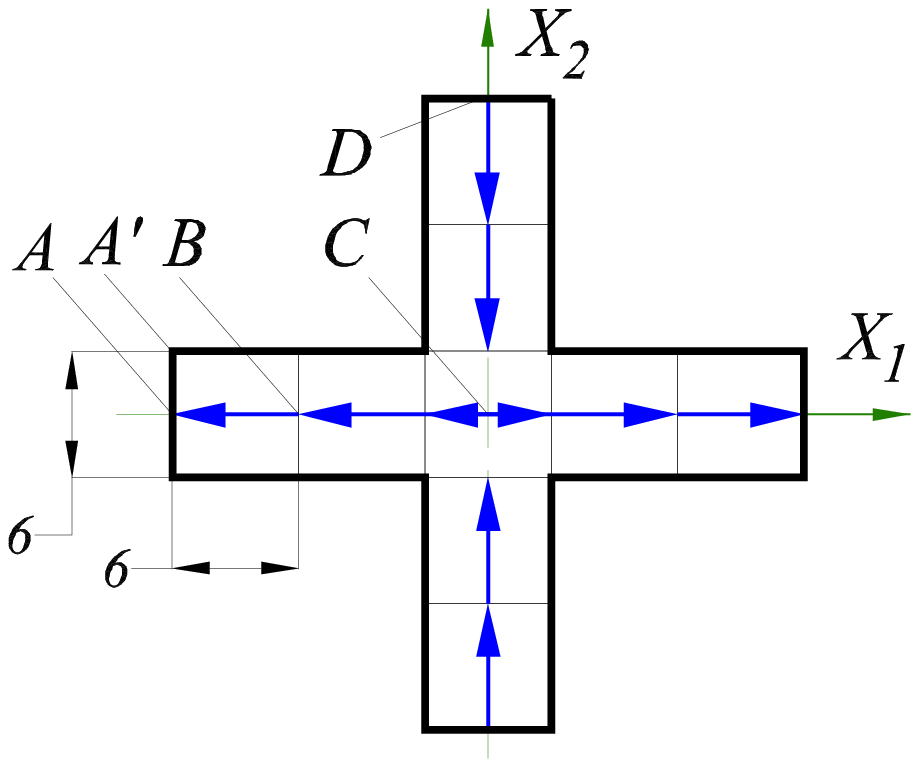}}
\put(-50,140){\scriptsize(\textbf{a})}
\put(-150,4){\scriptsize $\frac{1}{\mu \mu_0} |\mbbB^{\text{ext}}||\tilde{\mbbB}^{\text{rem}}| \times 10^3$}
\end{minipage}
\begin{minipage}[b]{.45\textwidth}
\includegraphics[trim=0cm -1cm  0cm  0cm, width=50mm]{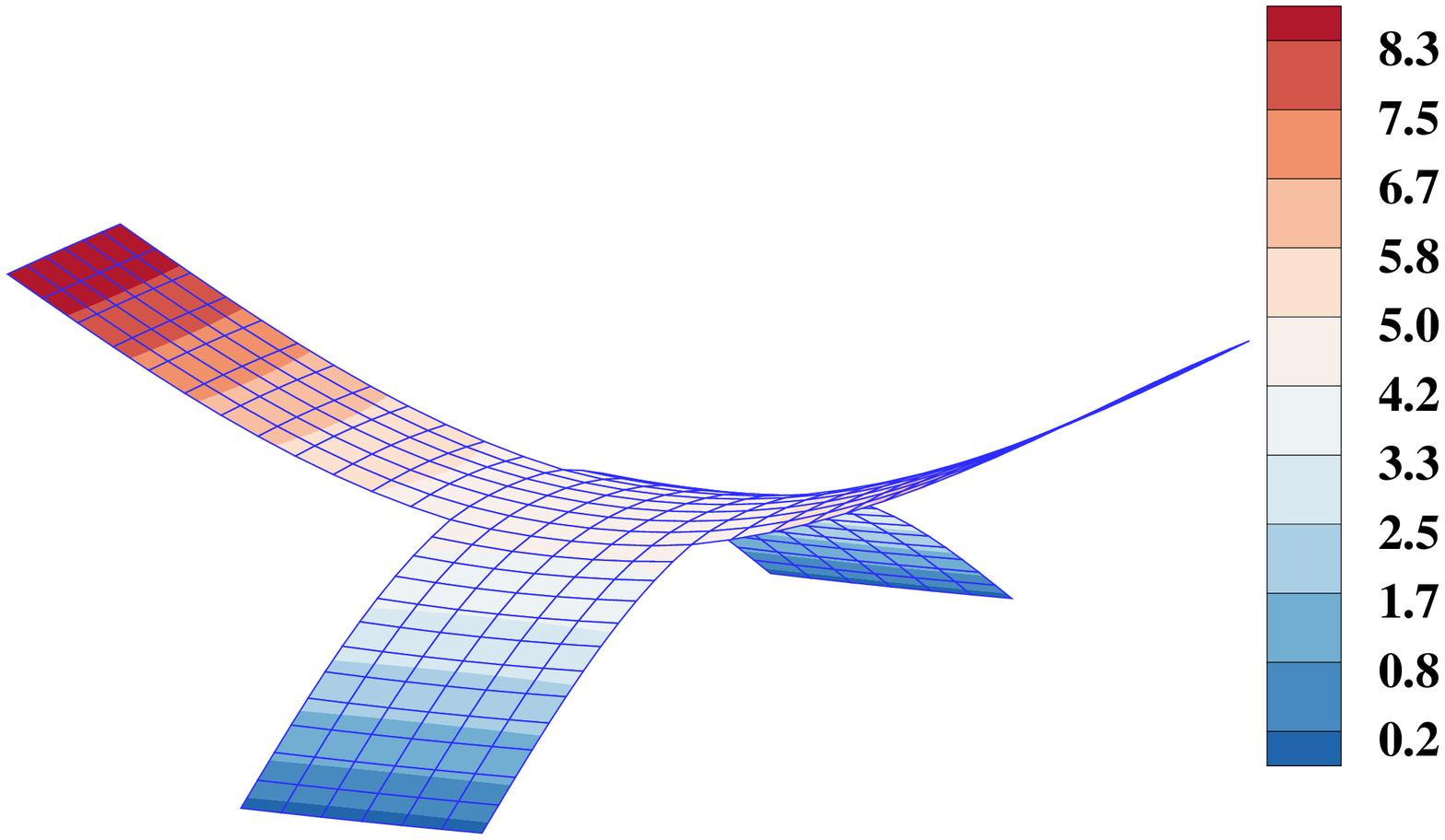} 
\put(-100,90){\includegraphics[width=2cm]{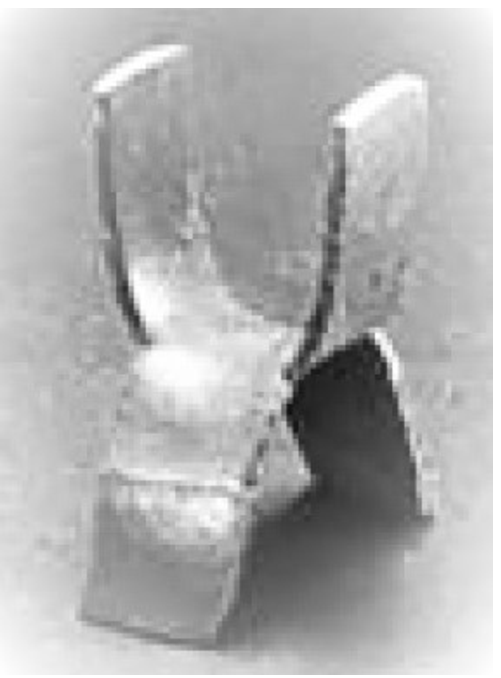}}
\put(-120,140){\scriptsize(\textbf{b})}
\put(-120,50){\scriptsize(\textbf{c})}
\end{minipage} \\
\vspace{10pt}
\begin{minipage}[b]{.32\textwidth}
\includegraphics[trim=0cm 1cm  0cm  0cm, width=40mm]{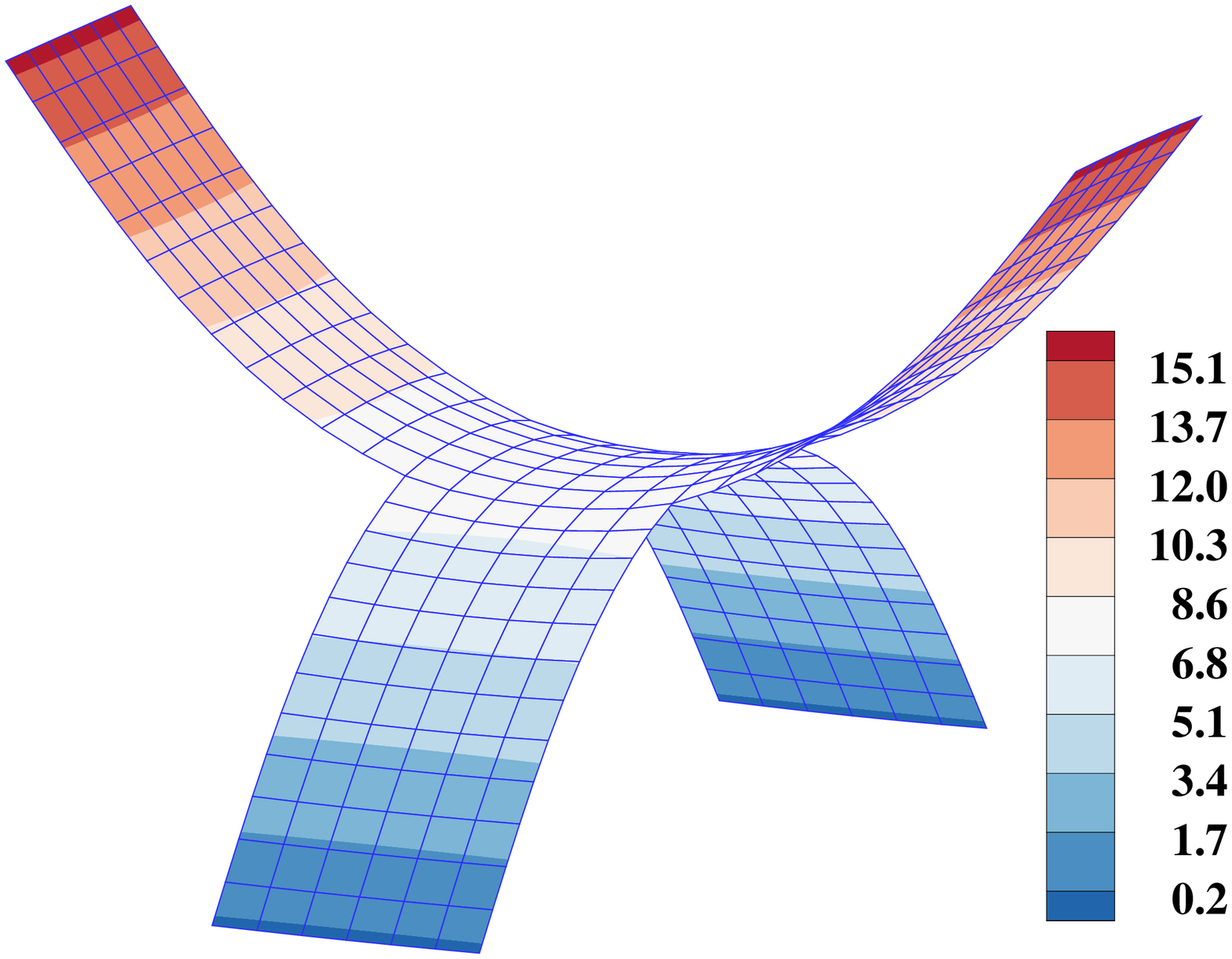} 
\put(-115,50){\scriptsize(\textbf{d})}
\end{minipage}   
\begin{minipage}[b]{.32\textwidth}
\includegraphics[trim=0cm 3mm  0cm  0cm, width=40mm]{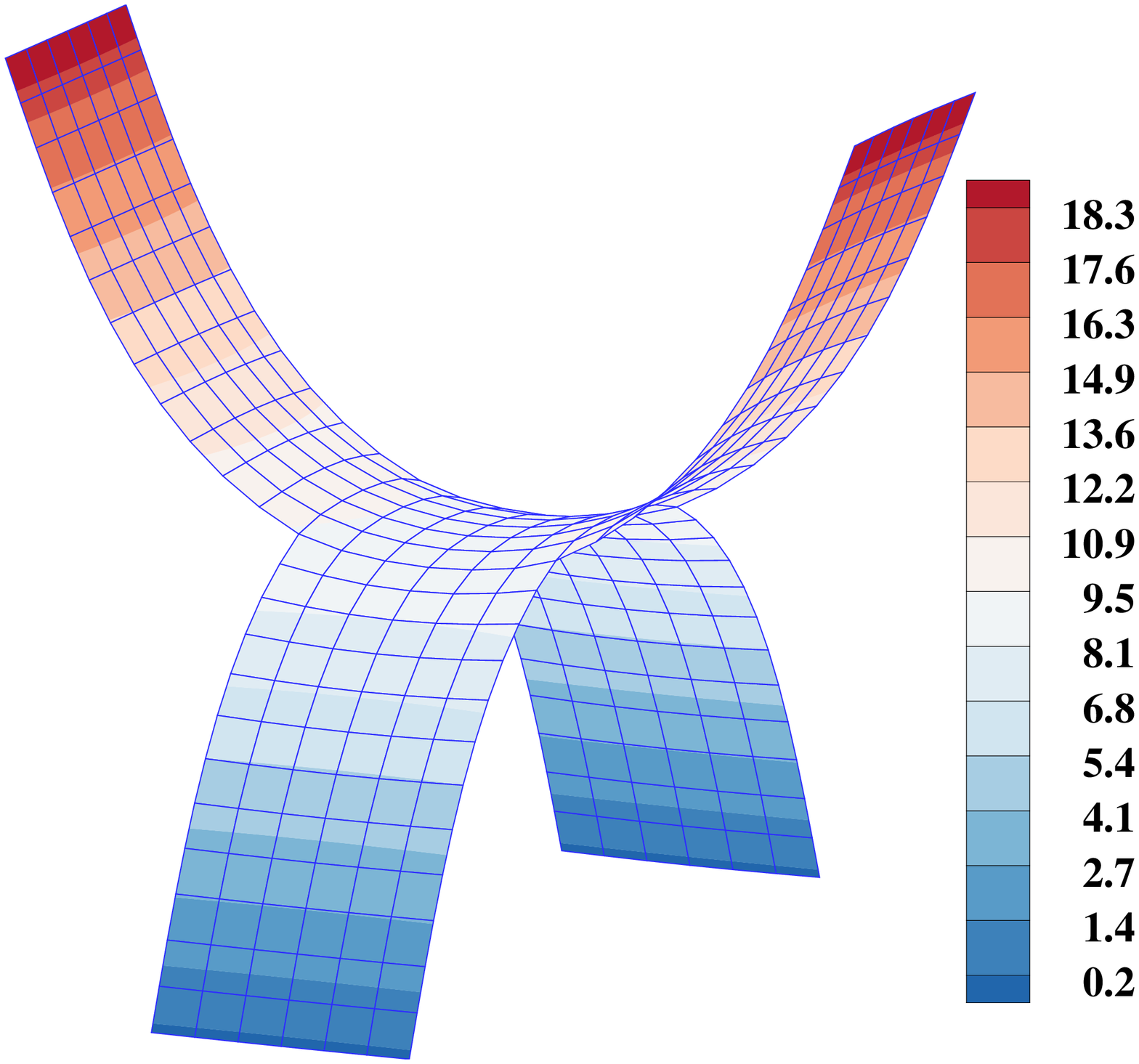} 
\put(-120,50){\scriptsize(\textbf{e})}
\end{minipage}
\begin{minipage}[b]{.32\textwidth}
\includegraphics[trim=0cm 0cm  0cm  0cm, width=45mm]{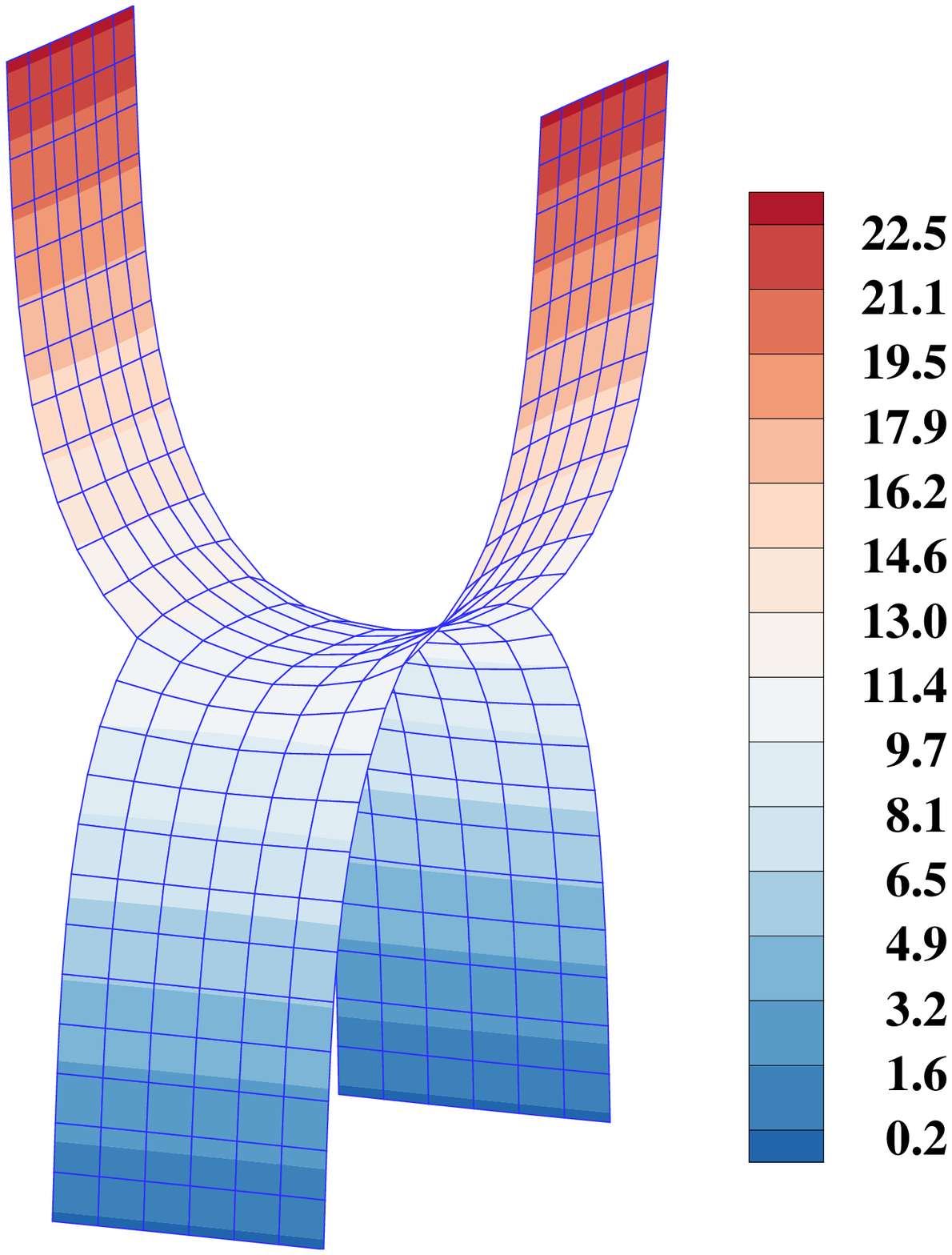} 
\put(-110,50){\scriptsize(\textbf{f})}
\end{minipage}
\vspace{0pt}
\caption{Deformation of a thin cross,
(\textbf{a}): load-displacement curves, 
(\textbf{b}): experiment \cite{Kuang2021}, 
(\textbf{c},\textbf{d},\textbf{e},\textbf{f}): sequences of deformation (with the contours of $u_3$ in mm) for $|\mbbB^\text{ext}| \in \{2,5, 10, 40 \}$ (mT)}
\label{Ex3Fig1} 
\end{figure*}
\subsection{Magnetostrictive response of an H-shaped structure}
\label{H-shaped}
\vspace{-5pt}
In this example, the mechanical response of an H-shaped thin structure to magnetic stimuli is simulated. 
The geometry of the structure, consisting of fifteen blocks, is displayed in Fig.~\ref{Ex4Fig1}(\textbf{a}).
The dimensions, mechanical, and magnetic properties of the blocks are identical to those given in the example \ref{Cross}.
The maximum applied magnetic loading is  
$\mathbbm{B}^{\text{ext}}_{\text{max}}= -50 \mathbbm{e}_3$ (mT).
By considering the symmetry properties of the geometry, merely $1/4$ of the geometry is analyzed.

Numerical experiments indicate that a  $4 \times 4$ mesh of shell elements in each block yields converging results. 
In other words, the number of elements along $AA'$, $AC$ and $CD$ is $2$, $14$, and $10$, respectively.
Fig.~\ref{Ex4Fig1}(\textbf{a}) demonstrates the variations of the displacement component $u_3$ at some points against the normalized loading 
$\frac{10^3}{\mu \mu_0}|\mbbB^{\text{ext}}||\tilde{\mbbB}^{\text{rem}}|$.
By assuming zero $u_3$ displacement at the point $D$, the maximum value of $u_3=24.65$ mm at the point $A$ is achieved.
The fully deformed shape of the H-shaped structure from the experimental observations of Kuang et al.~\cite{Kuang2021} is illustrated in Fig.~\ref{Ex4Fig1}(\textbf{b}).
Moreover, the deformed shapes of the body under four different values of the external magnetic flux are displayed in Figs.~\ref{Ex4Fig1}(\textbf{c,d,e,f}).
A comparison of figures~\ref{Ex4Fig1}(\textbf{b}) and \ref{Ex4Fig1}(\textbf{e}) shows that the deformed structure obtained by the present shell formulation is qualitatively similar to that reported in the experiments of   Ref.~\cite{Kuang2021}.
\begin{figure*}
\centering
\begin{minipage}[b]{.5\textwidth}
\includegraphics[width=85mm]{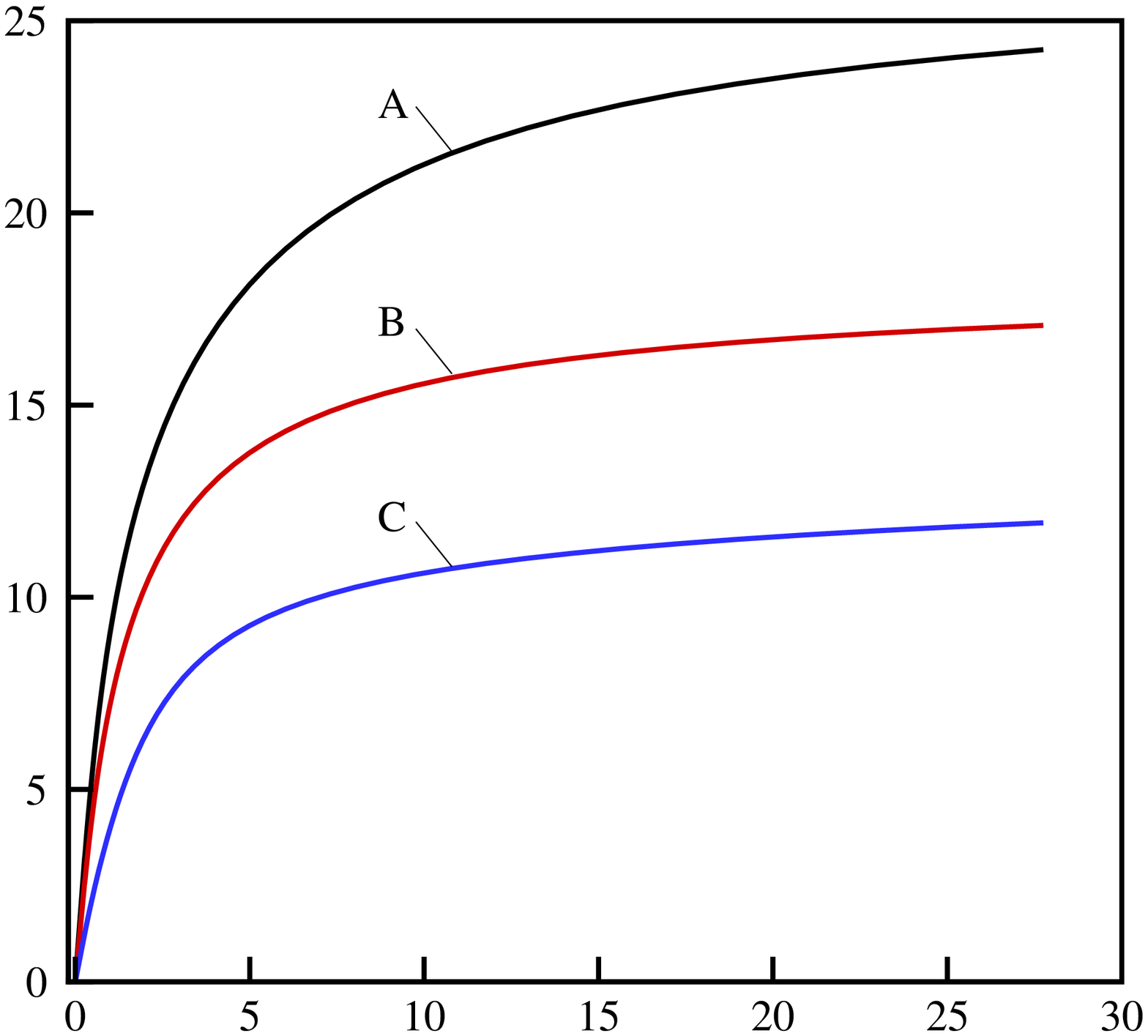} 
\put(-160,30){\includegraphics[width=50mm]{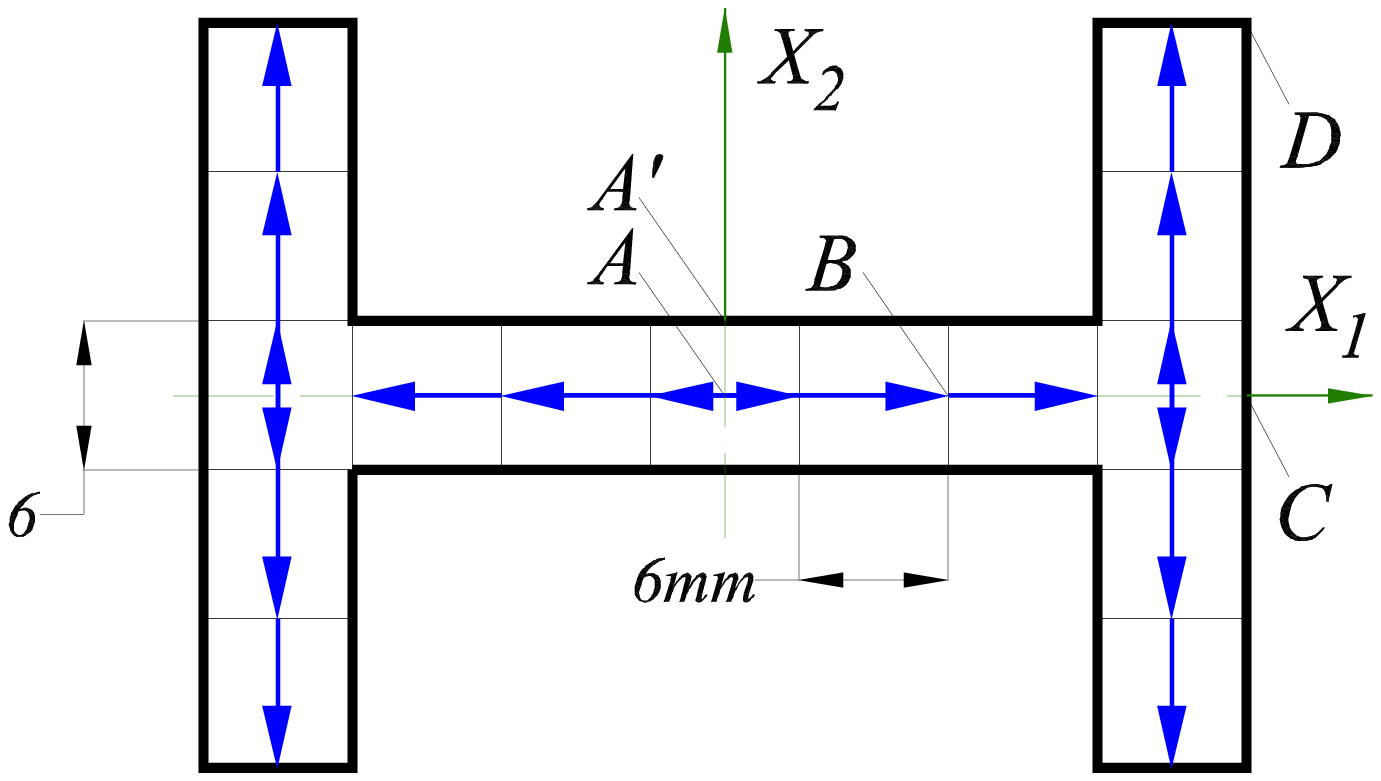}}
\put(-50,150){\scriptsize(\textbf{a})}
\put(-150,4){\scriptsize $\frac{1}{\mu \mu_0} |\mbbB^{\text{ext}}||\tilde{\mbbB}^{\text{rem}}| \times 10^3$}
\end{minipage}
\begin{minipage}[b]{.45\textwidth}
\includegraphics[trim=-3cm -4cm  0cm  0cm, width=52mm]{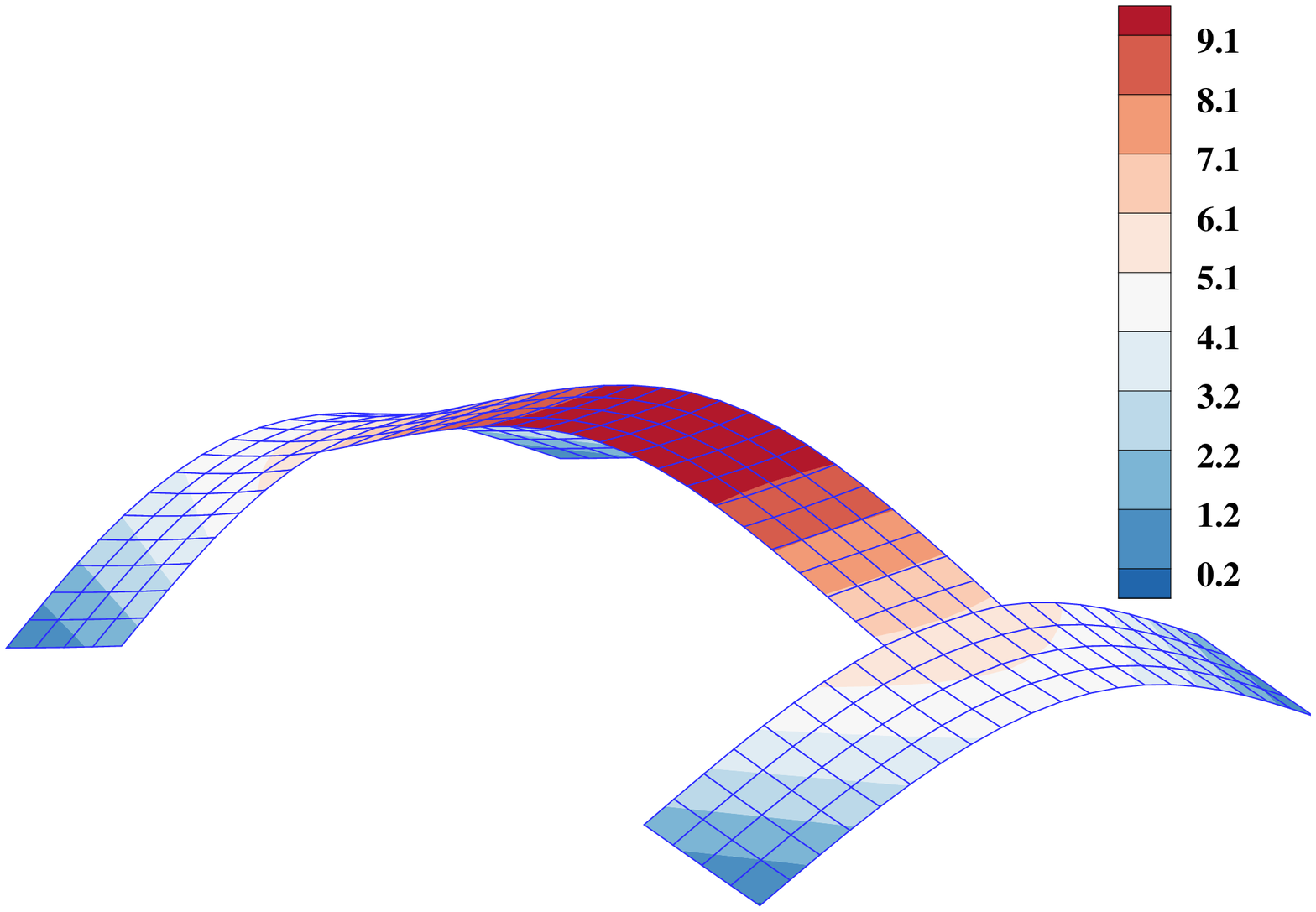} 
\put(-100,120){\includegraphics[width=25mm]{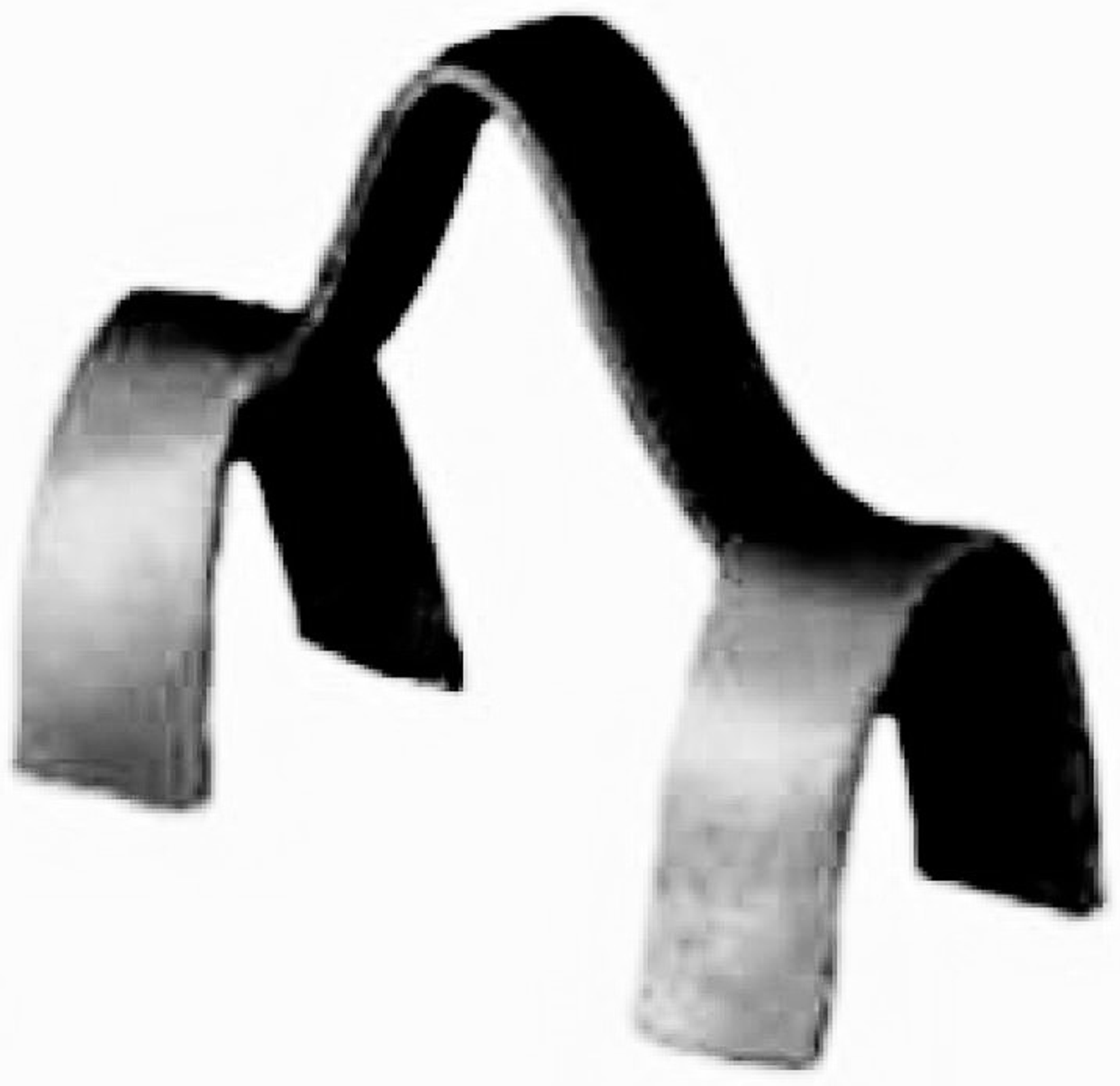}}
\put(-120,130){\scriptsize(\textbf{b})}
\put(-120,75){\scriptsize(\textbf{c})}
\end{minipage} \\
\vspace{10pt}
\begin{minipage}[b]{.32\textwidth}
\includegraphics[trim=0cm 1cm  0cm  0cm, width=47mm]{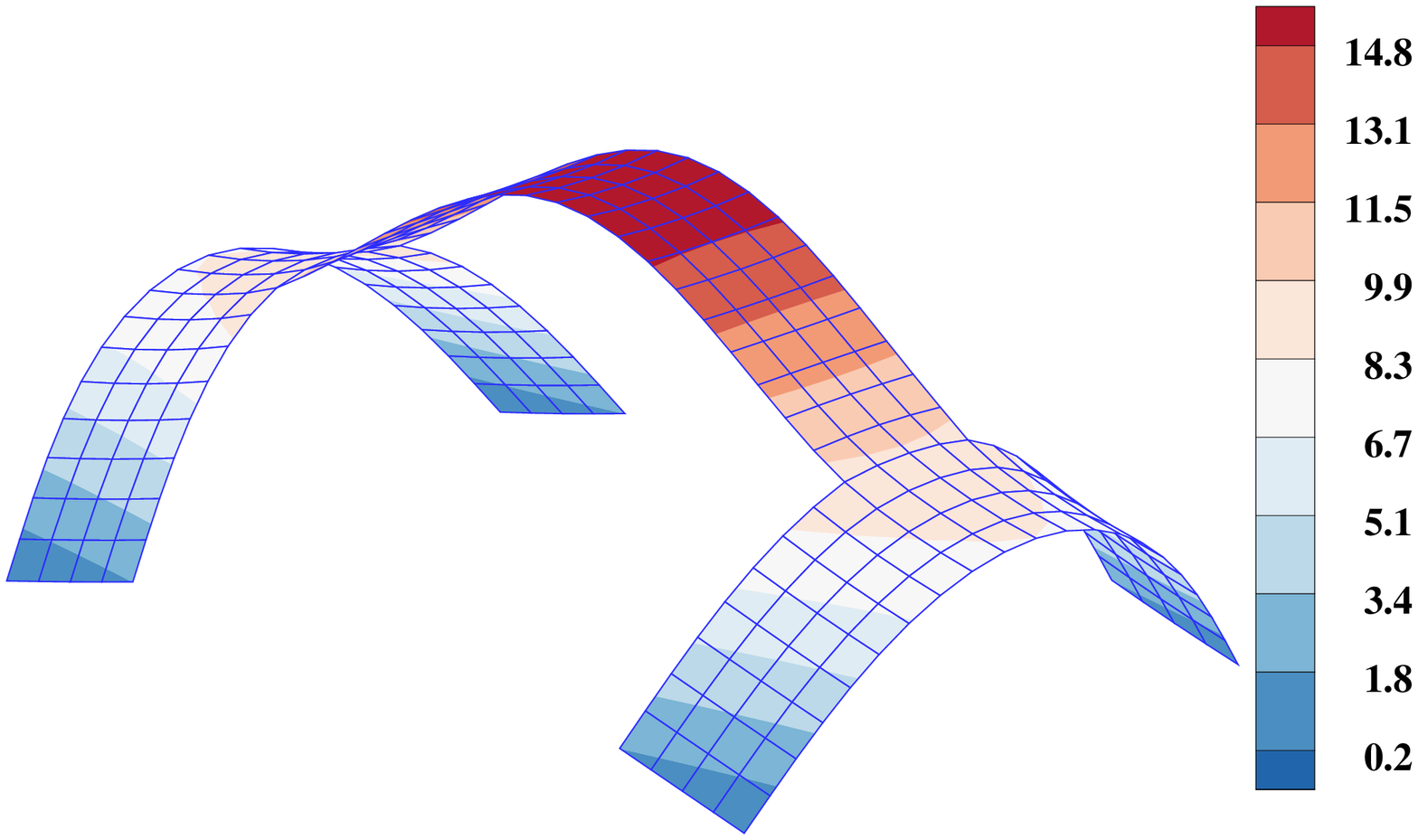} 
\put(-120,90){\scriptsize(\textbf{d})}
\end{minipage}   
\begin{minipage}[b]{.32\textwidth}
\includegraphics[trim=0cm 3mm  0cm  0cm, width=47mm]{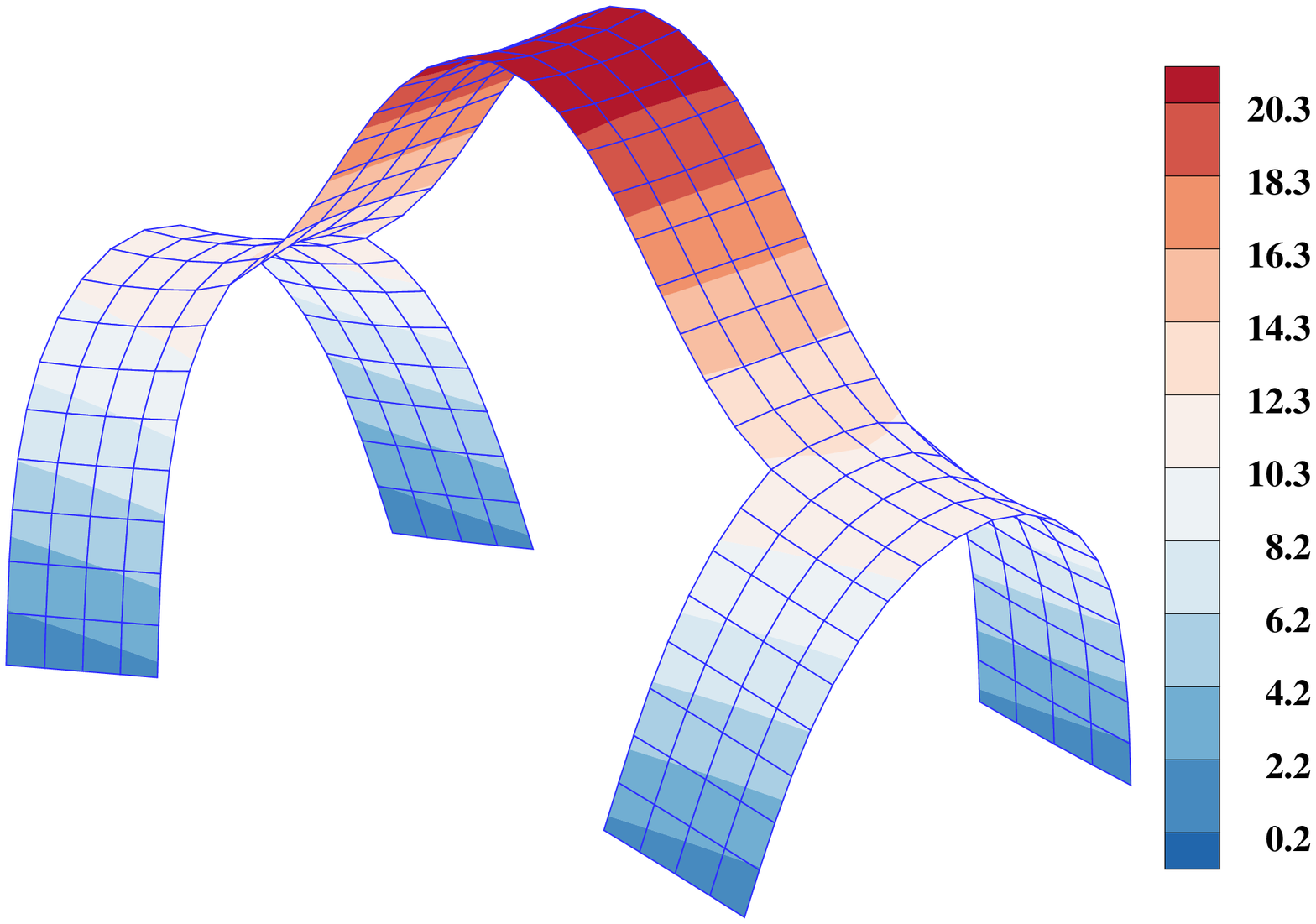} 
\put(-125,90){\scriptsize(\textbf{e})}
\end{minipage}
\begin{minipage}[b]{.32\textwidth}
\includegraphics[trim=0cm 0cm  0cm  0cm, width=47mm]{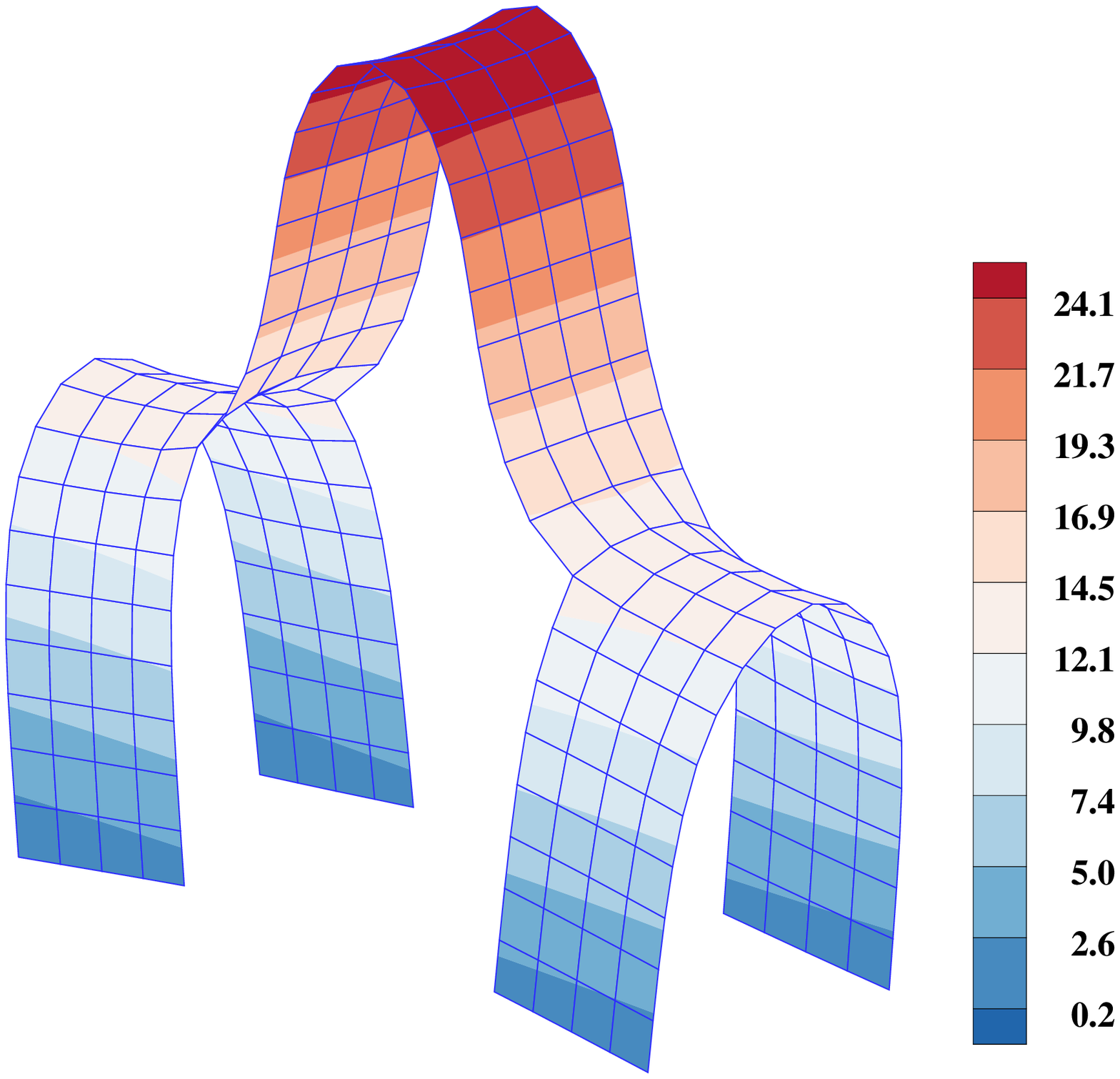} 
\put(-120,90){\scriptsize(\textbf{f})}
\end{minipage}
\vspace{0pt}
\caption{Response of an H-shaped structure,
(\textbf{a}): load-displacement curves, 
(\textbf{b}): experiment \cite{Kuang2021}, 
(\textbf{c},\textbf{d},\textbf{e},\textbf{f}): sequences of deformation (with the contours of $u_3$ in mm) for $|\mbbB^\text{ext}| \in \{2,5, 15, 50 \}$ (mT)}
\label{Ex4Fig1} 
\end{figure*}
\subsection{Deformation of a cylinder (magnetic pump)}
\label{pump}
\vspace{-5pt}
The finite elastic response of a cylindrical shell to magnetic loading is simulated in this example. 
As will be shown below, the deformation pattern in the cylinder is so that it may be used as a macro- or micro-fluidic magnetic pump in practical applications. 
In a relatively similar context, an electro-active polymer-based micro-fluidic pump can be seen in Yan et al. \cite{yan2015}.
In the present case, it is assumed that the cylinder has been made of the same blocks as described in the example \ref{Cross}.
To construct the geometry, 24 blocks in the circumferential direction and 20 ones along the axis of the cylinder are used.
Therefore, the mean radius and length of the cylinder are $R=22.9$ and $L=120$ (mm), respectively.
The remnant magnetic flux $\tilde{\mbbB}^{\text{rem}}$ is assumed to be tangent to the cylinder surface and perpendicular to the $X_2$ axis.
Moreover, it has a positive component along the $X_3$ axis.
The magnetic flux 
$\mathbbm{B}^{\text{ext}}_{\text{max}}= 150 \mathbbm{e}_3$ (mT) 
acts on the cylinder.
Moreover, both ends of the cylinder are considered to be clamped.
Symmetry considerations allow us to simulate $1/4$ of the full geometry.

Numerical simulations show that a mesh of $24 \times 20$ elements is sufficient to obtain convergent results. 
Variations of the displacement components $u_1$ and $u_3$ against $\frac{10^3}{\mu \mu_0}|\mbbB^{\text{ext}}||\tilde{\mbbB}^{\text{rem}}|$ are plotted in Fig.~\ref{Ex5Fig1}(\textbf{a}).
The coordinates of the material points $A$, $B$ and $C$, lying in the $XZ$-plane, are ($0,R$), ($R,0$), and $\frac{1}{\sqrt{2}}(R,R)$, respectively.
The maximum (horizontal) displacement occurs at the point $B$ and is about $16.75$ mm.
The deformed shapes of the cylinder under four different values of the external magnetic flux are demonstrated in Figs.~\ref{Ex4Fig1}(\textbf{b,c,d,e}).
It is observed that under the applied magnetic flux, the cylinder contracts at its middle section. This is the reason why it can be used as a magnetic pump in real applications. 
\begin{figure*}
\centering
\begin{minipage}[b]{.5\textwidth}
\includegraphics[width=85mm]{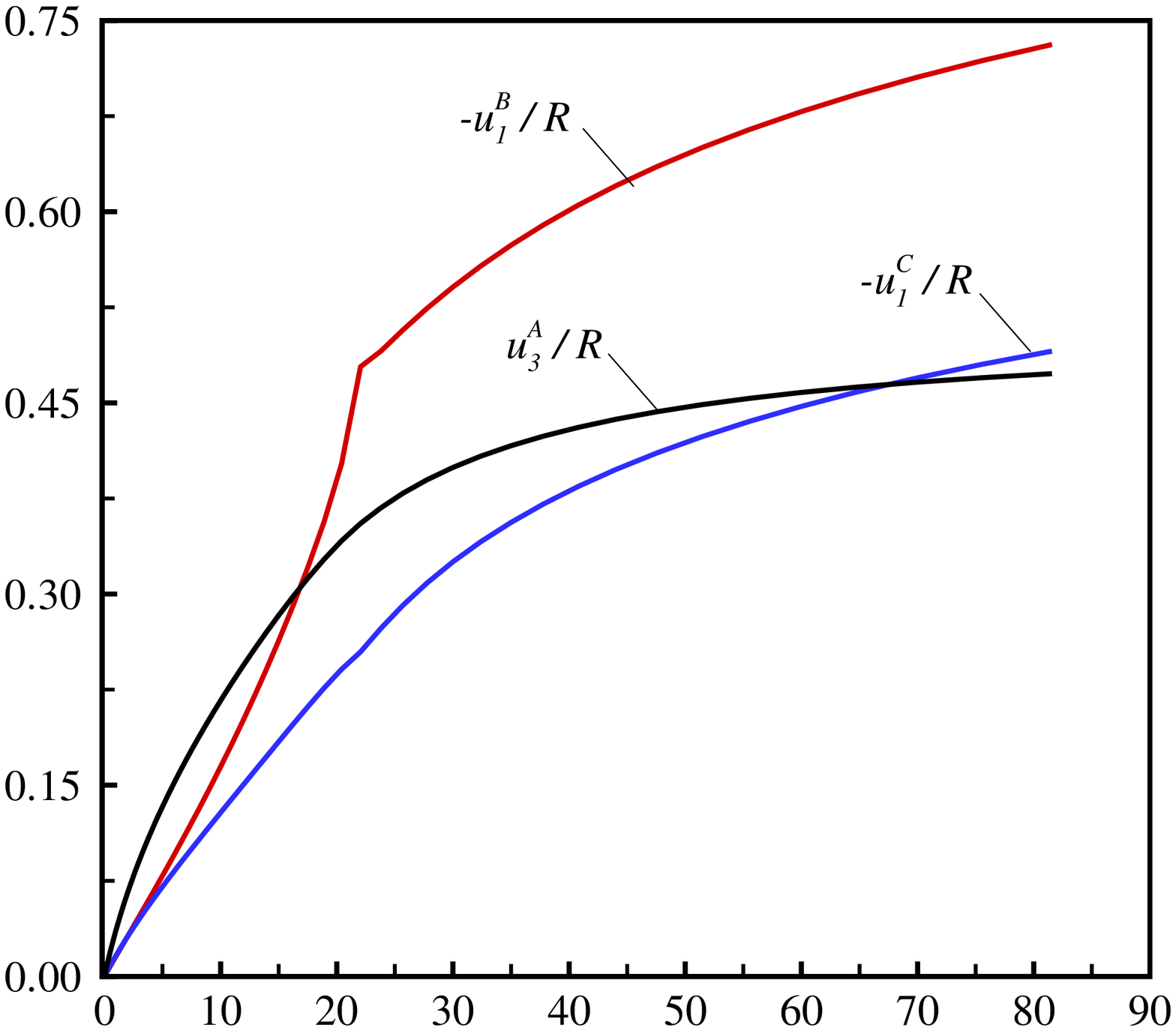} 
\put(-205,30){\includegraphics[width=75mm]{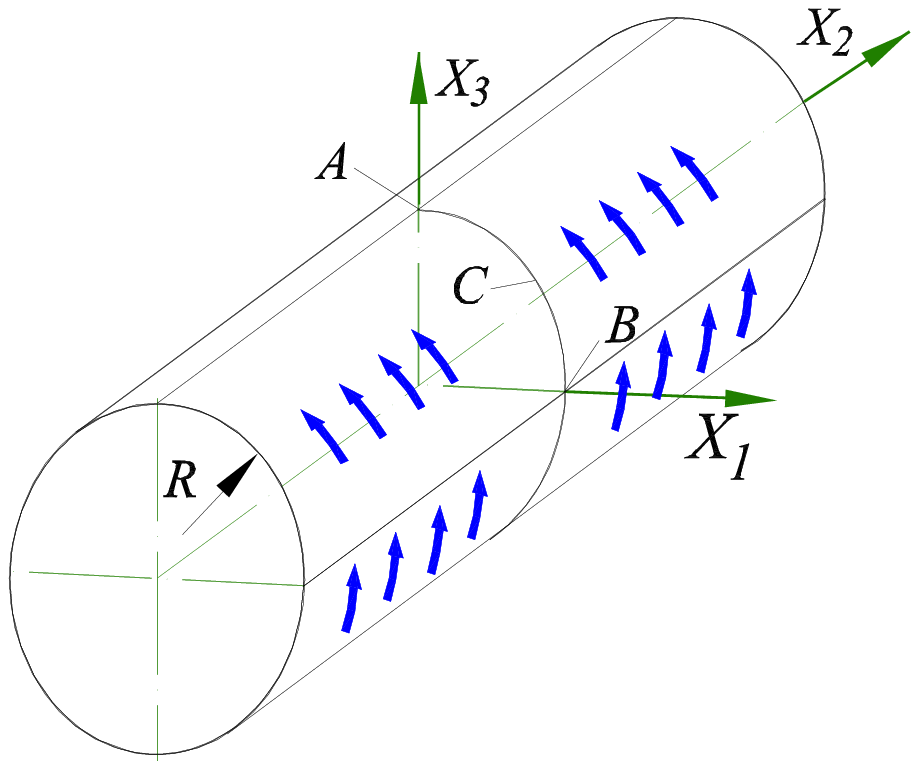}}
\put(-50,155){\scriptsize(\textbf{a})}
\put(-150,4){\scriptsize $\frac{1}{\mu \mu_0} |\mbbB^{\text{ext}}||\tilde{\mbbB}^{\text{rem}}| \times 10^3$}
\end{minipage}
\begin{minipage}[b]{.45\textwidth}
\includegraphics[trim=-4cm -5cm  0cm  0cm, width=55mm]{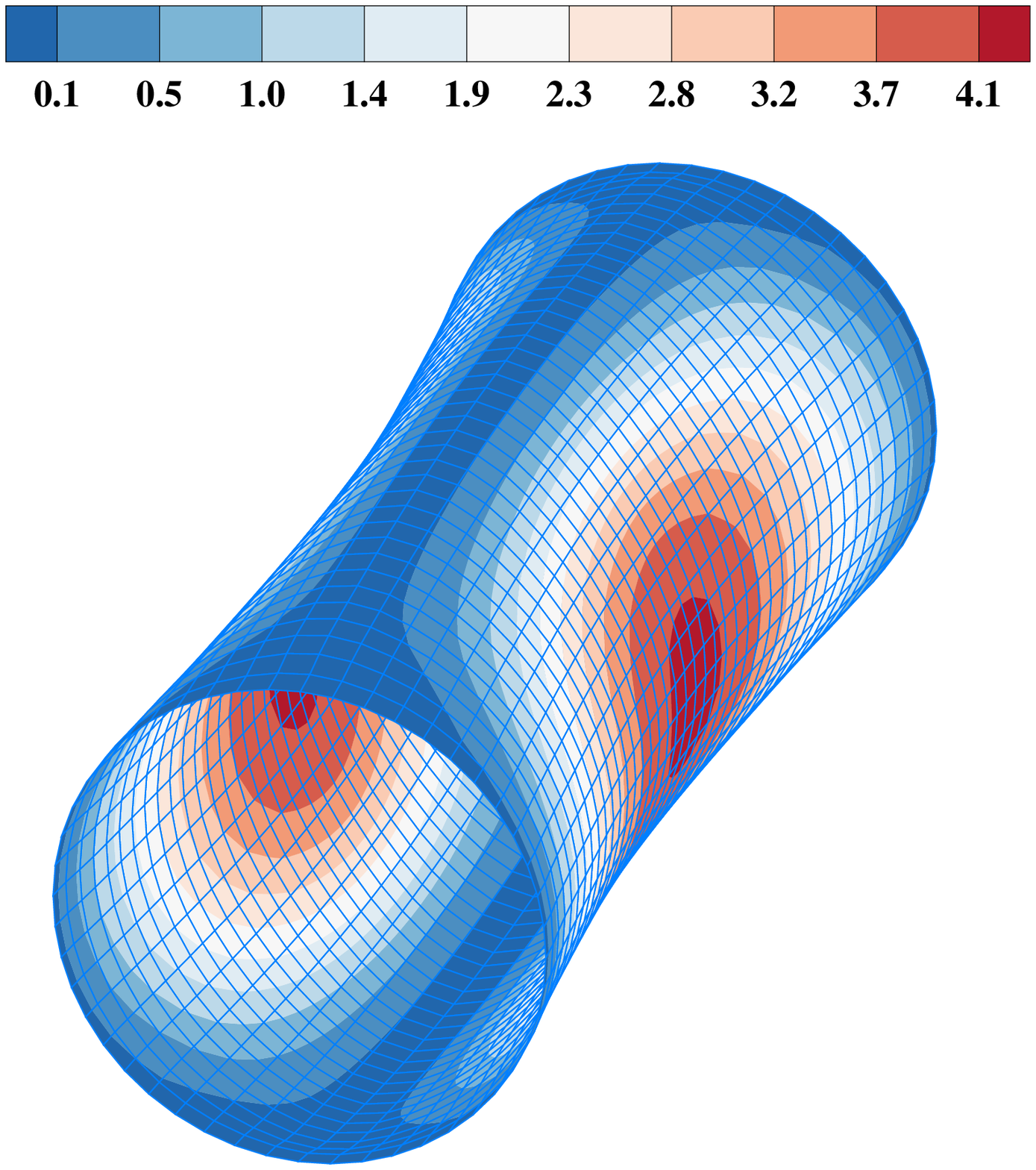} 
\put(-120,110){\scriptsize(\textbf{b})}
\end{minipage} \\
\vspace{10pt}
\begin{minipage}[b]{.32\textwidth}
\includegraphics[trim=0cm 0cm  0cm  0cm, width=48mm]{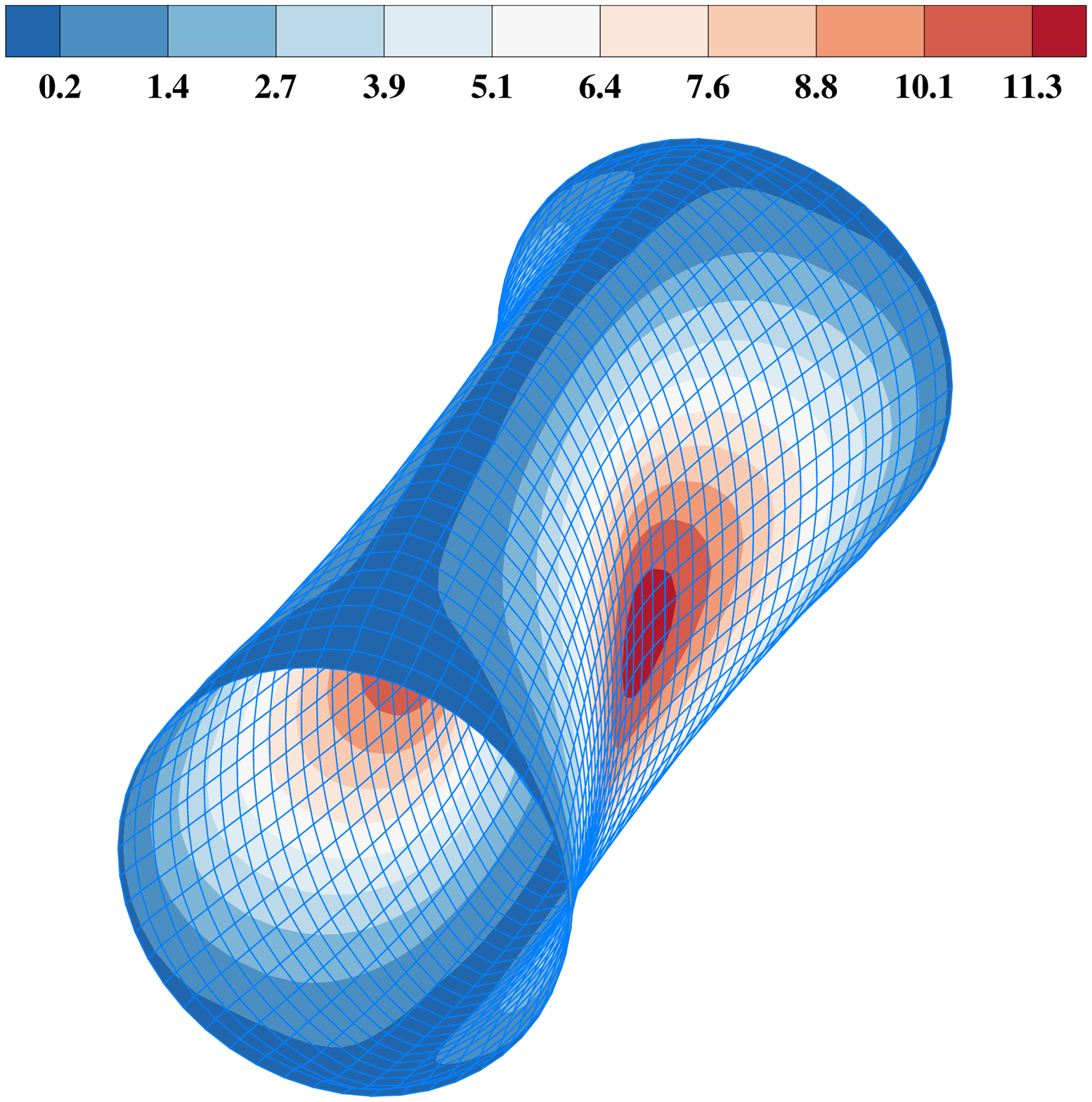} 
\put(-120,80){\scriptsize(\textbf{c})}
\end{minipage}   
\begin{minipage}[b]{.32\textwidth}
\includegraphics[trim=0cm 0mm  0cm  0cm, width=48mm]{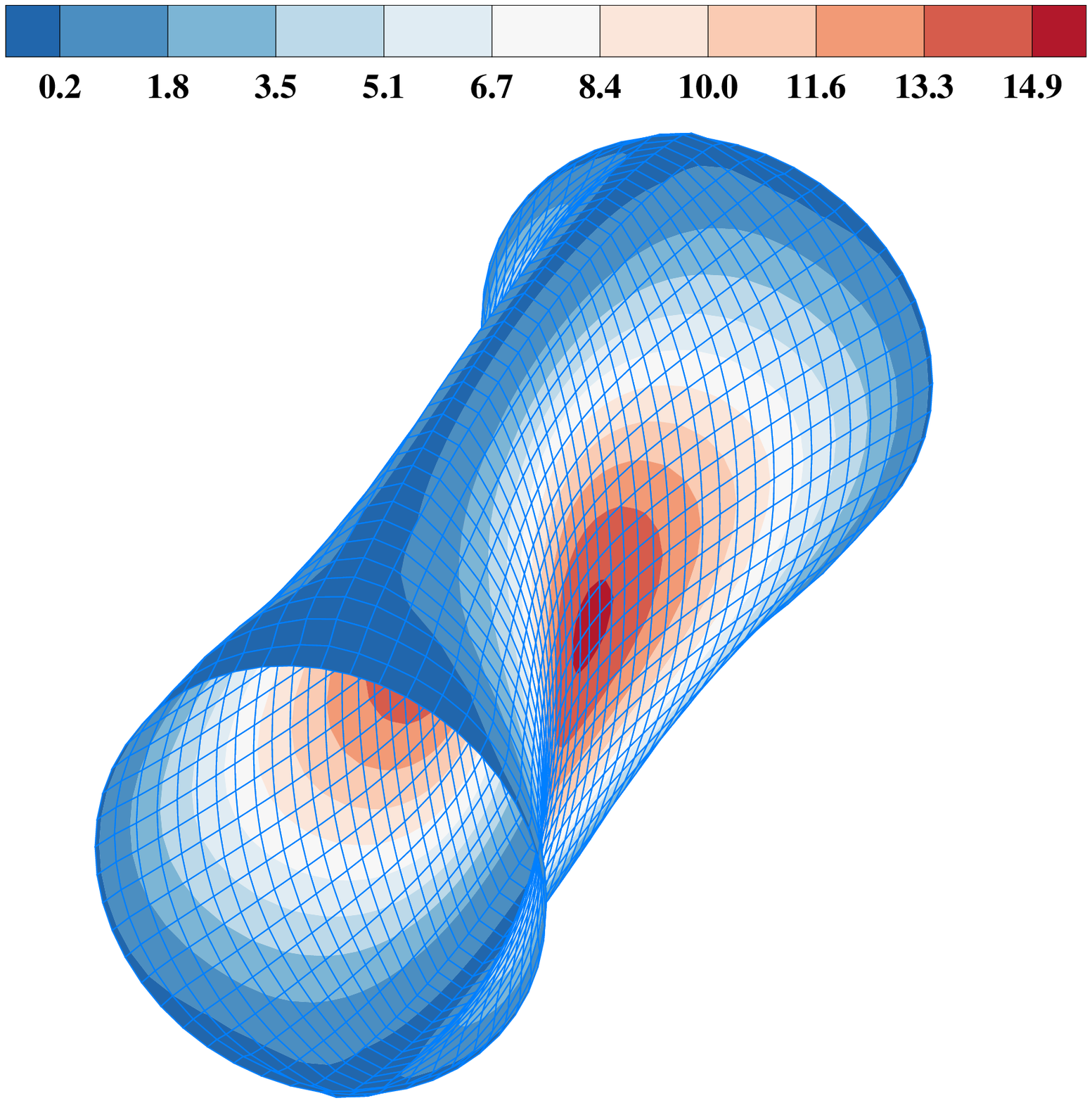} 
\put(-120,80){\scriptsize(\textbf{d})}
\end{minipage}
\begin{minipage}[b]{.32\textwidth}
\includegraphics[trim=0cm 0cm  0cm  0cm, width=48mm]{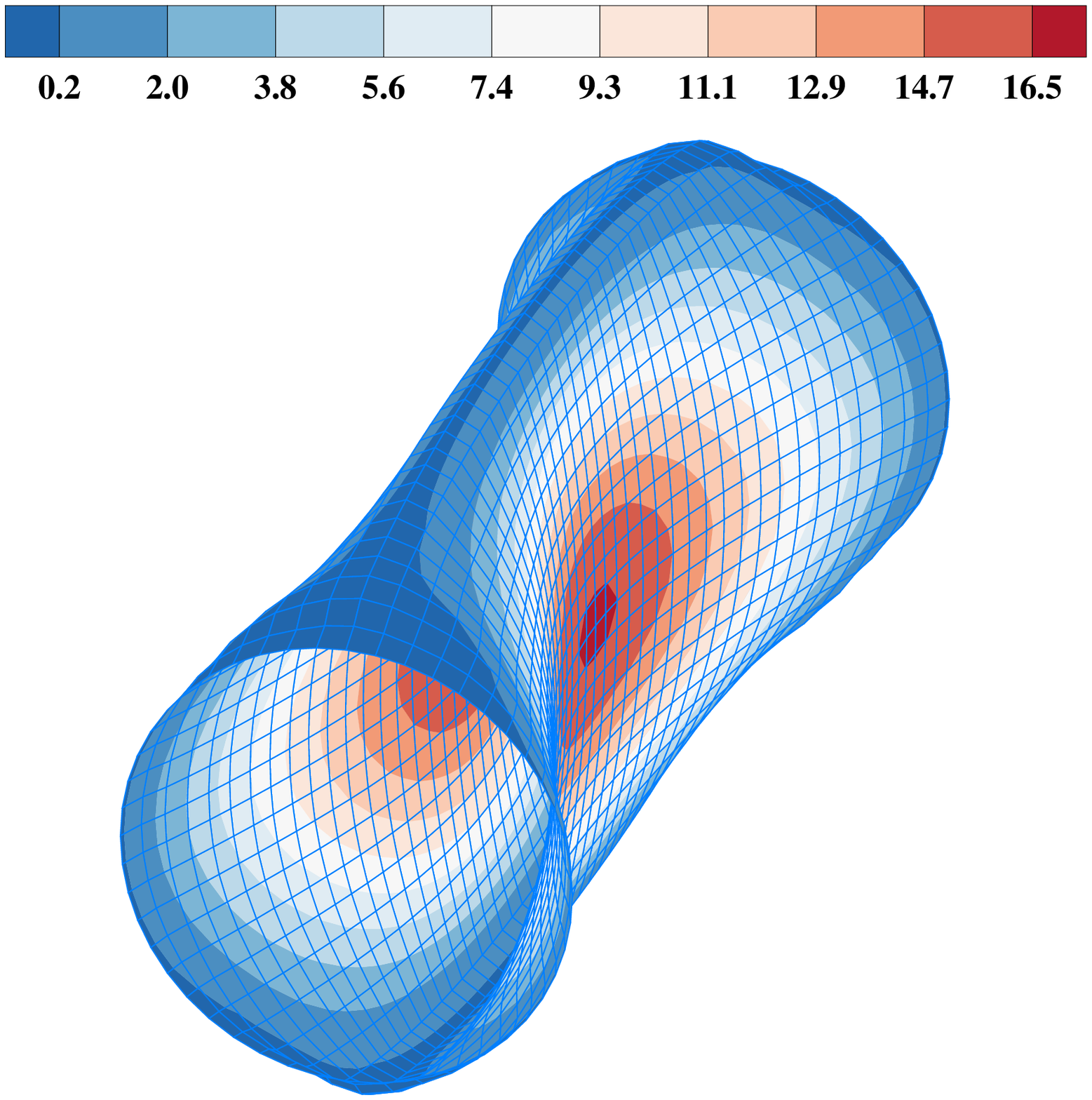} 
\put(-120,80){\scriptsize(\textbf{e})}
\end{minipage}
\vspace{0pt}
\caption{Deformation of a cylinder under magnetic loading,
(\textbf{a}): load-displacement curves, 
(\textbf{b},\textbf{c},\textbf{d},\textbf{e}): sequences of deformation (with the contours of $|u_1|$ in mm) for $|\mbbB^\text{ext}| \in \{20,50, 100, 150 \}$ (mT)}
\label{Ex5Fig1} 
\end{figure*}
\subsection{A magnetic gripper}
\label{gripper}
\vspace{-5pt}
The elastic response of a spherical gripper is simulated in this example. 
Soft grippers made of magneto-active materials have the potential as actuating components in soft robotics. 
For instance, Ju et al.~\cite{ju2021} and Carpenter et al. \cite{carpenter2021} demonstrated additively manufactured magneto-active grippers while Kadapa and Hossain~\cite{kadapa2022} simulated the viscoelastic influences of underlying polymeric materials.  
In our case, the gripper is composed of $12$ equal arms.
In the undeformed configuration, the arms cover the surface of an incomplete sphere of radius $R$.
It is assumed that the mechanical and magnetic properties, and the thickness of the HMSM are the same as those given in the example \ref{Cross}.
The geometry of a single arm is shown in Fig.~\ref{Ex6Fig1}(\textbf{a}).
The arc $DE$ lies in the $X_1X_2$ plane, its length is $12$ mm, and covers 
$30^{\circ}$ of a full circle.
Therefore, the mean radius of the arm is $R=\frac{12}{\pi / 6} =22.92$ mm.
The arc $AC$ lies in the $X_1X_3$ plane and its length is $60$ mm.
The angle between the radius $OA$ and the $X_3$-axis is $15^{\circ}$, and the geometry is symmetric w.r.t. the $X_1X_2$ plane.
Moreover, the topmost arc of the arm is assumed to be clamped.
As shown in the figure, let $\mbbe_{\varphi}$ be the standard meridian unit tangent vector to the sphere. 
It is assumed $\tilde{\mbbB}^{\text{rem}}$ is along 
$\mbbe_{\varphi}$ for $X_3>0$, and along 
$-\mbbe_{\varphi}$ for $X_3<0$.
It is noted that applying $\mathbbm{B}^{\text{ext}}$ in $\mathbbm{e}_3$ direction opens the arms of the gripper.
Here, the maximum magnetic loading
$\mathbbm{B}^{\text{ext}}_{\text{max}}= 10 \mathbbm{e}_3$ (mT) 
acts on the arms.

Numerical experiments indicate that a $6 \times 30$ mesh of shell elements in the arm provides convergence in the results. 
The displacement components $u_1$ and $u_3$ at the points $B$ and $C$  versus $\frac{10^3}{\mu \mu_0}|\mbbB^{\text{ext}}||\tilde{\mbbB}^{\text{rem}}|$
are plotted in Fig.~\ref{Ex6Fig1}(\textbf{a}).
For a single arm under the maximum external magnetic flux of $10$ mT, the maximum value of the displacement component $u_3$ is obtained to be about $43.9$ mm.
The deformed shapes of the gripper under four different values of the external magnetic flux are illustrated in Figs.~\ref{Ex6Fig1}(\textbf{b,c,d,e}).
It is noted that the maximum value of the external magnetic flux to avoid intersection between the arms is $6.8$ mT.
In this case, the maximum $u_3$ component of displacement is about $41.6$ mm.
\begin{figure*}
\centering
\begin{minipage}[b]{.5\textwidth}
\includegraphics[width=85mm]{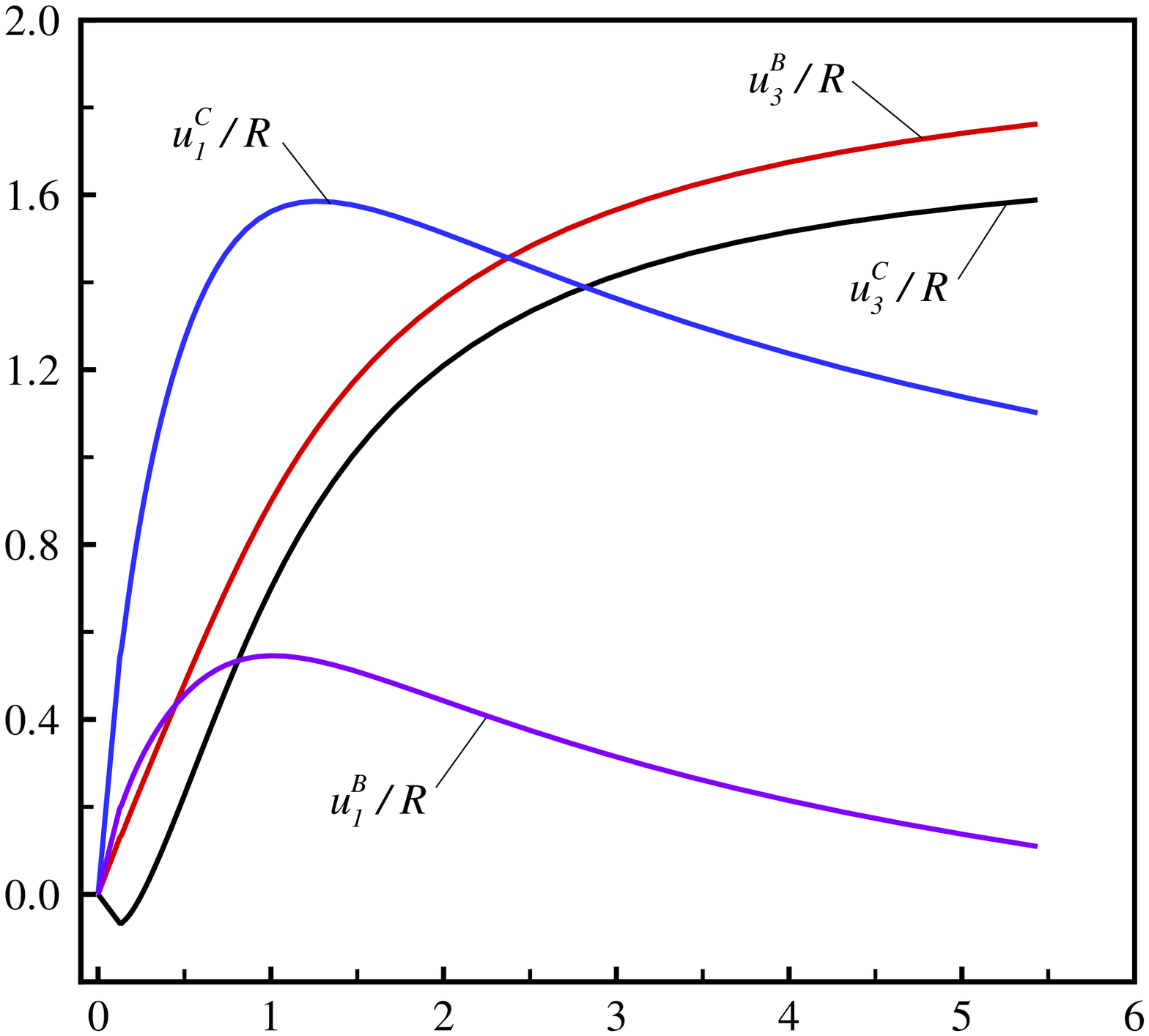} 
\put(-170,55){\includegraphics[width=65mm]{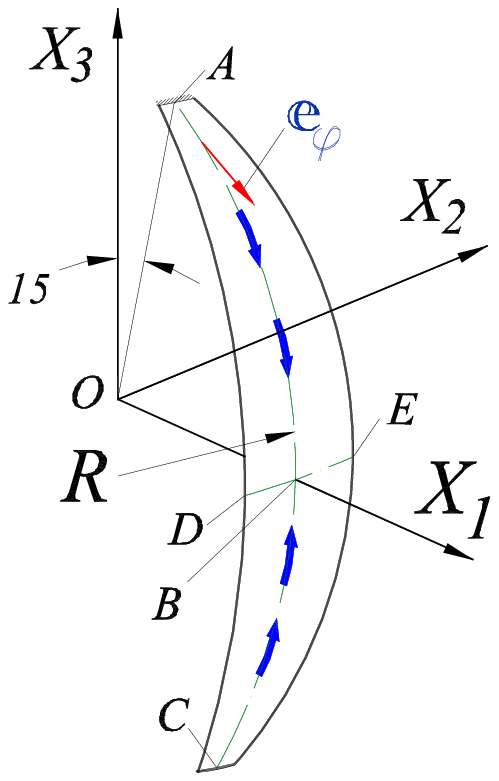}}
\put(-130,180){\scriptsize(\textbf{a})}
\put(-150,4){\scriptsize $\frac{1}{\mu \mu_0} |\mbbB^{\text{ext}}||\tilde{\mbbB}^{\text{rem}}| \times 10^3$}
\end{minipage}
\begin{minipage}[b]{.4\textwidth}
\includegraphics[trim=-2cm -5cm  0cm  0cm, width=55mm]{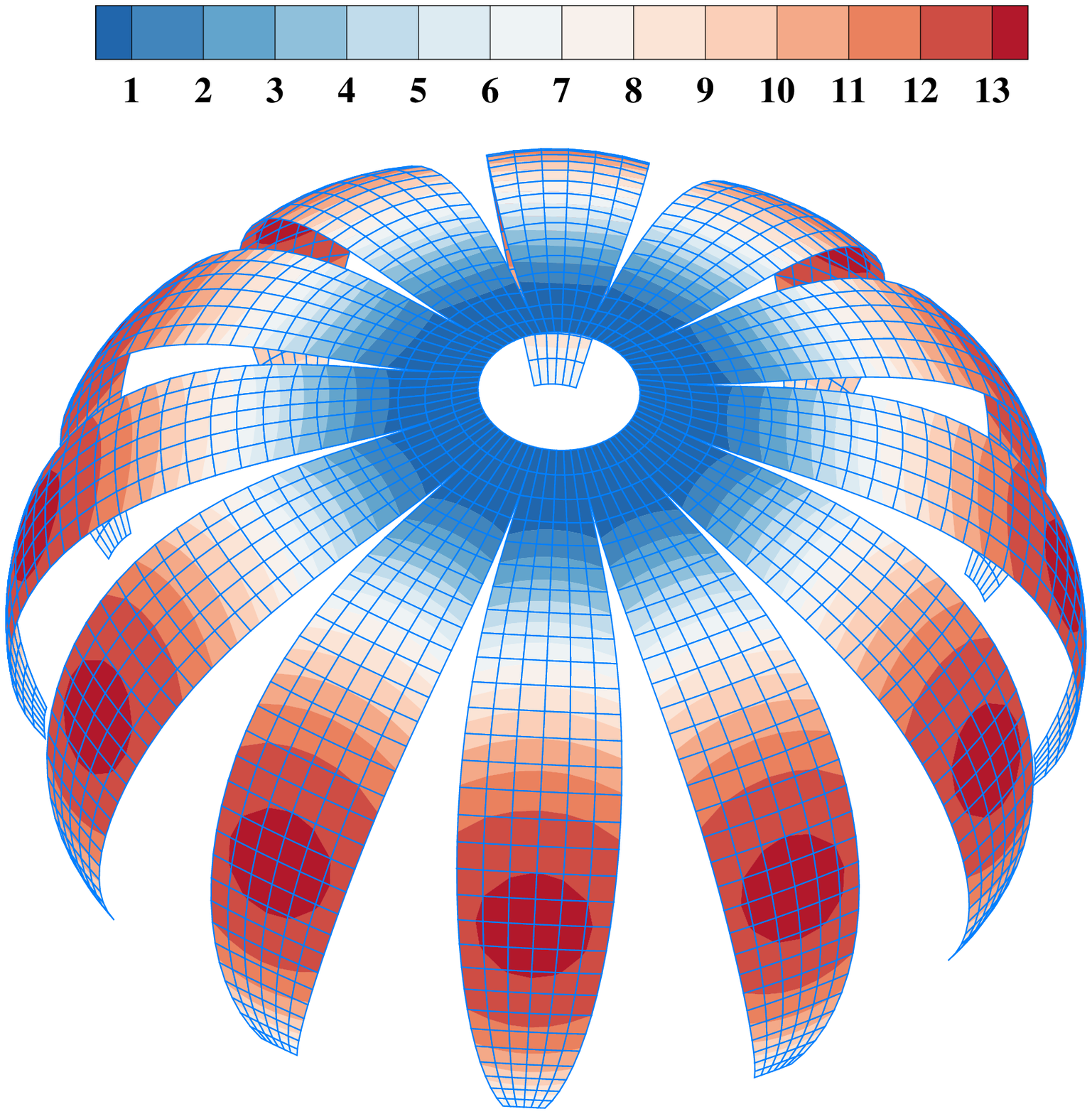} 
\put(-140,120){\scriptsize(\textbf{b})}
\end{minipage} \\
\vspace{10pt}
\begin{minipage}[b]{.32\textwidth}
\includegraphics[trim=0cm 0cm  0cm  0cm, width=49mm]{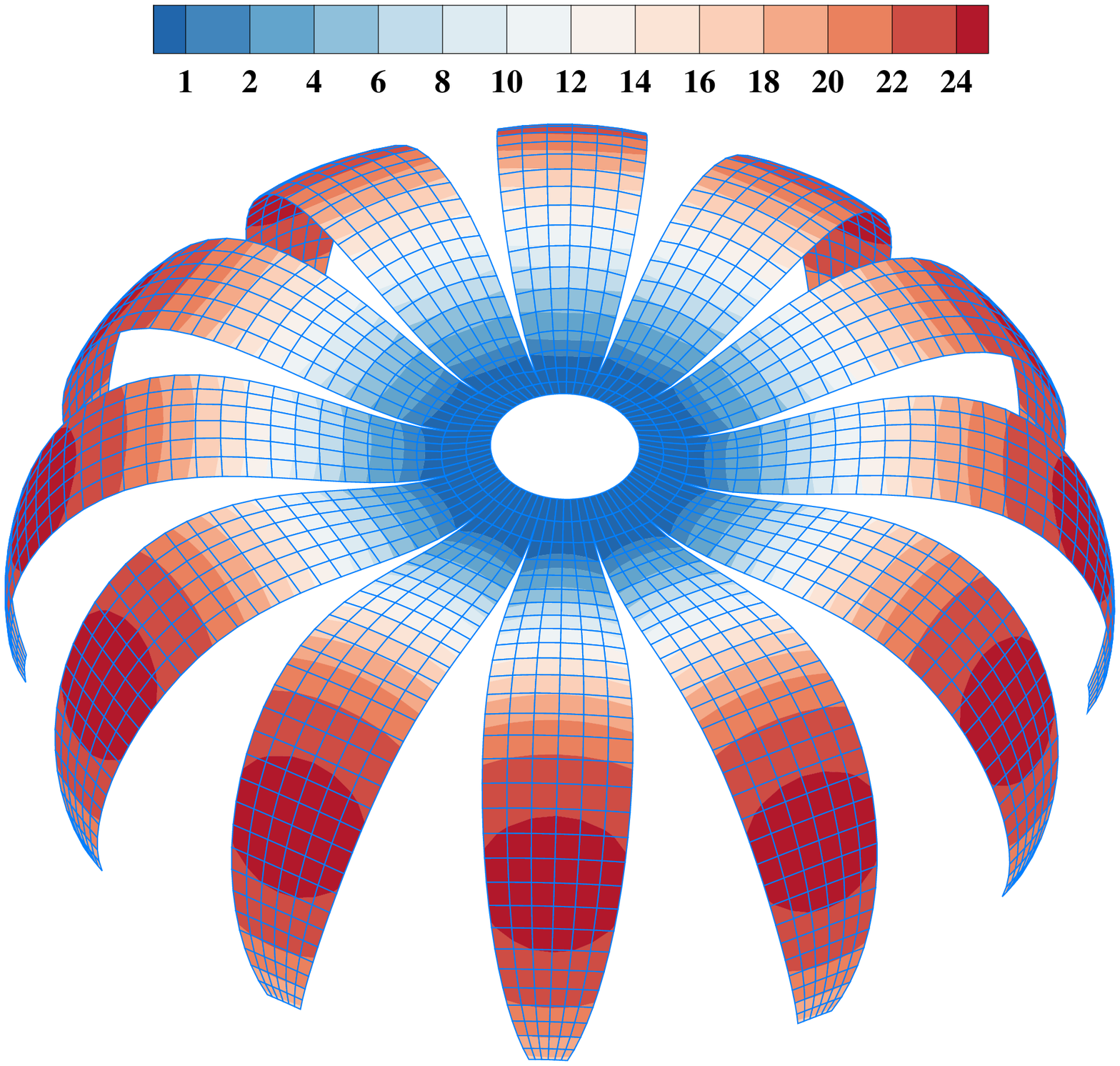} 
\put(-140,95){\scriptsize(\textbf{c})}
\end{minipage}   
\begin{minipage}[b]{.32\textwidth}
\includegraphics[trim=0cm 0mm  0cm  0cm, width=50mm]{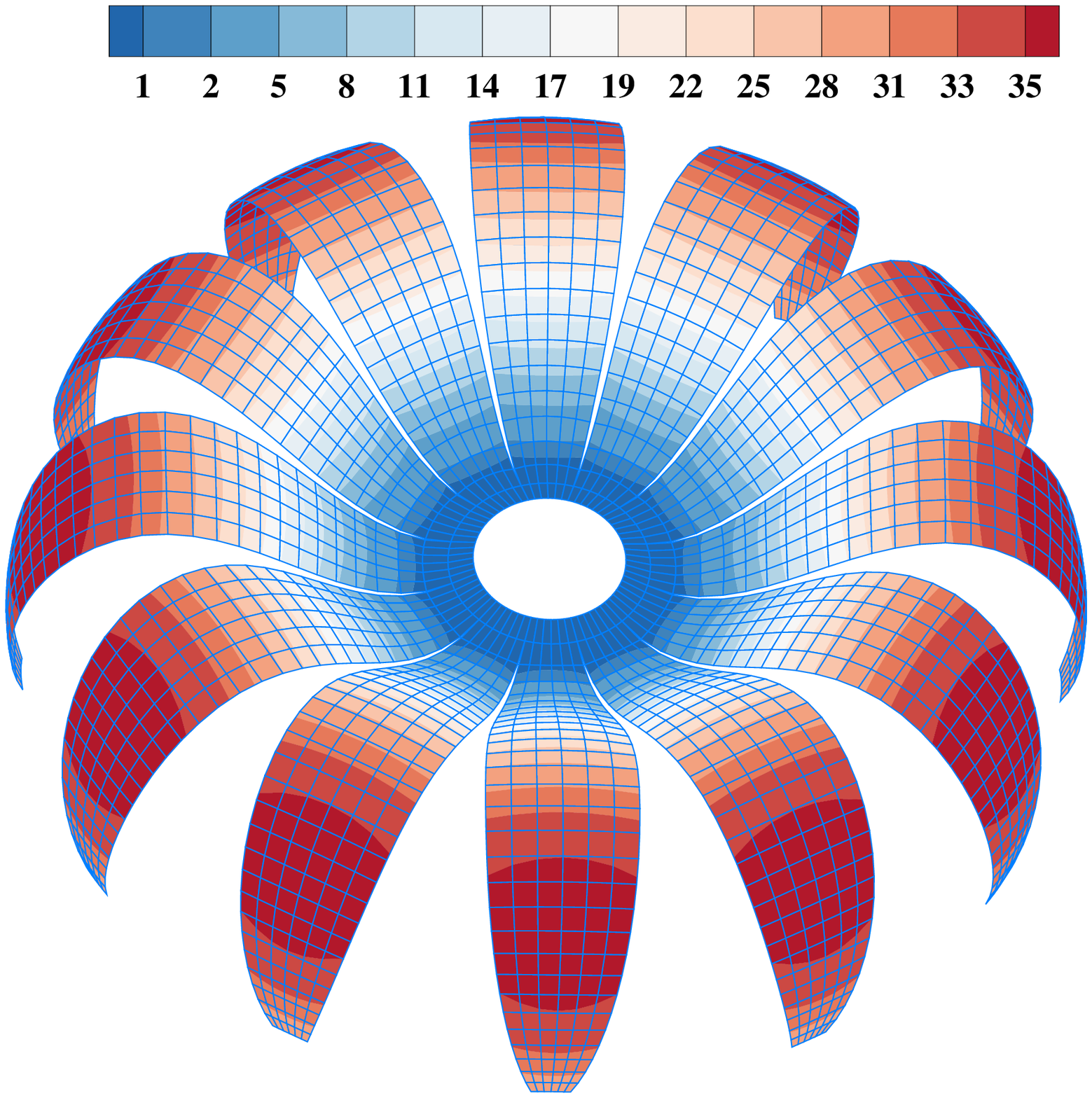} 
\put(-140,95){\scriptsize(\textbf{d})}
\end{minipage}
\begin{minipage}[b]{.32\textwidth}
\includegraphics[trim=0cm 0cm  0cm  0cm, width=50mm]{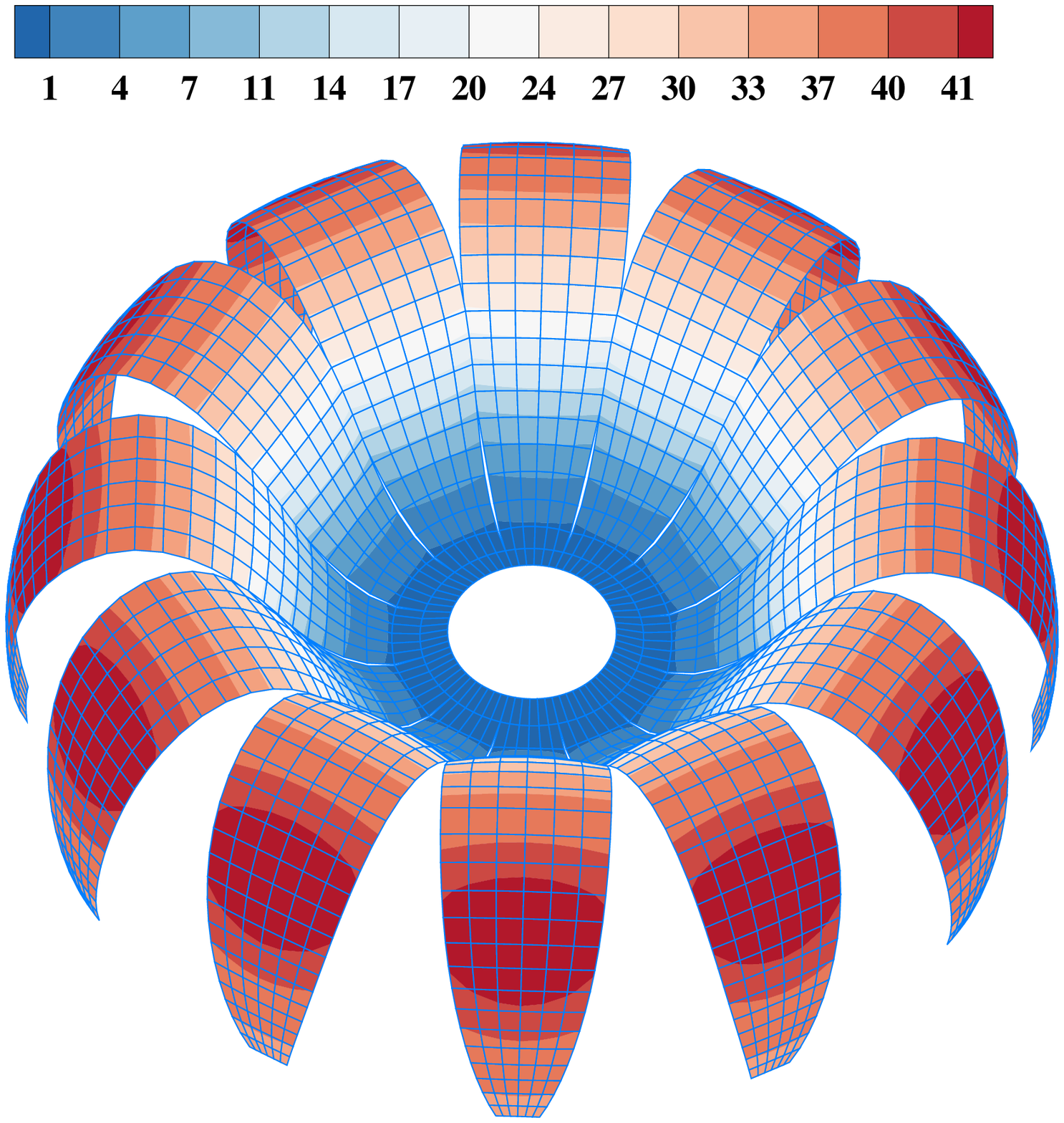} 
\put(-130,95){\scriptsize(\textbf{e})}
\end{minipage}
\vspace{0pt}
\caption{Deformation of a spherical gripper with $12$ arms, 
(\textbf{a}): load-displacement curves, 
(\textbf{b},\textbf{c},\textbf{d},\textbf{e}): sequences of deformation (with the contours of $u_3$ in mm) for $|\mbbB^\text{ext}| \in \{1,2,4,6.8\}$ (mT)}
\label{Ex6Fig1} 
\end{figure*}
\vspace{10pt}
\section{Summary}
\label{conc}
\vspace{-5pt}
In this research, a 10-parameter micropolar shell model for simulating the finite elastic deformation of thin hard-magnetic soft structures was formulated.
The idea of employing the micropolar theory comes from the fact that magnetic stimulation induces a body couple on these materials, which in turn leads to asymmetric Cauchy stress tensor. 
Since the governing equations at finite strains, including magnetic effects, cannot be solved analytically, a nonlinear finite element formulation for simulating the problems of arbitrary thin geometry, boundary conditions, and loading cases was also presented.
Six different numerical examples were solved to assess the applicability of the present formulation.
It was shown that the results of the proposed formulation are in good agreement with the available experimental and numerical ones.
The viscoelastic and thermal effects will be taken into account in the forthcoming contributions.
%
\vspace{10pt}
\section*{Declaration of competing interest}
\label{dci}
\vspace{-5pt}
The authors declare no competing interests.


\vspace{10pt}
\section*{Acknowledgements}
\label{dci}
\vspace{-5pt}
M. Hossain acknowledges the funding through an  Engineering and Physical Sciences Research Council (EPSRC)  Impact Acceleration Award (EP/R511614/1). He also acknowledges the support by EPSRC through the Supergen ORE Hub (EP/S000747/1), which has awarded funding for the Flexible Fund project Submerged bi-axial fatigue analysis for flexible membrane Wave Energy Converters (FF2021-1036).




\vspace{10pt}
\begin{thebibliography}{9}

\bibitem{Ren2019}  
Ren Z,
Hu W, 
Dong X, 
Sitti M.
Multi-functional soft-bodied jellyfish-like swimming. 
\textit{Nature Commun.} 
2019;10:2703.


\bibitem{Wu2020}  
Wu S, 
Hu W, 
Ze Q,   
Sitti M, 
Zhao R.
Multifunctional magnetic soft composites: a review. 
\textit{Multifuct Mater.} 
2020;3:042003.

\bibitem{Bastola2020} 
Bastola AK, 
Paudel M,  
Li L, 
Li W.
Recent progress of magnetorheological elastomers: a review.
\textit{Smart Mater Struct.} 
2020;29:123002.

\bibitem{Bastola2021} 
Bastola AK,  
Hossain M.
The shape-morphing performance of magnetoactive soft materials.  
\textit{Mat Des.}
2021;211:110172.
%
%

\bibitem{Lucarini2022}  
Lucarini S,   
Hossain M,  
Garcia-Gonzalez D. 
Recent advances in hard-magnetic soft composites: synthesis, characterisation, computational modelling, and applications. 
\textit{Compos Struct.}
2022;200:210001.

\bibitem{Yarali2022}  
Yarali E, 
Banishadi M, 
Zolfagharian A, 
Chavoshi M, 
Arefi F, 
Hossain M, 
Bastola B, 
Ansari M, 
Foyouzat A, 
Dabbagh A, 
Ebrahimi M, 
Mirzaali MJ, 
Bodaghi M. 
Magneto/electroresponsive polymers toward manufacturing, characterization, and biomedical/soft robotic applications.
\textit{Appl Mater Today.}
2022;26:101306.


\bibitem{Saxena2013}  
Saxena P, 
Hossain M,
Steinmann P. 
A theory of finite deformation magneto-viscoelasticity. 
\textit{Int J Solids Struct.}
2013;50:3886--3897.

\bibitem{Miehe2016} 
Ethiraj G, 
Miehe C.
Multiplicative magneto-elasticity of magnetosensitive polymers incorporating micromechanically-based network kernels.
\textit{Int J Eng Sci.} 
2016;102;93--119.

\bibitem{Mehnert2017}  
Mehnert M,  
Hossain M,  
Steinmann P. 
Towards  a  thermo-magneto-mechanical  coupling  framework  for  magneto-rheological  elastomers.
\textit{Int J Solids Struct.} 
2017;128:117--132.


\bibitem{Mukherjee2020}  
Mukherjee  D,  
Bodelot L, 
Danas K. 
Microstructurally-guided explicit continuum  models  for  isotropic  magnetorheological  elastomers  with  iron particles.
\textit{Int J Non-Linear Mech. }
2020;120:103380.


\bibitem{Bustamante2021} 
Bustamante R,   
Shariff MHBM,     
Hossain M.
Mathematical  formulations for  elastic  magneto-electrically  coupled  soft  materials  at  finite  strains: Time-independent processes.
\textit{Int J Eng Sci.}
2021;159:103429. 


\bibitem{Hu2022}  
Hu X, 
Zhu H, 
Chen S, 
Yu H, 
Qu S.
Magnetomechanical behavior of soft magnetoactive membranes.
\textit{Int J Solids Struct.}
2022;234--235:111310.

\bibitem{akbari2020}
Akbari E,  
Khajehsaeid H.
A continuum magneto-mechanical model for magnetorheological elastomers.
\textit{Smart Mater Struct.} 
2020;30:015008.

\bibitem{Schumann2020}  
Sch\"umann M,  
Borin DY, 
Morich J,   
Odenbach S.
Reversible   and non-reversible  motion  of  NdFeB-particles  in  magnetorheological  elastomers. 
\textit{J  Intell  Mater  Syst  Struct.}   
2020;32:3--15.

\bibitem{Lee2020}  
Lee M,  
Park T, 
Kim C,  
Park SM,
Characterization of a magneto-active  membrane  actuator  comprising  hard  magnetic  particles  with  varying crosslinking degrees.
\textit{Mater Des.}
2020;195:108921.


\bibitem{Lum2016}  
Lum GZ,  
Ye Z,  
Dong X,  
Marvi H,  
Erin O, 
Hu W,  
Sitti M. 
Shape-programmable magnetic soft matter.
\textit{Proc Natl Acad Sci.}
2016;113:6007--6015.

\bibitem{Wu2019}  
Wu S, 
Ze Q, 
Zhang R, 
Hu N,  
Cheng Y,  
Yang F, 
Zhao R. 
Symmetry-breaking actuation mechanism for soft robotics and active metamaterials.
\textit{ACS Appl Mater Interfaces.}
2019;11:41649--41658.


\bibitem{Kim2018}  
Kim Y,  
Yuk H,  
Zhao R,  
Chester SA,  
Zhao X. 
Printing  ferromagnetic domains for untethered fast-transforming soft materials.
\textit{Nature}, 
2018;558:274--279.


\bibitem{Alapan2020}  
Alapan Y, 
Karacakol AC, 
Guzelhan SN, 
Isik I, 
Sitti M.   
Reprogrammable shapemorphing of magnetic soft machines.
\textit{Sci Adv.}
2020;6:eabc6414.



\bibitem{Kuang2021} 
Kuang X, 
Wu S, 
Ze Q,  
Yue L,  
Jin Y, 
Montgomery SM,  
Yang F, 
Qi HJ,  
Zhao R.
Magnetic dynamic polymers for modular assembling and reconfigurable morphing architectures.
\textit{Adv Mater.}
2021;2102113.


\bibitem{Wang2021}  
Wang L,  
Zheng D, 
Harker P,  
Patel AB,  
Guo CF,  
Zhao X.  
Evolutionary design of magnetic soft continuum robots.
\textit{Proc Natl Acad Sci.}
2021;118:21.


\bibitem{Wu2021}  
Wu S, 
Hamel CM, 
Ze Q, 
Yang F, 
Qi HJ, 
Zhao R. 
Evolutionary algorithm-guided  voxel-encoding  printing  of  functional  hard-magnetic  soft  active materials. 
\textit{Adv Intell Syst.} 
2021;2:2000060.

\bibitem{Kalina2017} 

Kalina KA, 
Brummund J, 
Metsch P, 
Kaestner M, 
Borin DY, 
Linke JM, 
Odenbach S.
Modeling of magnetic hystereses in soft MREs filled with NdFeB particles.
\textit{Smart Mater Struct.} 
2017;26:105019.


\bibitem{Zhao2019}  

Zhao R,   
Kim Y,    
Chester AS,    
Sharma P, 
Zhao X.  
Mechanics of hard-magnetic soft materials. 
\textit{J Mech Phys Solids.}
2019;124:244--263.


\bibitem{Garcia2019} 
Garcia-Gonzalez D.  
Magneto-visco-hyperelasticity  for  hard-magnetic  soft materials:  theory  and  numerical  applications.
\textit{Smart  Mater  Struct.}  
2019;28:085020.

\bibitem{Mukherjee2021}  

Mukherjee D,
Rambausek M,
Danas K.
An explicit dissipative model for isotropic  hard  magnetorheological  elastomers.  
\textit{J Mech Phys Solids.}  
2021;151:104361.


\bibitem{rambausek2022}  


Rambausek M,
Mukherjee D,
Danas K.
A computational framework for magnetically hard and soft viscoelastic magnetorheological elastomers.
\textit{Comput  Methods Appl Mech  Eng.}  
2021;391:114500.


\bibitem{Zhang2020}  

Zhang R,  
Wu S,  
Qiji Z,
Zhao Z. 
Micromechanics  study  on  actuation efficiency of hard-magnetic soft active materials.
\textit{J Appl Mech.} 
2020;87:091008.


\bibitem{Garcia2021a} 
Garcia-Gonzalez D, 
Hossain M. 
A  microstructural-based  approach  to model magneto-viscoelastic materials at finite strains.
\textit{Int J Solids Struct.}
2021;208--209:119--132.


\bibitem{Garcia2021b} 
Garcia-Gonzalez D, 
Hossain M.
Microstructural modelling of hard-magnetic soft materials:
Dipole--dipole interactions versus Zeeman effect.
\textit{Extreme Mech Lett.}
2021;48:101382.

\bibitem{Ye2021}  

Ye H,  
Li Y,  
Zhang T.
Magttice:  A  lattice  model  for  hard-magnetic  soft materials.
\textit{Soft Matter}.
2021;17:3560--3568.

\bibitem{DH2022IJSS} 

Dadgar-Rad F,  
Hossain M. 
Finite deformation analysis of hard-magnetic soft materials based on micropolar continuum theory. 
\textit{Int J Solids Struct.}
2022;251:111747.




\bibitem{Yan2021}  

Yan D, 
Abbasi A, 
Reis PM. 
A comprehensive framework for hard-magnetic beams: reduced-order theory, 3D simulations, and experiments.
\textit{Int J Solids Struct.}  
2021; https://doi.org/10.1016/j.ijsolstr.2021.111319.




\bibitem{Wang2020}  

Wang L, 
Kim Y, 
Guo GF, 
Zhao X. 
Hard-magnetic elastica.
\textit{J Mech Phys Solids}.
2020;142:104045.


\bibitem{Rajan2021}  

Rajan A, 
Arockiarajan A.
Bending of hard-magnetic soft beams: A finite elasticity approach with
anticlastic bending.
\textit{Eur J Mech A/Sol.}
2021;90:104374.


\bibitem{Chen2020a} 

Chen W,    
Yan Z,  
Wang L.
Complex  transformations  of  hard-magnetic  soft  beams  by  designing
residual magnetic flux density.
\textit{Soft Matter}, 
2020;16:6379--6388. 

\bibitem{Chen2021} 

Chen W, 
Wang L, 
Yan Z, 
Luo B.
Three-dimensional large-deformation model of hard-magnetic soft beams.
\textit{Compos Struct.}
2021;266:113822.




\bibitem{Reis2021}

Yan D, 
Pezzulla M, 
Cruveiller L, 
Abbasi A, 
Reis PM.
Magneto-active elastic shells with tunable buckling strength.
\textit{Nature Commun.} 
2021;12:2831.

\bibitem{DH2022EML} 

Dadgar-Rad  F,  
Hossain  M.
Large viscoelastic deformation of hard-magnetic soft beams. 
\textit{Extreme Mech Lett.}
2022;54:101773.


\bibitem{Dorfmann2014} 
Dorfmann A,  
Ogden RW.  
\textit{Nonlinear Theory of Electroelastic and Magnetoelastic Interactions}.
Springer; 2014.

\bibitem{Kafadar1971} 
Kafadar CB, 
Eringen AC. 
Micropolar media--I the classical theory.
\textit{Int J Eng Sci.}
1971;9:271--307.

\bibitem{Eringen1976} 
Eringen AC,
Kafadar CB. 
Polar  field  theories.
In:  Eringen AC (Ed.),
\textit{Continuum Physics}, vol. IV. 
Academic Press; 1976; 1--73.


\bibitem{Eringen1999} 
Eringen AC. 
\textit{Microcontinuum Field Theories, vol. I, Foundations and Solids}.
Springer; 1999.




\bibitem{Borst1993} 
de Borst  R.
A generalization of  $J_2$-flow theory for polar continua.
\textit{Comput Methods Appl Mech Eng.}
1993;103:347--362.

\bibitem{Steinmann1994}  
Steinmann  P. 
A micropolar theory of finite deformation and finite rotation multiplicative elastoplasticity.
\textit{Int J Solids Struct.}
1994;31:1063--1084.

\bibitem{Tsakmakis2005}  
Grammenoudis  P, 
Tsakmakis  C.
Finite element implementation of large deformation micropolar plasticity exhibiting isotropic and kinematic hardening effects.
\textit{Int J Numer Meth Eng.}
2005;62:1691--1720.


\bibitem{Grammenoudis2007a} 

Grammenoudis  P, 
Tsakmakis  C. 
Micropolar plasticity theories and their classical limits. Part I: Resulting model.
\textit{Acta Mech.}
2007;189:151--175.
\bibitem{Grammenoudis2007b} 

Grammenoudis  P, 
Sator C,
Tsakmakis  C.
Micropolar plasticity theories and their classical limits. Part II: Comparison of responses predicted by the limiting and a standard classical model. 
\textit{Acta Mech.}
2007;189:177--191.


\bibitem{Bauer2012a} 

Bauer S, 
Dettmer WG, 
Peric D, 
Sch\"afer M.  
Micropolar hyper-elastoplasticity:  constitutive model, consistent linearization, and simulation of 3D scale effects.  
\textit{Int J Numer Meth Eng.}
2012;91:39--66.





\bibitem{Borst2022} 
de Borst  R,
Alizede Sabet  S,
Hageman T.
Non-associated Cosserat plasticity.
\textit{Int J Mech Sci.}
2022; https://doi.org/10.1016/j.ijmecsci.2022.107535.


\bibitem{Ramezani2008}  
Ramezani S,
Naghdabadi R, 
Sohrabpour S.  
Non-linear finite element implementation of micropolar hypo-elastic materials.
\textit{Comput Methods Appl Mech Eng.}
2008;197:4149--4159.
\bibitem{Ramezani2009}  
Ramezani S,
Naghdabadi R, 
Sohrabpour S. 
Constitutive equations for micropolar hyper-elastic materials. 
\textit{Int J Solids Struct.}
2009;46:2765--2773.

\bibitem{Pietraszkiewicz2009}  
Pietraszkiewicz W, 
Eremeyev VA. 
On natural strain measures of the non-linear micropolar continuum.
\textit{Int J Solids Struct.}
2009;46:774--787.

\bibitem{Bauer2010} 
Bauer S, 
Sch\"afer M, 
Grammenoudis P,
Tsakmakis C.
Three-dimensional  finite  elements  for  large  deformation  micropolar  elasticity. 
\textit{Comput Methods Appl Mech Eng.}  
2010;199:2643--2654.



\bibitem{Bauer2012b} 
Bauer S,
Dettmer WG, 
Peric D, 
Sch\"afer M. 
Micropolar hyper-elasticity: constitutive model, consistent linearization and simulation of 3D scale effects. 
\textit{Comput Mech.}
2012;50:383--396.

\bibitem{Erdelj2020} 

Erdelj SG, 
Jeleni\'c  G, 
Ibrahimbegovi\'c  A.
Geometrically non-linear 3D finite-element analysis of micropolar continuum.
\textit{Int J Solids Struct.}
2020;202:745--764.


\bibitem{Eremeyev2005}
Eremeyev  VA.
Nonlinear  micropolar  shells:  Theory  and  applications.
In: Pietraszkiewicz W,  Szymczak C. (Eds.), 
\textit{Shell structures: Theory and Applications}. 
Taylor \& Francis; 2005; 11--18.


\bibitem{Eremeyev2017}
V.A. Eremeyev, H. Altenbach, 
Basics of mechanics of micropolar shells.
In: Altenbach H, Eremeyev VA. (Eds.), 
\textit{Shell-like Structures}. 
Springer; 2017; 63--111.



\bibitem{Sargsyan2020}
Sargsyan A, 
Sargsyan S. 
Geometrically nonlinear models of static deformation of micropolar elastic thin plates and shallow shells.
\textit{Z Angew Math Mech.}
2020;e202000148.


\bibitem{Yoder2018}  
Yoder M, 
Thompson L, 
Summers J.
Size   effects   in   lattice   structures   and   a   comparison   to   micropolar elasticity.
\textit{Int J Solids Struct.}
2018;143:245--261.


\bibitem{Mayeur2011}  

Mayeur  JR, 
McDowell  DL, 
Bammann  DJ.
Dislocation-based micropolar single crystal plasticity:
Comparison of multi- and single criterion theorie.
\textit{J  Mech  Phys  Solids.} 
2011;59:398--422.

\bibitem{Zapata2020}  
Guar\'in-Zapata  N, 
Gomez  J, 
Valencia  C, 
Dargush  GF, 
Hadjesfandiari  AR.
Finite element modeling of micropolar-based phononic crystals.
\textit{Wave Motion}
2020;92:102406.


\bibitem{Spadoni2012}  
Spadoni A, 
Ruzzene M. 
Elasto-static micropolar behavior of a chiral auxetic lattice.
\textit{J  Mech  Phys  Solids.}
2012;60:156--171.


\bibitem{Goda2014} 
Goda I, 
Assidi  M, 
Ganghoffer  GF. 
A 3D elastic micropolar model of vertebral trabecular bone
from lattice homogenization of the bone microstructure.
\textit{Biomech Model Mechanobiol.}
2014;13:53--83.


\bibitem{Suh2020}  
Suh SS, 
Sun  W-C, 
O'Connor  DT. 
A phase field model for cohesive fracture in micropolar continua.
\textit{Comput Methods Appl Mech Eng.} 
2020;369:113181.

\bibitem{Sansour98}  
Sansour  C. 
Large  strain  deformations   of  elastic  shells, constitutive  modelling  and  finite  element  analysis.
\textit{Comput Methods Appl Mech Eng.}
1998;161:1--18.


\bibitem{SansKoll2000} 
Sansour  C, 
Kollmann  FG.  
Families of 4-node and 9-node finite elements for a finite deformation shell theory, an assessment of hybrid stress, hybrid strain and enhanced strain elements. 
\textit{Comput Mech.}
2000;24:435--447.






\bibitem{SimoArmero92}  
Simo  JC, 
Armero  F. 
Geometrically non-linear enhanced strain mixed methods and the method of incompatible modes. 
\textit{Int J Numer Methods Eng.} 
1992;33:1413--1449.


\bibitem{Simo93}
Simo  JC, 
Armero  F,
Taylor  RL.  
Improved versions of assumed enhanced strain tri-linear elements for 3D finite deformation problems. 
\textit{Comput Methods Appl Mech Eng.}
1993;110:359--386.

\bibitem{KW96} 
Korelc J, 
Wriggers P. 
Consistent gradient formulation for a stable enhanced strain method for large deformations,   
\textit{Eng Comput.} 
1996;13:103--123.

\bibitem{Glaser97}
Glaser  S, 
Armero  F. 
On the formulation of enhanced strain finite elements in finite deformations. 
\textit{Eng Comput.} 
1997;14:759--791.

\bibitem{I19}  
Itskov M.  
\textit{Tensor Algebra and Tensor Analysis for Engineers.}
Springer; 2019.




\bibitem{RamezBeam}  
Ramezani S,
Naghdabadi R, 
Sohrabpour S. 
Analysis of micropolar elastic beams. 
\textit{Eur J Mech A/Sol.} 
2009;28:202--208.

\bibitem{Wriggers2008}  
  Wriggers    P.  
\textit{Nonlinear Finite Element Methods}.
Springer; 2008.


\bibitem{Argyris1982} 
Argyris J. 
An excursion into large rotations.
\textit{Comput Methods Appl Mech Eng.} 
1982;32:85--155.

\bibitem{yan2015}
Yan B, 
Li  B, 
Kunecke  F, 
Gu  Z, 
Guo L.
Polypyrrole-based implantable electroactive pump for controlled drug microinjection.
\textit{ACS Appl Mater Interfaces.}
2018;7:14563--14568.

\bibitem{ju2021} 
Ju Y, 
Hu R, 
Xie Y, 
Yao J, 
Li X, 
Lv  Y, 
Han X, 
Cao Q, 
Li L.
Reconfigurable magnetic soft robots with multimodal locomotion.
\textit{Nano Energy}.
2021;87:106169.
\bibitem{carpenter2021}

Carpenter  JA, 
Eberle  TB, 
Schuerle  S, 
Rafsanjani  A, 
Studart  AR.
Facile manufacturing route for magneto-responsive soft actuators.
\textit{Adv Intell Syst.}
2021;3:20000283.
\bibitem{kadapa2022}
Kadapa C, 
Hossain M.
A unified numerical approach for soft to hard magneto-viscoelastically coupled polymers.
\textit{Mech Mat.} 
2022;166:104207.





%











































%

%







%
%
%
%
%
%



















































%
%
%
%


%
%




%
%
%
%
%
%
%








%






%
\end{thebibliography}
%
%

\end{document}